\documentclass[11pt, leqno]{article}
\usepackage{amssymb,latexsym}

\newcommand{\N}{{\mathbb N}}
\newcommand{\C}{{\mathbb C}}
\newcommand{\R}{{\mathbb R}}
\newcommand{\Q}{{\mathbb Q}}
\newcommand{\Z}{{\mathbb Z}}
\newcommand{\CP}{{\mathbb CP}}

\newtheorem{theorem}{Theorem}[section]

\newtheorem{corollary}[theorem]{Corollary}

\newtheorem{definition}[theorem]{Definition}
\newtheorem{example}[theorem]{Example}
\newtheorem{remark}[theorem]{Remark}

\newtheorem{lemma}[theorem]{Lemma}

\newtheorem{proposition}[theorem]{Proposition}
\def\Int{{ \mathrm{Int}}}

\begin{document}

\title{Gromov-Witten invariants \\ and pseudo symplectic capacities.}

\author{\\Guangcun Lu
\thanks{Partially supported by the NNSF 10371007 of China and the
Program for New Century Excellent Talents of the Education
Ministry of China.}\\
Department of Mathematics,  Beijing Normal University\\
Beijing 100875,   P.R.China\\
(gclu@bnu.edu.cn)}
\date{}
\maketitle
 \vspace{-0.1in}

\abstract{We introduce the concept of  pseudo symplectic
capacities which is a mild generalization of that of symplectic
capacities. As a generalization of the Hofer-Zehnder capacity we
construct a Hofer-Zehnder type pseudo symplectic capacity  and
estimate it in terms of Gromov-Witten invariants. The (pseudo)
symplectic capacities of Grassmannians and some product symplectic
manifolds are computed. As applications we first derive some
general nonsqueezing theorems that generalize and unite many
previous versions, then prove the Weinstein conjecture for
cotangent bundles over a large class of symplectic uniruled
manifolds (including the uniruled manifolds in algebraic geometry)
and also show that  any closed symplectic submanifold of
codimension two in any symplectic manifold has a small
neighborhood whose Hofer-Zehnder capacity is less than a given
positive number. Finally, we give two results on symplectic
packings in Grassmannians and on Seshadri constants.}
\vspace{-0.1in}

\tableofcontents

\section{Introduction and main results}
Gromov-Witten invariants and symplectic capacities are two kinds
of important symplectic invariants in symplectic geometry. Both
have many important applications. In particular, they are related
to the famous Weinstein conjecture and Hofer geometry (cf.\
\cite{En, FrGiSchl, FrSchl, HZ2, HV2, LaMc1,LaMc2, LiuT, Lu1, Lu2,
Lu3, Lu5, Lu7, Lu9, Mc2,Mc3, McSl, Po1,Po2, Po3, Schl, Schw, V1,
V2,V3, V4, We2} etc.). For some problems, Gromov-Witten invariants
are convenient and effective, but for other problems symplectic
capacities are more powerful. In the study of different problems
different symplectic capacities were defined. Examples of
symplectic capacities are the Gromov width ${\cal W}_G$
(\cite{Gr}), the Ekeland-Hofer capacity $c_{EH}$ (\cite{EH}), the
Hofer-Zehnder capacity $c_{HZ}$ (\cite{HZ1}) and Hofer's
displacement energy $e$ (\cite{H1}), the Floer-Hofer capacity
$c_{FH}$ (\cite{He}) and Viterbo's generating function capacity
$c_V$ (\cite{V3}). Only ${\cal W}_G$, $c_{HZ}$ and $e$ are defined
for all symplectic manifolds. In \cite{HZ1} an axiomatic
definition of a symplectic capacity was given. The Gromov width
${\cal W}_G$ is the smallest symplectic capacity. The
Hofer-Zehnder capacity is used in the study of many symplectic
topology questions. The reader can refer to \cite{HZ2,McSa1,V2}
for more details. But to the author's knowledge the relations
between Gromov-Witten invariants and symplectic capacities have
not been explored explicitly in the literature. Gromov-Witten
invariants are defined for closed symplectic manifolds
(\cite{FO,LiT,R,Sie}) and some non-closed symplectic manifolds
(cf.\ \cite{Lu4,Lu8}) and have been computed for many closed
symplectic manifolds. However, it is difficult to compute $c_{HZ}$
for a closed symplectic manifold. So far the only examples are
closed surfaces, for which $c_{HZ}$ is the area (\cite{Sib}),  and
complex projective space $(\CP^n,\sigma_n)$ with the standard
symplectic structure $\sigma_n$ related to the Fubini-Study
metric: Hofer and Viterbo proved $c_{HZ}(\CP^n,\sigma_n)=\pi$ in
\cite{HV2}. Perhaps the invariance of Gromov-Witten invariants
under deformations of the symplectic form is the main reason why
it is easier to compute them than Hofer-Zehnder capacities. Unlike
Gromov-Witten invariants, symplectic capacities do not depend on
homology classes of the symplectic manifolds in question. We
believe that this is a reason why they are difficult to compute or
estimate, and it is based on this observation that we introduced
the concept of pseudo symplectic capacities in the early version
\cite{Lu5} of this paper.

\subsection{Pseudo symplectic capacities}

In \cite{HZ1} a map $c$ from the class ${\cal C}(2n)$ of all
symplectic manifolds of dimension $2n$ to $[0, +\infty]$ is called
a symplectic capacity if
it satisfies the following properties: \\
({\bf monotonicity}) If there is a symplectic embedding
$(M_1,\omega_1)\to (M_2,\omega_2)$ of codimension zero then
$c(M_1,\omega_1)\le c(M_2,\omega_2)$;\\
  ({\bf conformality})
  $c(M,\lambda\omega)=|\lambda|c(M,\omega)$ for every
$\lambda\in\R\setminus\{0\}$;\\
({\bf nontriviality}) $c(B^{2n}(1),\omega_0)=\pi=c(Z^{2n}(1),\omega_0)$.\\
Here $B^{2n}(1)$ and $Z^{2n}(1)$ are the closed unit ball and
closed cylinder in the standard space $(\R^{2n},\omega_0)$, i.e.,
for any $r>0$,
$$
B^{2n}(r)=\{(x,y)\in\R^{2n}\;|\; |x|^2+ |y|^2\le r^2\}
$$
 and
 $$
Z^{2n}(r)=\{(x,y)\in\R^{2n}\;|\; x_1^2+ y_1^2\le r^2\}.
$$
Note
that the first property implies that $c$ is a symplectic
invariant.

Let $H_\ast(M; G)$ denote the singular homology of $M$ with
coefficient group $G$. For an integer $k \ge 1$ we denote by
${\cal C}(2n,k)$ the set of all tuples
$(M,\omega;\alpha_1,\cdots,\alpha_k)$ consisting of a
$2n$-dimensional connected symplectic manifold $(M,\omega)$ and
nonzero homology classes $\alpha_i\in H_\ast(M;G)$, $i=1,\cdots,
k$. We denote by $pt$ the homology class of a point.

\begin{definition}  \label{def:1.1}
{\rm A map $c^{(k)}$ from ${\cal C}(2n,k)$ to $[0, +\infty]$ is
called a $G_k$-{\bf pseudo symplectic capacity} if it satisfies
the following conditions.

\noindent{\rm P1.} {\bf Pseudo monotonicity:} If there is a
symplectic embedding $\psi:(M_1,\omega_1)\to (M_2,\omega_2)$ of
codimension zero, then for any $\alpha_i\in H_\ast(M_1;
G)\setminus\{0\}$, $i=1,\cdots,k$,
$$
c^{(k)}(M_1,\omega_1;\alpha_1,\cdots,\alpha_k)\le
  c^{(k)}(M_2,\omega_2;\psi_\ast(\alpha_1),\cdots,\psi_\ast(\alpha_k));
$$
\noindent{\rm P2.} {\bf Conformality:} $c^{(k)}(M,\lambda\omega;
\alpha_1,\cdots,
\alpha_k)=|\lambda|c^{(k)}(M,\omega;\alpha_1,\cdots,\alpha_k)$ for
every $\lambda\in\R\setminus\{0\}$ and all homology classes
$\alpha_i\in H_\ast(M; G) \setminus \{0\}$, $i=1,\cdots,k$;

\smallskip
\noindent{\rm P3.} {\bf Nontriviality:}
$c^{(k)}(B^{2n}(1),\omega_0; pt,\cdots, pt) =\pi$

\hspace{67.5mm}$=c^{(k)}(Z^{2n}(1),\omega_0; pt,\cdots, pt)$. }
\end{definition}

The pseudo monotonicity is the reason that a pseudo symplectic
capacity in general fails to be a symplectic invariant. If $k>1$
then a $G_{k-1}$-pseudo symplectic capacity $c^{(k-1)}$ is
naturally defined by
$$
c^{(k-1)}(M,\omega;
\alpha_1,\cdots,\alpha_{k-1}):=c^{(k)}(M,\omega; pt,
\alpha_1,\cdots,\alpha_{k-1}) ,
$$
and any $c^{(k)}$ induces a true symplectic capacity
$$
c^{(0)}(M,\omega) := c^{(k)}(M,\omega;pt,\cdots,pt).
$$

In this paper we shall concentrate on the case $k=2$ since in this
case there are interesting examples. More precisely, we shall
define a typical $G_2$-pseudo symplectic capacity of Hofer-Zehnder
type and give many applications. In view of our results we expect
that pseudo symplectic capacities will become a powerful tool in
the study of symplectic topology. Hereafter we assume $G = \Q$ and
often write $H_*(M)$ instead of $H_*(M;\Q)$.

\subsection{Construction of a pseudo symplectic capacity}

We begin with  recalling the Hofer-Zehnder capacity from
\cite{HZ1}. Given a symplectic manifold $(M,\omega)$, a smooth
function $H \colon M\to\R$ is called {\bf admissible}\, if there
exist an nonempty open subset $U$ and a compact subset
$K\subset M\setminus\partial M$ such that \\
(a) $H|_U=0$ and $H|_{M\setminus K}=\max H$;\\
(b) $0\le H\le\max H$;\\
(c) $\dot x=X_H(x)$ has no nonconstant fast periodic solutions. \\
Here $X_H$ is defined by $\omega(X_H, v)=dH(v)$ for $v \in TM$,
and ``fast'' means ``of period less than $1$''. Let ${\cal
H}_{ad}(M,\omega)$ be the set of admissible Hamiltonians on
$(M,\omega)$. The Hofer-Zehnder symplectic capacity
$c_{HZ}(M,\omega)$ of $(M,\omega)$ is defined by
$$
c_{HZ}(M,\omega) = \sup \left\{\max H\,|\, H\in {\cal
H}_{ad}(M,\omega) \right\}.
$$
Note that one can require the compact subset $K=K(H)$ to be a
proper subset of $M$ in the definition above. In fact, it suffices
to prove that for any $H\in{\cal H}_{ad}(M,\omega)$ and
$\epsilon>0$ small enough there exists a $H_\epsilon\in {\cal
H}_{ad}(M,\omega)$ such that $\max H_\epsilon \ge \max H-\epsilon$
and that the corresponding compact subset $K(H_\epsilon)$ is a
proper subset in $M$. Let us take a smooth function
$f_\epsilon:\R\to\R$ such that $0\le f'_\epsilon(t)\le 1$ and
$f_\epsilon(t)=0$ as $t\le 0$, and $f_\epsilon(t)=\max H-\epsilon$
as $t\ge\max H-\epsilon$. Then the composition $f_\epsilon\circ H$
is a desired $H_\epsilon$.

The invariant $c_{HZ}$ has many applications. Three of them are:
(i) giving a new proof of a foundational theorem in symplectic
topology -- Gromov's nonsqueezing theorem; (ii) studying the Hofer
geometry on the group of  Hamiltonian symplectomorphisms of a
symplectic manifold; (iii) establishing the existence of closed
characteristics on or near an energy surface. As mentioned above,
the difficulties in computing or estimating $c_{HZ}(M,\omega)$ for
a given symplectic manifold $(M,\omega)$ make it hard to find
further applications of this invariant. Therefore, it seems to be
important to give a variant of $c_{HZ}$ which can be easily
estimated and still has the above applications. An attempt was
made in \cite{McSl}.  In this paragraph we shall define a pseudo
symplectic capacity of Hofer-Zehnder type. The introduction of
such a pseudo symplectic capacity was motivated by  various papers
(e.g.\ \cite{LiuT,McSl}).

\begin{definition}  \label{def:1.2}
{\rm For a connected symplectic manifold $(M,\omega)$ of dimension
at least $4$ and two nonzero homology classes
$\alpha_0,\alpha_\infty\in H_\ast(M; \Q)$, we call a smooth
function $H:M\to\R$ {\bf $(\alpha_0,\alpha_\infty)$-admissible}
  (resp.\ {\bf $(\alpha_0,\alpha_\infty)^\circ$-admissible}) if
there exist two compact submanifolds $P$ and $Q$ of $M$ with
connected smooth boundaries and of codimension zero such that the
following condition groups $(1)(2)(3)(4)(5)(6)$ (resp.\
$(1)(2)(3)(4)(5)(6^\circ))$ hold:
\begin{description}
\item[(1)] $P\subset \Int (Q)$ and $Q\subset \Int (M)$;
\item[(2)]$H|_P=0$ and $H|_{M \setminus \Int (Q)}=\max H$;
\item[(3)] $0\le H\le \max H$; \item[(4)] There exist cycle
representatives of $\alpha_0$ and $\alpha_\infty$,
            still denoted by $\alpha_0,\alpha_\infty$, such that
            ${\rm supp}(\alpha_0)\subset \Int (P)$ and
            ${\rm supp}(\alpha_\infty)\subset M\setminus Q$;
\item[(5)] There are no critical values in $(0,\varepsilon)\cup
(\max H-\varepsilon, \max H)$
    for a small $\varepsilon=\varepsilon(H)>0$;

  \item[(6)] The
Hamiltonian system $\dot x=X_H(x)$ on $M$ has no nonconstant fast
periodic solutions; \item[(6$^\circ$)] The Hamiltonian system
$\dot x=X_H(x)$ on $M$ has no nonconstant contractible fast
periodic solutions.
\end{description}}
\end{definition}

We respectively denote by
\begin{equation} \label{e:1.1}
  {\cal H}_{ad}(M,\omega;\alpha_0,\alpha_\infty)\quad{\rm and}\quad
{\cal H}_{ad}^\circ(M,\omega;\alpha_0,\alpha_\infty)
\end{equation}
the set of all $(\alpha_0,\alpha_\infty)$-admissible and
$(\alpha_0,\alpha_\infty)^\circ$-admissible functions. Unlike
${\cal H}_{ad}(M,\omega)$ and ${\cal H}_{ad}^\circ(M,\omega)$,
  for some pairs $(\alpha_0,\alpha_\infty)$ the
sets in (\ref{e:1.1}) might be empty. On the other hand, one
easily shows that both sets in (\ref{e:1.1}) are nonempty if
$\alpha_0$ and $\alpha_\infty$ are separated by some  hypersurface
$S\subset M$ in the following sense.

\begin{definition}  \label{def:1.3}
{\rm A hypersurface $S\subset M$ is called {\bf separating the
homology classes} $\alpha_0,\alpha_\infty\in H_\ast(M)$ if (i) $S$
separates $M$ in the sense that there exist two submanifolds $M_0$
and $M_\infty$ of $M$ with common boundary $S$ such that $M_0\cup
M_\infty =M$ and $M_0\cap M_\infty =S$, (ii) there exist cycle
representatives of $\alpha_0$ and $\alpha_\infty$ with supports
contained in $\Int (M_0)$ and $\Int (M_\infty)$ respectively,
(iii) $M_0$ is compact and $\partial M_0=S$. }
\end{definition}

Without special statements a hypersurface in this paper always
means a smooth compact connected orientable submanifold of
codimension one and without boundary. Note that if $M$ is closed
and a hypersurface $S\subset M$ separates the homology classes
$\alpha_0$ and $\alpha_\infty$, then $S$ also separates
$\alpha_\infty$ and $\alpha_0$.

We define
\begin{equation} \label{e:1.2}
\left\{\begin{array}{ll}
C_{HZ}^{(2)}(M,\omega;\alpha_0,\alpha_\infty):=\sup\left\{\max
H\,|
\, H\in {\cal H}_{ad}(M,\omega;\alpha_0,\alpha_\infty)\right\},\\
C_{HZ}^{(2\circ)}(M,\omega;\alpha_0,\alpha_\infty):=\sup\left\{\max
H\,| \, H\in {\cal
H}_{ad}^\circ(M,\omega;\alpha_0,\alpha_\infty)\right\}.
\end{array}
\right.
\end{equation}
Hereafter we make the conventions that $\sup\emptyset=0$ and
$\inf\emptyset=+\infty$. As shown in Theorem~\ref{t:1.5} below,
$C_{HZ}^{(2)}$ is a $G_2$-pseudo symplectic capacity. We call it
{\bf pseudo symplectic capacity of Hofer-Zehnder type}.
$C_{HZ}^{(2)}$ and $C_{HZ}^{(2\circ)}$ in (\ref{e:1.2}) have
similar dynamical implications as the Hofer-Zehnder capacity
$c_{HZ}$. In fact, as in \cite{HZ2,HV2} one shows that
$$
0<C_{HZ}^{(2)}(M,\omega;\alpha_0,\alpha_\infty)<+\infty \quad
(0<C_{HZ}^{(2\circ)}(M,\omega;\alpha_0,\alpha_\infty)<+\infty)
$$
implies that every stable hypersurface $S\subset M$ separating
$\alpha_0$ and $\alpha_\infty$ carries a (contractible in $M$)
closed characteristic, i.e., there is an embedded (contractible in
$M$) circle in $S$ all of whose tangent lines belong to the
characteristic line bundle
\[
\mathcal{L}_S \,=\, \left\{ (x, \xi) \in TS \mid \omega (\xi,
\eta) =0 \mbox{ for all } \eta \in T_xS \right\} .
\]
This leads to the following version of the Weinstein conjecture.
\vspace{2mm}

 {\bf $(\alpha_0, \alpha_\infty)$-Weinstein conjecture}: Every
hypersurface $S$ of contact type in a symplectic manifold
$(M,\omega)$ separating $\alpha_0$ and $\alpha_\infty$ carries a
closed characteristic.

\vspace{2mm}

\noindent{In terms of this language the main result Theorem 1.1 in
\cite{LiuT} asserts that the $(\alpha_0,\alpha_\infty)$-Weinstein
conjecture holds if some GW-invariant
$$
\Psi_{A, g,
m+2}(C;\alpha_0,\alpha_\infty,\beta_1,\cdots,\beta_m)
$$
does not
vanish, see 1.3 below.}

As before let $pt$ denote the generator of $H_0(M;\Q)$ represented
by a point. Then we have the true symplectic capacities
\begin{equation}\label{e:1.3}
\left\{\begin{array}{ll}
C_{HZ}(M,\omega):=C_{HZ}^{(2)}(M,\omega; pt, pt),\\
C_{HZ}^{\circ}(M,\omega):=C_{HZ}^{(2\circ)}(M,\omega;pt, pt).
\end{array}
\right.
\end{equation}
Recall that we have also the $\pi_1$-sensitive Hofer-Zehnder
capacity denoted by $\bar C_{HZ}$ in \cite{Lu1} and by
$c^\circ_{HZ}$ in \cite{Schw}. By definitions, it is obvious that
$C_{HZ}(M, \omega)\le c_{HZ}(M, \omega)$ and $C_{HZ}^\circ(M,
\omega)\le c^\circ_{HZ}(M,\omega)$ for any symplectic manifold
$(M,\omega)$. One naturally asks when $C_{HZ}$ (resp.\
$C_{HZ}^\circ$) is equal to $c_{HZ}$ (resp.\ $c^\circ_{HZ}$). The
following result partially answers this question.

\begin{lemma}\label{lem:1.4} Let a symplectic manifold $(M,\omega)$
satisfy one of the following conditions:
\begin{description}
 \item[(i)] $(M,\omega)$ is closed.
\item[(ii)] For each compact subset $K\subset M\setminus\partial
M$ there exists a compact submanifold $W\subset M$ with connected
boundary and of codimension zero such that $K\subset W$. Then
$$C_{HZ}(M, \omega)=c_{HZ}(M, \omega)\quad{\rm and}\quad C_{HZ}^\circ(M,
\omega)=c^\circ_{HZ}(M,\omega).$$
\end{description}
\end{lemma}

For arbitrary homology classes $\alpha_0,\alpha_\infty\in
H_\ast(M)$,
\begin{equation}\label{e:1.4}
\left\{\begin{array}{ll}
C_{HZ}^{(2)}(M,\omega;\alpha_0,\alpha_\infty)\le
  C_{HZ}^{(2\circ)}(M,\omega;\alpha_0,\alpha_\infty),\\
C_{HZ}^{(2)}(M,\omega;\alpha_0,\alpha_\infty)\le C_{HZ}(M,\omega),\\
  C_{HZ}^{(2\circ)}(M,\omega;\alpha_0,\alpha_\infty)\le
C^\circ_{HZ}(M,\omega).
\end{array}
\right.
\end{equation}

Both $C_{HZ}^{(2)}$ and $C_{HZ}^{(2\circ)}$ are important because
estimating or calculating them is easier than for $C_{HZ}$ and
$C_{HZ}^\circ$, and because they still share those properties
needed for applications. In Remark~\ref{rem:1.28} we will give an
example which illustrates that  sometimes $C_{HZ}^{(2)}$ gives
better results than $C_{HZ}$. Recall that the Gromov width ${\cal
W}_G$ is the smallest symplectic capacity so that
\begin{equation}\label{e:1.5}
{\cal W}_G\le C_{HZ}\le C^\circ_{HZ} .
  \end{equation}
\vspace{2mm}

\noindent{\bf Convention}:\hspace{2mm} $C$ stands for both
$C^{(2)}_{HZ}$ and ${C}^{(2\circ)}_{HZ}$ if there is no danger of
confusion. \vspace{2mm}

  The following theorem shows that
$C^{(2)}_{HZ}$ is indeed a pseudo symplectic capacity.

\begin{theorem}\label{t:1.5}$\!{\rm :}$
\begin{description}
\item[(i)] If $M$ is closed then for any nonzero homology classes
$\alpha_0,\alpha_\infty\in H_\ast(M; \Q)$,
$$
C(M,\omega;\alpha_0,
\alpha_\infty)=C(M,\omega;\alpha_\infty,\alpha_0).
$$
\item[(ii)] $C(M,\omega;\alpha_0,\alpha_\infty)$ is invariant
under those symplectomorphisms \linebreak $\psi\in{\rm
Symp}(M,\omega)$ which induce the identity on $H_\ast(M;\Q)$.

\item[(iii)]{\bf (Normality)} For any $r>0$ and nonzero
$\alpha_0,\alpha_\infty\in H_\ast(B^{2n}(r);\Q)$ or
$H_\ast(Z^{2n}(r);\Q)$,
$$
C(B^{2n}(r),\omega_0;\alpha_0,\alpha_\infty)=C(Z^{2n}(r),\omega_0;\alpha_0,\alpha_\infty)=\pi
r^2.
$$
\item[(iv)]{\bf (Conformality)}  For any nonzero real number
$\lambda$,
$$
C(M, \lambda\omega;\alpha_0,\alpha_\infty)= |\lambda|
C(M,\omega;\alpha_0,\alpha_\infty).
$$
\item[(v)]{\bf (Pseudo monotonicity)} For any symplectic embedding
$\psi: (M_1, \omega_1)\to (M_2,\omega_2)$ of codimension zero and
any nonzero $\alpha_0,\alpha_\infty\in H_\ast(M_1; \Q)$,
$$
C_{HZ}^{(2)}(M_1, \omega_1;\alpha_0,\alpha_\infty)\le
  C_{HZ}^{(2)}(M_2,\omega_2; \psi_\ast(\alpha_0),
\psi_\ast(\alpha_\infty)).
$$
Furthermore, if $\psi$ induces an injective homomorphism
$\pi_1(M_1)\to\pi_1(M_2)$ then
$$
C_{HZ}^{(2\circ)}(M_1, \omega_1;\alpha_0,\alpha_\infty)\le
C_{HZ}^{(2\circ)}(M_2,\omega_2; \psi_\ast(\alpha_0),
\psi_\ast(\alpha_\infty)).
$$
\item[(vi)] For any $m\in\N$
$$\left\{\begin{array}{ll}
&C(M,\omega;\alpha_0,\alpha_\infty)\le C(M,\omega;
m\alpha_0,\alpha_\infty),\\
&C(M,\omega;\alpha_0,\alpha_\infty)\le C(M,\omega; \alpha_0,
m\alpha_\infty),\\
&C(M,\omega;-\alpha_0,\alpha_\infty)= C(M,\omega; \alpha_0,
\alpha_\infty)=C(M,\omega;\alpha_0, -\alpha_\infty).
\end{array}\right.$$
\item[(vii)] If $\dim\alpha_0+\dim\alpha_\infty\le\dim M-2$ and
  $\alpha_0$ or  $\alpha_\infty$ can be represented by a  connected
closed submanifold, then
\[
C(M,\omega;\alpha_0,\alpha_\infty)>0.
\]
\end{description}
\end{theorem}

\begin{remark}\label{rem:1.6}
{\rm If $M$ is not closed, $C(M,\omega; pt, \alpha)$ and
$C(M,\omega; \alpha, pt)$ might be different. For example, let $M$
be the annulus in $\R^2$ of area $2$, and $\alpha$ be a generator
of $H_1(M)$. Then ${\cal W}_G(M, \omega) = C_{HZ}^{(2)} (M,
\omega; pt, \alpha) =2$, while $C_{HZ}^{(2)} (M, \omega; \alpha,
pt) =0$ since ${\cal H}_{ad}(M,\omega;\alpha, pt)=\emptyset$. This
example also shows that the dimension assumption
$\dim\alpha_0+\dim\alpha_\infty\le\dim M-2$ cannot be weakened.
}\end{remark}

\begin{proposition}  \label{prop:1.7}
  Let $W\subset \Int (M)$ be a smooth compact submanifold of codimension
  zero and with connected boundary such that
the homology classes $\alpha_0,\alpha_\infty\in
H_\ast(M;\Q)\setminus\{0\}$ have representatives supported in
$\Int (W)$ and $\Int (M)\setminus W$, respectively. Denote by
$\tilde\alpha_0\in H_\ast(W; \Q)$ and $\tilde\alpha_\infty\in
H_\ast(M\setminus W; \Q)$ the nonzero homology classes determined
by them. Then
\begin{equation}  \label{e:1.6}
C_{HZ}^{(2)}(W,\omega;\tilde\alpha_0, pt) \le
C_{HZ}^{(2)}(M,\omega;\alpha_0,\alpha_\infty),
\end{equation}
and we specially have
\begin{equation}  \label{e:1.7}
c_{HZ}(W,\omega)=C_{HZ}(W,\omega)\le C_{HZ}^{(2)}(M,\omega; pt,
\alpha)
\end{equation}
for any $\alpha\in H_\ast(M; \Q)\setminus\{0\}$ with
representative supported in $\Int (M) \setminus W$. If the
inclusion $W\hookrightarrow M$ induces an injective homomorphism
$\pi_1(W)\to\pi_1(M)$ then
\begin{equation}  \label{e:1.8}
C_{HZ}^{(2\circ)}(W,\omega;\tilde\alpha_0, pt) \le
C_{HZ}^{(2\circ)}(M,\omega;\alpha_0,\alpha_\infty)
\end{equation}
and corresponding to (\ref{e:1.7}) we have
\begin{equation}  \label{e:1.9}
c_{HZ}^\circ(W,\omega)=C_{HZ}^\circ(W,\omega)\le
C_{HZ}^{(2\circ)}(M,\omega; pt, \alpha).
\end{equation}
   Also
\begin{equation}  \label{e:1.10}
C_{HZ}^{(2)}(M\setminus W,\omega;\tilde\alpha_\infty, pt)\le
C_{HZ}^{(2)}(M,\omega;\alpha_\infty, \alpha_0),
\end{equation}
and
\begin{equation}  \label{e:1.11}
C_{HZ}^{(2\circ)}(M\setminus W,\omega;\tilde\alpha_\infty, pt)\le
C_{HZ}^{(2\circ)}(M,\omega;\alpha_\infty, \alpha_0)
\end{equation}
if the inclusion $M\setminus W\hookrightarrow M$ induces an
injective homomorphism $\pi_1(M\setminus W)\to\pi_1(M)$.
Furthermore, for any $\alpha\in H_\ast(M; \Q)\setminus\{0\}$ with
$\dim\alpha\le\dim M-1$,
\begin{equation}  \label{e:1.12}
{\cal W}_G(M,\omega) \le C(M,\omega; pt, \alpha).
\end{equation}
\end{proposition}

For closed symplectic manifolds,  Proposition~\ref{prop:1.7} can
be strengthened as follows.
\begin{theorem}  \label{t:1.8}
If in the situation of Proposition~\ref{prop:1.7} the symplectic
manifold $(M,\omega)$ is closed and $M \setminus W$ is connected,
then
\begin{equation}\label{e:1.13}
\qquad C_{HZ}^{(2)}(W,\omega;\tilde\alpha_0, pt)+
C_{HZ}^{(2)}(M\setminus W,\omega;\tilde\alpha_\infty, pt) \le
C_{HZ}^{(2)}(M,\omega;\alpha_0,\alpha_\infty).
\end{equation}
In particular, if $\alpha\in H_\ast(M;\Q)\setminus\{0\}$ has a
representative supported in $M\setminus W$ and thus determines a
homology class $\tilde\alpha\in H_\ast(M\setminus
W;\Q)\setminus\{0\}$, then
\[
c_{HZ}(W,\omega)+ C_{HZ}^{(2)}(M\setminus W,\omega;\tilde\alpha,
pt) \le C_{HZ}^{(2)}(M,\omega; pt, \alpha).
\]
If both inclusions $W\hookrightarrow M$ and $M\setminus
W\hookrightarrow M$ induce an injective homomorphisms
$\pi_1(W)\to\pi_1(M)$ and $\pi_1(M\setminus W)\to\pi_1(M)$, then
\begin{equation}\label{e:1.14}
\qquad C_{HZ}^{(2\circ)}(W,\omega;\tilde\alpha_0, pt)+
C_{HZ}^{(2\circ)}(M\setminus W,\omega;\tilde\alpha_\infty, pt) \le
C_{HZ}^{(2\circ)}(M,\omega;\alpha_0,\alpha_\infty),
\end{equation}
and specially
\[
c_{HZ}^\circ(W,\omega)+ C_{HZ}^{(2\circ)}(M\setminus
W,\omega;\tilde\alpha, pt) \le C_{HZ}^{(2\circ)}(M,\omega; pt,
\alpha)
\]
for any $\alpha\in H_\ast(M;\Q)\setminus\{0\}$ with a
representative supported in $M\setminus W$.
\end{theorem}

An inequality similar to (\ref{e:1.13}) was first proved for the
usual Hofer-Zehnder capacity by Mei-Yu Jiang \cite{Ji}. In the
following subsections we always take $G=\Q$.

\subsection{Estimating the pseudo capacity in terms of
Gromov-Witten invariants}

To state our main results we recall that for a given class $A\in
H_2(M;\Z)$ the Gromov-Witten invariant of genus $g$ and with $m+2$
marked points is a homomorphism
$$\Psi_{A, g, m+2}: H_\ast(\overline{\cal M}_{g, m+2};\Q)\times
H_\ast(M;\Q)^{m+2} \to \Q .
$$
We refer to the appendix and  \cite{FO,LiT,R,Sie} and \cite{Lu8}
for more details on Gromov-Witten invariants.

The Gromov-Witten invariants for general (closed) symplectic
manifolds were constructed by different methods, cf. \cite{FO,
LiT, R, Sie}, and \cite{LiuT} for a Morse theoretic set-up. It is
believed that these methods define the same symplectic
Gromov-Witten invariants, but no proof has been written down so
far. A detailed construction of the GW-invariants by the method in
\cite{LiuT}, including proofs of the composition law and reduction
formula, was given in \cite{Lu8} for a larger class of symplectic
manifolds including all closed symplectic manifolds. The method by
Liu-Tian was also used in \cite{Mc2}. Without special statements,
the Gromov-Witten invariants in this paper are the ones
constructed by the method in \cite{LiuT}. The author strongly
believes that they agree with those constructed in \cite{R}.

\begin{definition}  \label{def:1.9}
{\rm Let $(M,\omega)$ be a closed symplectic manifold and let
$\alpha_0,\alpha_\infty\in H_\ast(M;\Q)$. We define
$$
{\rm GW}_g(M,\omega;\alpha_0,\alpha_\infty)\in (0, +\infty]
$$
as the infimum of the $\omega$-areas $\omega(A)$ of the homology
classes $A\in H_2(M;\Z)$ for which the Gromov-Witten invariant
$\Psi_{A, g,
m+2}(C;\alpha_0,\alpha_\infty,\beta_1,\cdots,\beta_m)\ne 0$ for
some homology classes $\beta_1,\cdots,\beta_m\in H_\ast(M;\Q)$ and
$C\in  H_\ast(\overline{\cal M}_{g, m+2};\Q)$ and an integer $m
\ge 1$. We define
$$
{\rm GW}(M,\omega;\alpha_0,\alpha_\infty):= \inf \left\{ {\rm
GW}_g(M,\omega;\alpha_0,\alpha_\infty)\,|\, g\ge 0 \right\}\in [0,
+\infty].
$$
}
\end{definition}

The positivity ${\rm GW}_g(M,\omega;\alpha_0,\alpha_\infty)>0$
follows from the compactness of the space of $J$-holomorphic
stable maps (cf.\ \cite{FO,LiT,R,Sie}). Here we have used  the
convention  $\inf\emptyset=+\infty$ below (\ref{e:1.2}). One
easily checks that both ${\rm GW}_g$ and ${\rm GW}$ satisfy the
pseudo monotonicity and conformality in Definition~\ref{def:1.1}.
As Professor Dusa McDuff suggested, one can consider closed
symplectic manifolds only and replace the nontriviality condition
in Definition~\ref{def:1.1} by
$$
c^{(2)}(\CP^n,\sigma_n; pt, pt)=c^{(2)}(\CP^1\times
T^{2n-2},\sigma_1\oplus\omega_0; pt, [pt\times T^{2n-2}])=\pi;
$$
then both ${\rm GW}_0$ and ${\rm GW}$ are pseudo symplectic
capacities in view of (\ref{e:1.19}) and (\ref{e:1.22}) below. The
following result is the core of this paper. Its proof is given in
$\S3$  based on \cite{LiuT} and the key Lemma~\ref{lem:3.3}.

\begin{theorem}  \label{t:1.10}
For any closed symplectic manifold $(M,\omega)$ of dimension $\dim
M \ge 4$ and homology classes $\alpha_0,\alpha_\infty\in
H_\ast(M;\Q)\setminus\{0\}$ we have
\begin{equation}  \label{e:1.15}
C_{HZ}^{(2)}(M,\omega;\alpha_0,\alpha_\infty)\le {\rm
GW}(M,\omega;\alpha_0,\alpha_\infty)
\end{equation}
and
\begin{equation}  \label{e:1.16}
C_{HZ}^{(2\circ)}(M,\omega;\alpha_0,\alpha_\infty)\le {\rm
GW}_0(M,\omega;\alpha_0,\alpha_\infty).
\end{equation}
\end{theorem}

\begin{remark}  \label{rem:1.11}
{\rm By the reduction formula (\ref{e:7.3}) for Gromov-Witten
invariants recalled in the appendix,
\begin{eqnarray*}
\Psi_{A, g, m+3}([\pi_{m+3}^{-1}(K)];\alpha_0,\alpha_\infty,
\alpha,\beta_1,\cdots, \beta_m)\\
=PD(\alpha)(A)\cdot\Psi_{A, g,
m+2}([K];\alpha_0,\alpha_\infty, \beta_1,\cdots,\beta_m)
\end{eqnarray*}
for any $\alpha\in H_{2n-2}(M,\Z)$ and $[K]\in
H_\ast(\overline{\mathcal M}_{g, m+2},\Q)$. Here $2n=\dim M$. It
easily follows that ${\rm
GW}_g(M,\omega;\alpha_0,\alpha_\infty)<+\infty$ implies that
  ${\rm GW}_g(M,\omega;\alpha_0,\alpha)$,
${\rm GW}_g(M,\omega;\alpha,\alpha_\infty)$ and ${\rm
GW}_g(M,\omega;\alpha,\beta)$ are finite for any $\alpha,\beta\in
H_{2n-2}(M,\Z)$ with $PD(\alpha)(A)\ne 0$ and $PD(\beta)(A)\ne 0$.
In particular, it is easily proved that for any integer $g\ge 0$
\begin{equation}\label{e:1.17}
\qquad {\rm GW}_g(M,\omega; pt, PD([\omega]))=\inf\{ {\rm
GW}_g(M,\omega; pt,\alpha)\,|\,\alpha\in H_\ast(M,\Q)\}
\end{equation} }
\end{remark}

\begin{corollary}  \label{cor:1.12}
If ${\rm GW}_g(M,\omega;\alpha_0,\alpha_\infty)<+\infty$ for some
integer $g\ge 0$ then the $(\alpha_0,\alpha_\infty)$-Weinstein
conjecture holds in $(M,\omega)$.
\end{corollary}

Many results in this paper are based on the following special case
of Theorem~\ref{t:1.10}.

\begin{theorem}  \label{t:1.13}
For any closed symplectic manifold $(M,\omega)$ of dimension at
least four and a nonzero homology class $\alpha\in H_\ast(M;\Q)$,
$$
C_{HZ}^{(2)}(M,\omega; pt, \alpha)\le {\rm GW}(M,\omega;
pt,\alpha)
$$
and
$$
C_{HZ}^{(2\circ)}(M,\omega; pt, \alpha)\le {\rm GW}_0(M,\omega;
pt,\alpha).
$$
\end{theorem}

\begin{definition}  \label{def:1.14}
{\rm Given a nonnegative integer $g$, a closed symplectic manifold
$(M,\omega)$ is called $g$-{\bf symplectic uniruled} if $\Psi_{A, g,
m+2}(C; pt,\alpha,\beta_1,\cdots,\beta_m)\ne 0$ for some homology
classes $A\in H_2(M;\Z)$, $\alpha, \beta_1,\cdots, \beta_m\in
H_\ast(M;\Q)$ and $C\in H_\ast(\overline{\cal M}_{g, m+2};\Q)$ and
an integer $m \ge 1$. If $C$ can be chosen as a point $pt$ we say
$(M,\omega)$ to be {\bf strong} $g$-{\bf symplectic uniruled}.
Moreover, $(M,\omega)$ is called {\bf symplectic uniruled} $($resp.\
{\bf strong symplectic uniruled}$)$ if it is $g$-symplectic uniruled
$($resp.\ strong $g$-symplectic uniruled$)$ for some  integer $g\ge
0$. }
\end{definition}

It was proved in (\cite{Ko}) and (\cite{R}) that (projective
algebraic) uniruled manifolds are strong $0$-symplectic uniruled
\footnote{This is the only place in which we assume that our
GW-invariants agree with the ones in \cite{R}. In a future paper we
shall use the method in \cite{LiuT} and the techniques in \cite{Lu8}
to prove this fact.}. In Proposition~\ref{prop:7.3} we shall prove
that for a closed symplectic manifold $(M,\omega)$, if there exist
homology classes $A\in H_2(M;\Z)$ and $\alpha_i\in H_\ast(M;\Q)$,
$i=1,\cdots, k$, such that the Gromov-Witten invariant $\Psi_{A, g,
k+1}(pt; pt,\alpha_1,\cdots, \alpha_k)\ne 0$ for some integer $g\ge
0$, then there exists a homology class $B\in H_2(M;\Z)$ with
$\omega(B)\le\omega(A)$ and $\beta_i\in H_\ast(M;\Q)$, $i=1,2$, such
that the Gromov-Witten invariant $\Psi_{B, 0, 3}(pt;
pt,\beta_1,\beta_2)\ne 0$. Therefore, every strong symplectic
uniruled manifold is strong $0$-symplectic uniruled. Actually, we
shall prove in Proposition~\ref{prop:7.5} that the product of any
closed symplectic manifold and a strong symplectic uniruled
manifolds is strong symplectic uniruled. Moreover, the class of
$g$-symplectic uniruled manifolds is closed under deformations of
symplectic forms because Gromov-Witten invariants are symplectic
deformation invariants. For a $g$-symplectic uniruled manifold
$(M,\omega)$, i.e., ${\rm GW}_g(M,\omega; pt, PD([\omega]))
<+\infty$, the author observed in \cite{Lu3} that if a hypersurface
of contact type $S$ in $(M,\omega)$ separates $M$ into two parts
$M_+$ and $M_-$, then there exist two classes $PD([\omega])_+$ and
$PD([\omega])_-$ in $H_{2n-2}(M,\R)$ with cycle representatives
supported in $M_+$ and $M_-$ respectively such that $PD([\omega])_+
+ PD([\omega])_-=PD([\omega])$ and that at least one of the numbers
${\rm GW}_g(M,\omega; pt, PD([\omega])_+)$ or ${\rm GW}_g(M,\omega;
pt, PD([\omega])_-)$ is finite. Theorem~\ref{t:1.13} (or
(\ref{e:1.15})) implies that at least one of the following two
statements holds:
\begin{equation}\label{e:1.18}
\left.\begin{array}{ll}
  & C_{HZ}^{(2)}(M,\omega; pt,
PD([\omega])_+)\le {\rm GW}_g(M,\omega; pt,
PD([\omega])_+)<+\infty\quad{\rm or}\\
&C_{HZ}^{(2)}(M,\omega; pt, PD([\omega])_-)\le {\rm
GW}_g(M,\omega; pt, PD([\omega])_-)<+\infty.
\end{array}\right.
\end{equation}
On the other hand (\ref{e:1.12}) shows that
$C_{HZ}^{(2)}(M,\omega; pt, PD([\omega])_+)$ and\linebreak
$C_{HZ}^{(2)}(M,\omega; pt, PD([\omega])_-)$ are always positive.
Consequently, $S$ carries a nontrivial closed characteristic,
i.e., the $(pt, pt)$-Weinstein conjecture holds in symplectic
uniruled manifolds (\cite{Lu3}).

The Grassmannians and their products with any closed symplectic
manifold are symplectic uniruled. For them we have

\begin{theorem}  \label{t:1.15}
For the Grassmannian $G(k, n)$ of $k$-planes in $\C^n$ we denote
by $\sigma^{(k,n)}$ the canonical symplectic form for which
$\sigma^{(k,n)}(L^{(k,n)})=\pi$ for the generator $L^{(k,n)}$ of
$H_2(G(k,n);\Z)$. Let the submanifolds $X^{(k,n)}\approx G(k,
n-1)$ and $Y^{(k,n)}$ of $G(k,n)$ be given by $\{ V\in G(k,n) \mid
w_0^\ast v=0 \mbox{ for all } v \in V \}$ and $\{V\in G(k,n) \mid
v_0\in V\}$ for some fixed $v_0, w_0\in\C^n \setminus \{0\}$
respectively. Their homology classes $[X^{(k,n)}]$ and
$[Y^{(k,n)}]$ are independent of the choices of $v_0, w_0\in\C^n
\setminus \{0\}$ and $\deg [X^{(k,n)}]=2k(n-k-1)$ and $\deg
[Y^{(k,n)}]=2(k-1)(n-k)$. Then
\[
{\cal W}_G(G(k, n),\sigma^{(k,n)})= C_{HZ}^{(2)}(G(k,n),
\sigma^{(k,n)}; pt, \alpha)=\pi
\]
for $\alpha=[X^{(k,n)}]$ or $\alpha=[Y^{(k,n)}]$ with $k \le n-2$.
\end{theorem}

In particular, if $k=1$ and $n\ge 3$ then $[Y^{(1,n)}]=pt$ and
$(G(1, n), \sigma^{(1,n)})=(\CP^{n-1}, \sigma_{n-1})$, where
$\sigma_{n-1}$ the unique ${\rm U}(n)$-invariant K\"ahler form on
$\CP^{n-1}$ whose integral over the line $\CP^1\subset\CP^{n-1}$
is equal to $\pi$. In this case Theorem~\ref{t:1.15} and
Lemma~\ref{lem:1.4} yield:
\begin{equation}\label{e:1.19}
\left.\begin{array}{ll}
 c_{HZ}(\CP^{n-1},
\sigma_{n-1})=C_{HZ}(\CP^{n-1},
\sigma_{n-1})\\
\hspace{30mm}:=C_{HZ}^{(2)}(\CP^{n-1}, \sigma_{n-1}; pt, pt)=\pi.
\end{array}\right.
\end{equation}
Hofer and Viterbo \cite{HV2} firstly proved that $c_{HZ}(\CP^n,
\sigma_n)=\pi$. Therefore, Theorem~\ref{t:1.15} can be viewed as a
generalization of their result. If $k=1$, on one hand the volume
estimate gives ${\cal W}_G \left( \CP^{n-1}, \sigma_{n-1} \right)
\le \pi$, and on the other hand there exists an explicit
symplectic embedding $B^{2n-2}(1) \hookrightarrow \left(
\CP^{n-1}, \sigma_{n-1} \right)$, see \cite{Ka, HV2}. So we have
${\cal W}_G \left( \CP^{n-1}, \sigma_{n-1} \right) = \pi$. For $k
\ge 2$, however, the remarks below Theorem~\ref{t:1.35} show that
the identity ${\cal W}_G(G(k, n),\sigma^{(k,n)})=\pi$ does not
follow so easily. Karshon and Tolman [KaTo1] independently
computed ${\cal W}_G(G(k, n),\sigma^{(k,n)})$ in a different
  method.

\begin{theorem}  \label{t:1.16}
For any closed symplectic manifold
  $(M,\omega)$,
\begin{equation}  \label{e:1.20}
C(M\times G(k,n), \omega\oplus (a\sigma^{(k,n)}); pt,
[M]\times\alpha) \le |a|\pi
\end{equation}
for any $a\in\R\setminus\{0\}$ and $\alpha=[X^{(k,n)}]$ or
$\alpha=[Y^{(k,n)}]$ with $k\le n-2$. Moreover, for the product
$$(W,\Omega)=\bigl(G(k_1, n_1)\times\cdots\times G(k_r,n_r),
(a_1\sigma^{(k_1,n_1)})\oplus\cdots\oplus
(a_r\sigma^{(k_r,n_r)})\bigr)$$
   we have
\begin{equation}  \label{e:1.21}
C(W, \Omega; pt, \alpha_1\times\cdots\times\alpha_r) \le
(|a_1|+\cdots +|a_r|)\pi \end{equation}
  for any
$a_i\in\R\setminus\{0\}$ and $\alpha_i=[X^{(k_i,n_i)}]$ or
$[Y^{(k_i, n_i)}]$. Furthermore,
\begin{equation}  \label{e:1.21.5}
{\cal W}_G(G(k_1, n_1)\times\cdots\times G(k_r,n_r),
\sigma^{(k_1,n_1)}\oplus\cdots\oplus\sigma^{(k_r,n_r)})=\pi.
\end{equation}
\end{theorem}

For the projective space $\CP^n=G(1, n+1)$ we have:\vspace{2mm}

\begin{theorem}  \label{t:1.17}
Let $(M,\omega)$ be a closed symplectic manifold and $\sigma_n$
the unique ${\rm U}(n+1)$-invariant K\"ahler form on $\CP^n$ whose
integral over the line $\CP^1\subset\CP^n$ is equal to $\pi$. Then
\begin{equation}  \label{e:1.22}
C(M\times\CP^n,\omega\oplus(a\sigma_n); pt, [M\times pt])= |a|\pi
\end{equation}
  for any
$a\in\R\setminus\{0\}$. Moreover, for any $r>0$ and the standard
ball $B^{2n}(r)$ of radius $r$ and the cylinder
$Z^{2n}(r)=B^2(r)\times\R^{2n-2}$ in $(\R^{2n},\omega_0)$, we have
\begin{equation}  \label{e:1.23}
C(M\times B^{2n}(r),\omega\oplus\omega_0)= C(M\times
Z^{2n}(r),\omega\oplus\omega_0)=\pi r^2
\end{equation}
for $C=C_{HZ}$, $C^\circ_{HZ}$, $c_{HZ}$ and $c^\circ_{HZ}$.
\end{theorem}

\begin{remark}  \label{rem:1.18}
{\rm Combining the arguments in \cite{McSl,Lu1} one can prove a
weaker version of (\ref{e:1.23})  for any weakly monotone
noncompact geometrically bounded symplectic manifold $(M,\omega)$
and any $r>0$, namely
$$C_{HZ}^{\circ}(M\times B^{2n}(r),\omega\oplus\omega_0)\le
  C_{HZ}^{\circ}(M\times Z^{2n}(r),\omega\oplus\omega_0)\le\pi r^2.
$$
This generalization can be used to find periodic orbits of a
charge subject to a magnetic field (cf.\ \cite{Lu2}). }
\end{remark}

 From Theorem~\ref{t:1.13} and Lemma~\ref{lem:1.4} we obtain

\begin{corollary}  \label{cor:1.19}
For any closed symplectic manifold $(M,\omega)$ of dimension at
least $4$ we have
$$
c_{HZ}(M,\omega)\le {\rm GW}(M,\omega; pt, pt),\quad
c_{HZ}^{\circ}(M,\omega)\le{\rm GW}_0(M,\omega; pt, pt).$$
\end{corollary}

Thus $c_{HZ}(M,\omega)$ is finite if the Gromov-Witten invariant
$$
\Psi_{A, g, m+2}(C; pt, pt,\beta_1,\cdots,\beta_m)
$$
does not vanish for some homology classes $A\in H_2(M;\Z)$,
$\beta_1,\dots,\beta_m\in H_\ast(M;\Q)$ and $C\in
H_\ast(\overline{\cal M}_{g, m+2};\Q)$ and integers $g\ge 0$ and
$m>0$. Notice that ${\rm GW}_0(M,\omega; pt, pt)$ is needed here.
For example, consider
$$
(M,\omega)=(\CP^1\times\CP^1, \sigma_1\oplus\sigma_1).
$$
The following Theorem~\ref{t:1.21} and its proof show that
$c_{HZ}(M,\omega)=c_{HZ}^{\circ}(M,\omega)=2\pi$ and ${\rm
GW}_0(M,\omega; pt, pt)=2\pi$. However, one easily proves that
\begin{eqnarray*}
{\rm GW}_0(M,\omega; pt, PD([\omega]))\!\!\!\!\!\!\!\!&&={\rm
GW}_0(M,\omega; pt,
[pt\times\CP^1])\\
 &&={\rm GW}_0(M,\omega; pt, [\CP^1\times
pt])=\pi.
\end{eqnarray*}
So ${\rm GW}_0(M,\omega; pt, pt)$ is necessary.

\begin{example}  \label{ex:1.20}
{\rm (i) For a smooth complete intersection $(X,\omega)$ of
degree\linebreak $(d_1,\cdots, d_k)$ in $\CP^{n+k}$ with
$n=2\sum(d_i-1)-1$ or $3\sum(d_i-1)-3$, we have
$c_{HZ}^{\circ}(X,\omega)=C_{HZ}^{\circ}(X,\omega)<+\infty$.

\noindent{(ii)} For a rational algebraic manifold $(X,\omega)$, if
there exists a surjective morphism $\pi: X\to\CP^n$ such that
$\pi|_{X\setminus S}$ is one to one for some subvariety $S$ of $X$
with ${\rm codim}_{\C}\pi(S)\ge 2$ then
$c_{HZ}^\circ(X,\omega)=C_{HZ}^\circ(X,\omega)$ is finite.

(i) follows from the corollaries of Propositions 3 and 4 in
\cite{Be} and (ii) comes from Theorem 1.5 in \cite{LiuT}. We
conjecture that the conclusion also holds for the {\bf rationally
connected manifolds}\, introduced in \cite{KoMiMo}. }
\end{example}

In some cases we can get better results.

\begin{theorem}  \label{t:1.21}
For the standard symplectic form $\sigma_{n_i}$ on $\CP^{n_i}$ as
in Theorem~\ref{t:1.17} and any $a_i\in\R \setminus \{0\}$,
$i=1,\cdots, k$, we have
\[
C(\CP^{n_1}\times\cdots\times\CP^{n_k},
  a_1\sigma_{n_1}\oplus\cdots\oplus a_k\sigma_{n_k})
=(|a_1|+\cdots + |a_k|)\pi.
\]
for $C=c_{HZ}$ and $c_{HZ}^\circ$.
\end{theorem}

According to Example 12.5 of \cite{McSa1}
$$
{\cal W}_G(\CP^1\times\cdots\times\CP^1,
a_1\sigma_1\oplus\cdots\oplus a_k\sigma_1)
=\min\{|a_1|,\cdots,|a_k|\}\pi
$$
for any $a_i\in\R\setminus\{0\}$, $i=1,\cdots, k$. This,
Theorem~\ref{t:1.21} and (\ref{e:1.5}) show that $C_{HZ}$,
$C_{HZ}^\circ$, $c_{HZ}$ and $c_{HZ}^\circ$ are different from the
Gromov width ${\cal W}_G$.

\subsection{The Weinstein conjecture and periodic orbits near symplectic
submanifolds}

\noindent{\bf 1.4.1. Weinstein conjecture in cotangent bundles of
uniruled manifolds}.\quad By ``Weinstein conjecture'' we in the
sequel mean the $(pt, pt)$-Weinstein conjecture, i.e.: Every
separating hypersurface $S$ of contact type in a symplectic
manifold carries a closed characteristic.
While in some of the previous works on
the Weinstein conjecture, e.g.\ \cite{HV1}, the assumption that
$S$ is separating was also imposed, Weinstein's original
conjecture, \cite{We2}, does not assume that $S$ is separating.
So far this conjecture has been proved for many symplectic
manifolds, cf.\ \cite{C,FHV,
  FrSchl, H2,HV1,HV2,LiuT,Lu1,Lu2,Lu3,V1,V4,V5} and the recent
nice survey \cite{Gi} for more references. In particular, for the
Weinstein conjecture in cotangent bundles Hofer and Viterbo
\cite{HV1} proved that if a connected hypersurface $S$ of contact
type in the cotangent bundle of a closed manifold $N$ of dimension
at least 2 is such that the bounded component of $T^\ast
N\setminus S$ contains the zero section of $T^\ast N$, then it
carries a closed characteristic. In \cite{V5} it was proved that
the Weinstein conjecture holds in cotangent bundles of simply
connected closed manifolds. We shall prove

\begin{theorem}  \label{t:1.22}
Let $(M, \omega)$ be a closed connected symplectic manifold of
dimension at least $4$ and let $L\subset M$ be a Lagrangian
submanifold. Given a homology class $\tilde\alpha_0\in
H_\ast(L;\Q)\setminus\{0\}$ we denote by $\alpha_0\in
H_\ast(M;\Q)$ the class induced by the inclusion $L\hookrightarrow
M$. Assume that the Gromov-Witten invariant $\Psi_{A, g,
m+1}(C;\alpha_0,\alpha_1,\cdots,\alpha_m)$ does not vanish for
some homology classes $A\in H_2(M;\Z)$,
$\alpha_1,\dots,\alpha_m\in H_\ast(M;\Q)$ and $C\in
H_\ast(\overline{\cal M}_{g, m+1};\Q)$ and integers $m>1$ and
$g>0$. Then for every $c>0$,
$$
C_{HZ}^{(2)}(U_c,\omega_{\rm can}; \tilde\alpha_0, pt) < +\infty,
$$
and
$$
 C_{HZ}^{(2\circ)} (U_c,\omega_{\rm can};
\tilde\alpha_0, pt) < +\infty
$$
if $g=0$ and the inclusion $L\hookrightarrow M$ induces an
injective homomorphism $\pi_1(L)\to\pi_1(M)$. Here $U_c=\{(q,
v^\ast)\in T^\ast L\,|\, \langle v^\ast, v^\ast\rangle\le c^2\}$
is with respect to a Riemannian metric $\langle\cdot,\cdot\rangle$
on $T^\ast L$. Consequently, every hypersurface of contact type in
$(T^\ast L, \omega_{\rm can})$ separating $\tilde\alpha_0$ and
$pt$ carries a closed characteristic and a contractible one in the
latter case. In particular, if $(M,\omega)$ is a $g$-symplectic
uniruled manifold then for each $c>0$,
$$
  c_{HZ}(U_c,\omega_{\rm can})=C_{HZ}(U_c,\omega_{\rm can})<+\infty
  $$
and
\begin{equation}  \label{e:1.24}
c_{HZ}^{\circ}(U_c,\omega_{\rm can})
=C_{HZ}^{\circ}(U_c,\omega_{\rm can})<+\infty
  \end{equation}
if $g=0$ and the inclusion $L\hookrightarrow M$ induces an
injective homomorphism $\pi_1(L)\to\pi_1(M)$. If $(M,\omega)$
itself is strong symplectic uniruled then (\ref{e:1.24}) also
holds for $L=M \subset (T^*M, \omega_{\rm can})$.
\end{theorem}

Using a recent refinement by Macarini and Schlenk \cite{MaSchl} of
the arguments in [HZ2, Sections 4.1 and 4.2] we immediately
derive:
  if $L$ is a  Lagrangian submanifold in a $g$-symplectic uniruled
manifold and $S\subset (T^\ast L,\omega_{\rm can})$ a smooth
compact connected orientable hypersurface without boundary, then
for any thickening of $S$,
\begin{eqnarray*}
&&\hspace{3cm}\psi: I\times S\to U\subset (T^\ast L,\omega_{\rm can})\\
&&\mu\{t\in I\,|\,{\mathcal P}(S_t)\ne\emptyset\}=\mu(I)\quad{\rm
and}\quad \mu\{t\in I\,|\,{\mathcal
P}^\circ(S_t)\ne\emptyset\}=\mu(I)
\end{eqnarray*}
if $g=0$ and
the inclusion $L\hookrightarrow M$ induces an injective
homomorphism $\pi_1(L)\to\pi_1(M)$. Here $\mu$ denotes Lebesgue
measure, $I$ is an open neighborhood of $0$ in $\R$, and
${\mathcal P}(S_t)$ (resp.\ ${\cal P}^\circ(S_t)$) denotes the set
of all (resp.\ contractible in $U$) closed characteristics on
$S_t=\psi(S\times\{t\})$.

\begin{corollary}  \label{cor:1.23}
The Weinstein conjecture holds in the following manifolds:
\begin{description}
\item[(i)] symplectic uniruled manifolds of dimension at least
$4$;

\item[(ii)] the cotangent bundle $(T^\ast L,\omega_{\rm can})$ of
a closed Lagrangian submanifold $L$ in a $g$-symplectic uniruled
manifold of dimension at least $4$;

\item[(iii)] the product of a closed symplectic manifold and a
strong symplectic uniruled manifold;

\item[(iv)] the cotangent bundles of  strong symplectic uniruled
manifolds.
\end{description}
\end{corollary}

The result in (i) is actually not new. As observed in \cite{Lu3}
the Weinstein conjecture in symplectic uniruled manifolds can be
derived from Theorem 1.1 in \cite{LiuT}. With the present
arguments  it may be derived from (\ref{e:1.18}) and
Corollary~\ref{cor:1.12}. (ii) is a direct consequence of
Theorem~\ref{t:1.22}. (iii) can be derived from (i) and
Proposition~\ref{prop:7.5}. By (ii) and Proposition~\ref{prop:7.5}
the standard arguments give rise to (iv).

\bigskip
\noindent{\bf 1.4.2. Periodic orbits near symplectic
submanifolds}.\quad The existence of periodic orbits of autonomous
Hamiltonian systems near a closed symplectic submanifold has been
studied by several authors, see \cite{CiGiKe, GiGu, Ke} and the
references there for details. Using Proposition~\ref{prop:1.7} and
suitably modifying the arguments in \cite{Lu6} and \cite{Bi1} we
get

\begin{theorem}  \label{t:1.24}
Let $(M, \omega)$ be any symplectic manifold and let $N \subset M$ be a
connected closed symplectic submanifold of
codimension $2$.  Then for any $\varepsilon>0$ there exists
a smooth compact submanifold $W\subset M$ with connected boundary
and of codimension zero which is a neighborhood of $N$ in $M$ such
that
$$
c_{HZ}^\circ(W,\omega)=C_{HZ}^\circ(W,\omega)<\varepsilon.
$$
Consequently, for any  smooth compact connected orientable
hypersurface $S\subset W\setminus\partial W$ without boundary and
any thickening $\psi \colon S\times I\to U \subset W$ it holds
that
$$\mu (\{t\in I\,|\,{\cal
P}^\circ(S_t)\ne\emptyset\})=\mu(I).
$$
  Here  $\mu$, $I$,  $S_t$ and ${\cal P}^\circ(S_t)$  are  as
  above Corollary~\ref{cor:1.23}.
\end{theorem}

The first conclusion will be proved in $\S5$, and the second
follows from the first one and the refinement of the Hofer-Zehnder
theorem by Macarini and Schlenk \cite{MaSchl} mentioned above. The
second conclusion in Theorem~\ref{t:1.24} implies: For any smooth
proper function $H:W\to \R$ the levels $H=\epsilon$ carry
contractible in $U$ periodic orbits for almost all $\epsilon>0$
for which $\{H=\epsilon\}\subset {\rm Int}(W)$. Using Floer
homology and symplectic homology, results similar to
Theorem~\ref{t:1.24} were obtained in \cite{CiGiKe, GiGu} for some
closed symplectic submanifolds of positive codimension in
geometrically bounded, symplectically aspherical manifolds. Recall
that a symplectic manifold $(M,\omega)$ is said to be
symplectically aspherical if $\omega|_{\pi_2(M)}=0$ and
$c_1(TM)|_{\pi_2(M)}=0$.
It seems possible that our method can be
generalized to any closed symplectic submanifold of codimension
more than $2$.

\subsection{Nonsqueezing theorems}
We first give a general nonsqueezing theorem and then discuss some
corollaries and relations to the  various previously found
nonsqueezing theorems.

\begin{definition}  \label{def:1.25}
{\rm For a symplectic manifold $(M,\omega)$ we define
$\Gamma(M,\omega)\in [0, +\infty]$ by
$$
\Gamma(M,\omega)=\inf_{\alpha} C_{HZ}^{(2)}(M,\omega; pt,\alpha),
$$
where $\alpha\in H_\ast(M;\Q)$ runs over all nonzero homology
classes of degree $\deg\alpha \le \dim M-1$.}
\end{definition}

By (\ref{e:1.12}), for any  connected symplectic manifold
$(M,\omega)$  we have
\begin{equation}  \label{e:1.25}
{\cal W}_G(M,\omega)\le\Gamma(M,\omega).
\end{equation}
However, it is difficult to determine or estimate
$\Gamma(M,\omega)$. In some cases one can replace it by another
number.\vspace{2mm}

\begin{definition}  \label{def:1.26}
{\rm For a  closed connected symplectic manifold $(M,\omega)$ of
dimension at least $4$ we define ${\rm GW}(M,\omega)\in (0,
+\infty ]$ by
$$
{\rm GW}(M,\omega)=\inf{\rm GW}_g(M,\omega; pt,\alpha)
$$
where the infimum is taken over all nonnegative integers $g$ and
all homology classes $\alpha\in H_\ast(M;\Q)\setminus\{0\}$ of
degree $\deg\alpha\le\dim M-1$. }
\end{definition}

By (\ref{e:1.17}) we have ${\rm GW}(M,\omega)=\inf_g{\rm
GW}_g(M,\omega; pt, PD([\omega]))$. Note that ${\rm GW}(M,\omega)$
is finite if and only if $(M,\omega)$ is a symplectic uniruled
manifold. From Theorem~\ref{t:1.13} and (\ref{e:1.25}) we get

\begin{theorem}  \label{t:1.27}
For any symplectic uniruled manifold $(M,\omega)$ of dimension at
least $4$ we have
$${\cal W}_G(M,\omega)\le {\rm GW}(M,\omega).$$
\end{theorem}

Actually, for a uniruled manifold $(M,\omega)$, i.e., a K\"ahler
manifold covered by rational curves, the arguments in \cite{Ko,R}
show that ${\rm GW}(M,\omega)\le\omega(A)$, where $A=[C]$ is the
class of a rational curve $C$ through a generic $x_0\in M$ and
such that $\int_C\omega$ is minimal.\vspace{2mm}
\begin{remark}  \label{rem:1.28}
{\rm Denote by $(W,\Omega)$ the product
$$
(\CP^{n_1}\times\cdots\times\CP^{n_k},
a_1\sigma_{n_1}\oplus\cdots\oplus a_k\sigma_{n_k})
$$
in Theorem~\ref{t:1.21}. It follows from Theorem~\ref{t:1.13} and
the proof of Theorem~\ref{t:1.17} that
$${\rm GW}(W,\Omega)\le
\min\{|a_1|,\cdots, |a_k|\}\pi.$$ By (\ref{e:1.25}) and definition
of $\Gamma(W,\Omega)$, for any small $\epsilon>0$ there exists a
class $\alpha_\epsilon\in H_\ast(W,\Q)$ of degree ${\rm
deg}(\alpha_\epsilon)\le\dim W-1$ such that
$${\cal W}_G(W,\Omega)\le C_{HZ}^{(2)}(W,\Omega;
pt,\alpha_\epsilon)<\min\{|a_1|,\cdots, |a_k|\}\pi+\epsilon.$$ But
Theorem~\ref{t:1.21} shows that
$$c_{HZ}(W,\Omega)=C_{HZ}(W,\Omega)= (|a_1|+\cdots +|a_k|)\pi.$$
Therefore, if $k>1$ and $\epsilon>0$ is small enough then
$${\cal W}_G(W,\Omega)\le C_{HZ}^{(2)}(W,\Omega;
pt,\alpha_\epsilon)<C_{HZ}(W,\Omega).$$
  This shows that our pseudo symplectic capacity $C_{HZ}^{(2)}(W,\Omega;
pt,\alpha_\epsilon)$ can give a better upper bound for ${\cal
W}_G(W,\Omega)$ than the symplectic capacities $c_{HZ}(W,\Omega)$
and $C_{HZ}(W,\Omega)$.}
\end{remark}

Recall that Gromov's famous nonsqueezing theorem states that if
there exists a symplectic embedding $B^{2n}(r)\hookrightarrow
Z^{2n}(R)$, then $r\le R$. Gromov proved it by using $J$-holomorphic
curves, \cite{Gr}. Later on, proofs were given by Hofer and Zehnder
based on the calculus of variation and by Viterbo using generating
functions, \cite{V3}. As a direct consequence of Theorem~\ref{t:1.5}
and (\ref{e:1.23}) we get

\begin{corollary}  \label{cor:1.29}
For any closed symplectic manifold $(M,\omega)$ of dimension $2m$,
if there exists a symplectic embedding
$$
B^{2m+2n}(r)\hookrightarrow (M\times Z^{2n}(R),
\omega\oplus\omega_0),
$$
then $r\le R$.
\end{corollary}

Actually, Lalonde and McDuff proved Corollary~\ref{cor:1.29} for
any symplectic manifold $(M, \omega)$ in \cite{LaMc1}. Moreover,
one can derive from it the foundational energy-capacity inequality
in  Hofer geometry (cf.\ \cite{LaMc1,La2} and \cite[Ex.\
12.21]{McSa1}). From (\ref{e:1.23}) one can also derive the
following version of the non-squeezing theorem which was listed
below Corollary 5.8 of [LaMc2,II] and which can be used to prove
that the group of Hamiltonian diffeomorphisms of some compact
symplectic manifolds have infinite diameter with respect to
Hofer's metric.\vspace{2mm}

\begin{corollary}  \label{cor:1.30}
Let $(M,\omega)$ and $(N,\sigma)$ be  closed symplectic manifolds
of dimensions $2m$ and $2n$ respectively. If there exists a
symplectic embedding
$$
M \times B^{2n+2p}(r) \hookrightarrow ( M\times N\times B^{2p}(R),
\omega\oplus\sigma\oplus\omega_0^{(p)})
$$
or a symplectic embedding
$$
M \times B^{2n+2p}(r) \hookrightarrow (M \times \R^{2n} \times
B^{2p}(R), \omega\oplus\omega_0^{(n)} \oplus\omega_0^{(p)}) ,
$$
then $r\le R$. Here $\omega_0^{(m)}$ denotes the standard
symplectic structure on $\R^{2m}$.
\end{corollary}

The second statement can be reduced to the first one. From
Theorem~\ref{t:1.16} we get

\begin{corollary}  \label{cor:1.31}
For any closed symplectic manifold $(M,\omega)$ of dimension $2m$,
\[
{\cal W}_G \left( M \times G(k,n), \omega\oplus(a\sigma^{(k,n)})
\right) \le |a|\pi .
\]
\end{corollary}

The study of Hofer geometry requires various nonsqueezing
theorems. Let us recall the notion of quasicylinder introduced by
Lalonde and McDuff in \cite{LaMc2}.

\begin{definition}  \label{def:1.32}
{\rm For a closed symplectic manifold $(M,\omega)$ and a set $D$
diffeomorphic to a closed disk in $(\R^2,\omega_0=ds\wedge dt)$, the
manifold $Q=(M\times D,\Omega)$ endowed with the symplectic form
$\Omega$ is called a {\bf quasicylinder} if
\begin{description}
\item[(i)] $\Omega$ restricts to $\omega$ on each fibre
$M\times\{pt\}$; \item[(ii)] $\Omega$ is the product
$\omega\times\omega_0$ near the boundary
             $\partial Q=M\times\partial D$.
\end{description}
If $\Omega=\omega\times\omega_0$ on $Q$, the quasicylinder is called
{\bf split}. The {\bf area} of a quasicylinder $(M\times D,\Omega)$
is defined as the number $\Lambda=\Lambda(M\times D,\Omega)$ such
that
$${\rm Vol}(M\times D, \Omega)=\Lambda\cdot{\rm
Vol}(M,\omega).$$ }
\end{definition}

As proved in Lemma~2.4 of \cite{LaMc2}, the area $\Lambda(M\times
D,\Omega)$ is equal to $\int_{\{x\}\times D}\Omega$ for any $x\in
M$.

Following \cite{McSl} we replace $Q$ in Definition~\ref{def:1.32} by
the obvious \linebreak $S^2$-compactification $(M\times S^2,
\Omega)$. Here $\Omega$ restricts to $\omega$ on each fibre. It is
clear that $\Omega(A)=\Lambda(Q,\Omega)$ for $A=[pt\times S^2]\in
H_2(M\times S^2)$. But it is proved in Lemma 2.7 of \cite{LaMc2}
that $\Omega$ can be symplectically deformed to a product symplectic
form $\omega\oplus\sigma$. Therefore, it follows from the
deformation invariance of Gromov-Witten invariants that
$$\Psi_{A,0,3}(pt; pt,[M\times pt], [M\times pt])\ne 0.$$
By Theorem~\ref{t:1.13} we get
$$
C(M\times S^2, \Omega; pt, [M\times
pt])\le\Omega(A)=\Lambda(Q,\Omega).
$$
As in the proof of Theorem~\ref{t:1.17} we can derive from this

\begin{theorem}  \label{t:1.33}
{\bf (Area-capacity inequality)} For any quasicylinder
$(Q,\Omega)$
$$
c^\circ_{HZ}(Q,\Omega)=C_{HZ}^{\circ}(Q,\Omega)\le\Lambda(Q,\Omega).
$$
\end{theorem}

Area-capacity inequalities for ${\cal W}_G$, $c_{HZ}$ and
$c_{HZ}^\circ$ have been studied in \cite{FHV,HV2,LaMc1,Lu1,McSl}.
As in \cite{LaMc2,McSl} we can use Theorem~\ref{t:1.33} and
Lemma~\ref{lem:1.4} to deduce the main result in \cite{McSl}: For
an autonomous Hamiltonian $H \colon M \to \R$ on a closed
symplectic manifold $(M,\omega)$ of dimension at least $4$, if its
flow has no nonconstant contractible  fast periodic solution  then
the path $\phi^H_{t\in [0, 1]}$ in ${\rm Ham}(M,\omega)$ is
length-minimizing among all paths homotopic with fixed endpoints.

 From Theorem~\ref{t:1.33} and (\ref{e:1.5}) we obtain the
following non-squeezing theorem for quasi-cylinders.

\begin{corollary}  \label{cor:1.34}
For any quasicylinder $(M\times D,\Omega)$ of dimension $2m+2$,
\[
{\cal W}_G \left( M \times D, \Omega \right) \le \Lambda(M\times
D,\Omega).
\]
\end{corollary}

Our results also lead to the nonsqueezing theorem Proposition 3.27
in \cite{Mc2} for Hamiltonian fibrations $P\to S^2$.

\subsection{Symplectic packings and Seshadri constants}

{\bf 1.6.1. Symplectic packings}.\hspace{2mm} Suppose that
$B^{2n}(r)=\{z\in\R^{2n}\,|\, |z| < r\}$ is endowed with the
standard symplectic structure $\omega_0$ of $\R^{2n}$. For an
integer $k>0$, a {\bf symplectic} $k$-{\bf packing} of a
$2n$-dimensional symplectic manifold $(M,\omega)$ via $B^{2n}(r)$ is
a set of symplectic embeddings $\{\varphi_i\}^k_{i=1}$ of
$(B^{2n}(r), \omega_0)$ into $(M,\omega)$ such that ${\rm
Im}\varphi_i\cap {\rm Im}\varphi_j=\emptyset$ for $i\ne j$. If ${\rm
Vol}(M,\omega)$ is finite and $\Int (M)\subset \overline{\cup{\rm
Im}\varphi_i}$, then $(M, \omega)$ is said to have a {\bf full
symplectic $k$-packing}. Symplectic packing problems were studied
for the first time by Gromov in \cite{Gr} and later by McDuff and
Polterovich \cite{McPo}, Karshon \cite{Ka}, Traynor \cite{Tr}, Xu
\cite{Xu}, Biran \cite{Bi1,Bi2} and Kruglikov \cite{Kru}. As before,
let $\sigma_n$ denote the unique ${\rm U}(n+1)$-invariant K\"ahler
form on $\CP^n$ whose integral over $\CP^1$ is equal to $\pi$. For
every positive integer $p$, a full symplectic $p^n$-packing of
$(\CP^n, \sigma_n)$ was explicitly constructed by McDuff and
Polterovich \cite{McPo} and Traynor \cite{Tr}. A direct geometric
construction of a full symplectic $n+1$-packing of
$(\CP^n,\sigma_n)$ was given by Yael Karshon, \cite{Ka}. By
generalizing the arguments in \cite{Ka} we shall obtain

\begin{theorem}  \label{t:1.35}
Let the Grassmannian $(G(k, n), \sigma^{(k, n)})$ be as in
Theorem~\ref{t:1.15}. Then for every integer $1<k<n$ there exists
a symplectic $[n/k]$-packing of $(G(k, n), \sigma^{(k, n)})$ by
$B^{2k(n-k)}(1)$. Here $[n/k]$ denotes the largest integer less
than or equal to $n/k$.
\end{theorem}

This result shows that the Fefferman invariant of $(G(k, n),
\sigma^{(k, n)})$ is at least $[n/k]$. Recall that the Fefferman
invariant $F(M,\omega)$ of a $2n$-dimensional symplectic manifold
$(M,\omega)$ is defined as the largest integer $k$ for which there
exists a symplectic packing by $k$ open unit balls. Moreover, at
the end of \S6 we shall prove
\begin{equation}\label{e:1.26}
\quad  {\rm Vol}(G(k, n),\sigma^{(k,n)})=\frac{(k-1)!\cdots
2!\cdot 1!\cdot (n-k-1)! \cdots 2!\cdot 1!}{(n-1)!\cdots 2!\cdot
1!}\cdot\pi^{k(n-k)}.
\end{equation}
Note that ${\rm Vol}(B^{2k(n-k)}(1),
\omega_0)=\pi^{k(n-k)}/(k(n-k))!$. One easily sees that the
symplectic packings in Theorem~\ref{t:1.35} is not full in
general. On the other hand a full packing of each of the
Grassmannians $Gr^+(2,\R^5)$ and $Gr^+(2,\R^6)$ by two equal
symplectic balls was constructed in \cite{KaTo2}.

\noindent{\bf 1.6.2. Seshadri constants.} Our previous results can
also be used to estimate Seshadri constants, which are interesting
invariants in algebraic geometry. Recall that for a compact
complex manifold $(M, J)$  of complex dimension $n$ and an ample
line bundle $L\to M$, the {\bf Seshadri constant} of $L$ at a
point $x\in M$ is defined as the nonnegative real number
\begin{equation}\label{e:1.27}
\varepsilon(L, x):=\inf_{C\ni x}\frac{\int_C c_1(L)}{{\rm
mult}_xC},
\end{equation}
  where the infimum is taken over all irreducible holomorphic curves $C$
passing through the point $x$, and ${\rm mult}_xC$ is the
multiplicity of $C$ at $x$ (\cite{De}). The global Seshadri
constant is defined by
\[
\varepsilon(L) := \inf_{x\in M} \varepsilon(L, x).
\]
Seshadri's criterion for ampleness says that $L$ is ample if and
only if $\varepsilon(L)>0$. The cohomology class $c_1(L)$ can be
represented by a $J$-compatible K\"ahler form $\omega_L$ (the
curvature form for a suitable metric connection on $L$). Denote by
$L^n=\int_M\omega_L^n=n!{\rm Vol}(M,\omega_L)$. Then
$\varepsilon(L, x)$ has the elementary upper bound
\begin{equation}  \label{e:1.28}
\varepsilon(L, x) \le \sqrt[n]{L^n}.
\end{equation}
Biran and Cieliebak \cite[Prop.\ 6.2.1]{BiCi} gave a better upper
bound, i.e.
\[
\varepsilon(L) \le {\cal W}_G(M,\omega_L) .
\]
However, it is difficult to estimate ${\cal W}_G(M,\omega_L)$.
Together with Theorem~\ref{t:1.27} we get

\begin{theorem}  \label{t:1.36}
For a closed connected complex manifold of complex dimension at
least $2$,
\[
\varepsilon(L) \le {\rm GW}(M,\omega_L).
\]
\end{theorem}

\begin{remark}  \label{rem:1.37}
{\rm By Definition~\ref{def:1.26}, if ${\rm GW}(M,\omega_L)$ is
finite then $(M,\omega_L)$ is symplectic uniruled. So
Theorem~\ref{t:1.36} has only actual sense for uniruled $(M, J)$.
In this case our upper bound ${\rm GW}(M,\omega_L)$ is better than
$\sqrt[n]{L^n}$ in (\ref{e:1.28}). As an example, let us consider
the hyperplane $[H]$ in $\CP^n$. It is ample, and the Fubini-Study
form $\omega_{\rm FS}$ with $\int_{\CP^1}\omega_{\rm FS}=1$ is a
K\"ahler representative of $c_1([H])$. Let $p_1$ and $p_2$ denote
the projections of the product $\CP^n\times\CP^n$ to the first and
second factors. For an integer $m>1$ the line bundle $p_1^\ast[H]+
p_2^\ast(m[H])\to \CP^n\times\CP^n$ is ample and $c_1(p_1^\ast[H]+
p_2^\ast(m[H]))$ has a K\"ahler form representative $\omega_{\rm
FS}\oplus m\omega_{\rm FS}$. From the proof of
Theorem~\ref{t:1.16} it easily follows that
$${\rm GW}(\CP^n\times\CP^n,\omega_{\rm FS}\oplus m\omega_{\rm FS})\le 1.$$
(In fact, equality holds.) But a direct computation gives
$$
\sqrt[2n]{(p_1^\ast[H]+
p_2^\ast(m[H]))^{2n}}=\Bigl(\int_{\CP^n\times\CP^n} (\omega_{\rm
FS}\oplus m\omega_{\rm
FS})^{2n}\Bigr)^{\frac{1}{2n}}\sqrt[2n]{m}\cdot\sqrt[2n]{\frac{(2n)!}{n!n!}}
 >1 .
$$
}
\end{remark}

 From the above arguments and the subsequent proofs the reader can
see that some of our results are probably not optimal. In fact, it
is very possible that using our methods one can obtain better
results in some cases (\cite{Lu7} and \cite{Lu9}). We content
ourselves with illustrating the new ideas and methods.
\vspace{2mm}

The paper is organized as follows. In Section~2 we give the proofs
of Lemma~\ref{lem:1.4}, Theorems \ref{t:1.5}, \ref{t:1.8} and
Proposition~\ref{prop:1.7}. The proof of Theorem~\ref{t:1.10} is
given in Section 3. In Section~4 we prove Theorems~\ref{t:1.15},
\ref{t:1.16}, \ref{t:1.17} and \ref{t:1.21}. In Section 5 we prove
Theorems~\ref{t:1.22}, \ref{t:1.24}. Theorem~\ref{t:1.35} is
proved in Section 6. In the Appendix we discuss some related
results on the Gromov-Witten invariants of product manifolds.
\vspace{2mm}

\noindent{\bf Acknowledgements.} I am deeply grateful to Professor
Dusa McDuff for  kindly explaining to me some details of
\cite{McSl}, sending me a preliminary version of her preprint
\cite{Mc3}, correcting many errors, giving many comments and
improving some results in the preliminary version of this paper.
The author is  very thankful to  Dr.\ Felix Schlenk for kindly
checking this paper, giving many useful comments and correcting
many errors in mathematics and English. I wish to thank Professor
Claude Viterbo for many useful discussions on the symplectic
capacity theory during the author's visit at IHES in  spring 1999.
I thank Kai Cieliebak for sending me some preprints. I would like
to thank Professor L.\ Polterovich and the referee for their
advice and for improving the presentation, pointing out some
errors and gaps in the previous version of the paper and for
providing some enlightening techniques in some references. The
example in Remark~\ref{rem:1.6} is essentially of them. This
revised version was finished during my visit at ICTP of Italy and
IHES at Paris in 2004. I thank Professors D.T. L\^e and J. P.
Bourguignon for their invitations and for the financial support
and hospitality of both institutions. I also thank Leonardo
Macarini and Felix Schlenk for pointing out an error in the
original proof of Theorem~\ref{t:1.24}, and Paul Biran for his
friendly helps in revising proof of it.

\section{Proofs of Lemma \ref{lem:1.4}, Theorems~\ref{t:1.5}, \ref{t:1.8}
and Proposition~\ref{prop:1.7}}

We first give two lemmas. They are key to our proofs in this
section the next one. According to Lemma 4.4 on page 107 and
Exercise~9 on page~108 of \cite{Hi} we have:

\begin{lemma}\label{lem:2.1}
If $N$ is a connected smooth manifold and $W\subset {\it Int}(N)$ a
compact smooth submanifold with connected boundary and of
codimension zero, then $\partial W$ separates $N$ in the sense that
${\it Int}(N)\setminus\partial W$ has exactly two connected
components and the topological boundary of each component is
$\partial W$. In this case $\partial W$ has a neighborhood in $N$
which is a product $\partial W\times (-2, 2)$ with $\partial W$
corresponding to $\partial W\times\{0\}$. If $W$ is only contained
in $N$ then $\partial W$ has a neighborhood in $W$ which is a
product $\partial W\times (-2,0]$.
\end{lemma}

 From Lemma 12.27 in \cite{McSa1} we easily derive

\begin{lemma}  \label{lem:2.2}
Given a Riemannian metric $g$ on $M$, there exists $\rho = \rho
(g,M) >0$ such that for every smooth function $H$ on $M$ with
$$\sup_{x\in M}\|\nabla_g\nabla_g H(x)\|_g<\rho$$
  the Hamiltonian equation $\dot x=X_H(x)$ has
no nonconstant fast periodic solutions. In particular, the
conclusion holds if $\|H\|_{C^2} < \rho$. Here $\nabla_g$ is the
Levi-Civita connection of $g$ and norms are taken with respect to
$g$.
\end{lemma}

 From Darboux's theorem we obtain

\begin{lemma}\label{lem:2.3} Let $(M,\omega)$ be a
$2n$-dimensional symplectic manifold, and
$B^{2n}(r)\\=\{z\in\R^{2n} : |z|\le r\}$ with $r>0$. Then for any
$z_0\in{\rm Int}(M)$ and any small $\varepsilon>0$ there exist
$r>0$, a symplectic embedding $\varphi:(B^{2n}(2r),\omega_0)\to
(M,\omega)$ with $\varphi(0)=z_0$ and a smooth function
$H^\varphi_{r,\varepsilon}:M\to\R$ such
  that:
\begin{description}
\item[(i)] $H^\varphi_{r,\varepsilon}=0$ outside ${\rm
Int}(\varphi(B^{2n}(2r))$, and
$H^\varphi_{r,\varepsilon}=\varepsilon$ on $\varphi(B^{2n}(r))$.

\item[(ii)] $H^\varphi_{r,\varepsilon}$ is constant $h(s)$ along
  $\varphi(\{|z|=s\})$ for any $s\in [0,2r]$, where $h:[0, 2r]\to
  [0,\varepsilon]$ is a  nonnegative smooth function which is
strictly decreasing
  on $[r, 2r]$. Consequently,
$H^\varphi_{r,\varepsilon}(\varphi(z))>H^\varphi_{r,\varepsilon}(\varphi(z'))$
if $r\le |z|<|z'|\le 2r$,
  and $H^\varphi_{r,\varepsilon}$ has no critical values in
$(0,\varepsilon)$.

  \item[(iii)] $\dot x=X_{H^\varphi_{r,\varepsilon}}(x)$ has no
nonconstant fast periodic solutions.

\end{description}
\end{lemma}

\noindent{\bf Proof of Lemma 1.4.} {\bf Case (i).} We only need to
prove that
$$C_{HZ}(M,\omega;pt,pt)\ge c_{HZ}(M,\omega).$$
To this end it suffices to construct for any $H\in{\mathcal
H}_{ad}(M,\omega)$ an
$$
F\in {\cal H}_{ad}(M,\omega; pt, pt)
$$
such that $\max F\ge\max H$.
By the definition  there exist a nonempty open subset $U$ and a
compact subset $K\subset M\setminus\partial M$ such that: (a)
$H|_U=0$ and $H|_{M\setminus K}=\max H$, (b) $0\le H\le\max H$,
(c) $\dot x=X_H(x)$ has no nonconstant fast periodic solutions.
These imply that $U\subset{\rm Int}(K)$. By the illustrations
below the definition of $c_{HZ}$ in \S1.2 we may assume that
$M\setminus K\ne\emptyset$. Then both $U$ and $M\setminus K$ are
nonempty open sets because $M$ is a closed manifold.
For a given small $\varepsilon>0$ we may take symplectic
embeddings $\varphi$ and $\psi$ from $(B^{2n}(2r),\omega_0)$ to
$(M,\omega)$ such that
$$\varphi(B^{2n}(2r))\subset U\quad{\rm
and}\quad\psi(B^{2n}(2r))\subset M\setminus K.
$$
Let $H^\varphi_{r,\varepsilon}$ and $H^\psi_{r,\varepsilon}$ be
the corresponding functions as in Lemma~\ref{lem:2.3}. Since
$H^\varphi_{r,\varepsilon}$ (resp.\ $H^\psi_{r,\varepsilon}$) is
equal to zero outside $\varphi(B^{2n}(2r))$ (resp.
$\psi(B^{2n}(2r))$)  we can define a smooth function $\widetilde H
\colon M\to\R$ by
$$
  \widetilde H(x)=\left\{\begin{array}{llc}
  \max H+ H^\psi_{r,\varepsilon}(x) & {\rm if}& x\in
  M\setminus K,\\
H(x) & {\rm if} & x\in K\setminus U, \\
-H^\varphi_{r,\varepsilon}(x) & {\rm if} &x\in U.
\end{array}\right.
$$
   Define $F=\widetilde H+\varepsilon$. Then
  $\max F=\max H+2\varepsilon\ge\max H$, $\min F=0$ and  $\dot
x=X_F(x)$ has no nonconstant fast periodic orbits in $M$.

  Since $M$ is a
closed manifold,  $M\setminus{\rm Int}(\psi(B^{2n}(r)))$ is a
compact submanifold with boundary $\psi(\partial B^{2n}(r))$. It
follows that $F\in{\mathcal H}_{ad}(M,\omega;pt,pt)$ with
$P(F)=\varphi(B^{2n}(r))$ and $Q(F)=M\setminus{\rm
Int}(\psi(B^{2n}(r)))$. The desired result follows.

Going through the above proof we see that if $H\in{\mathcal
H}^\circ_{ad}(M,\omega)$, i.e.,\ $\dot x=X_H(x)$ has no
nonconstant contractible fast periodic solutions, then
$F\in{\mathcal H}^\circ_{ad}(M,\omega;pt,pt)$. This implies that
$C_{HZ}^\circ(M,\omega)=c_{HZ}^\circ(M,\omega)$.

\noindent{\bf Case (ii).} The arguments are similar. We only point
out different points. Let $H\in {\mathcal H}_{ad}(M,\omega)$. For
a compact subset $K(H)\subset M\setminus\partial M$ we find by
assumption a compact submanifold $W$ with connected  boundary and
of codimension zero such that $K(H)\subset W$. Since $K(H)$ is
compact and disjoint from $\partial M$ we can assume that $K(H)$
is also disjoint from $\partial W$. By Lemma~\ref{lem:2.1} we can
choose embeddings
$$
\Phi \colon [-5, 0]\times \partial W\to M
$$
such that $\Phi(\{0\}\times\partial W)=\partial W$ and that
$$
  \Phi([-5, 0]\times \partial
W)\subset W\quad{\rm and}\quad K(H)\cap\Phi([-5, 0]\times \partial
W)=\emptyset.
$$
  For each $t\in [-5, 0]$ the set
$$W_t:=W\setminus\Phi((t, 0]\times\partial W)$$
is a compact submanifold of $M$ which is diffeomorphic to $W$. By
shrinking $\varepsilon>0$
  in Case (i) if necessary, one easily constructs a smooth  function
$H_\varepsilon:M\to\R$ such that
  \begin{description}
  \item[(a)] $H_\varepsilon=0$ in ${\rm Int}(W_{-4})$ and
  $H_\varepsilon=\varepsilon$ outside $W_{-1}$;
  \item[(b)] $0\le H_\varepsilon\le\varepsilon$ and each $c\in
  (0,\varepsilon)$ is a regular value of $H_\varepsilon$;

  \item[(c)] $H_{\varepsilon}$ is constant $f(s)$ along
  $\Phi(\{s\}\times\partial W)$ for any $s\in [-5, 0]$, where $f:[-5, 0]\to
  [0,\varepsilon]$ is a  nonnegative smooth function which is strictly
increasing
  in $[-4, -1]$.
  \item[(d)] $\dot x=X_{H_\varepsilon}(x)$ has no
nonconstant fast periodic solutions.
\end{description}
Let $H^\varphi_{r,\varepsilon}$ be as in Case (i). We can define a
smooth function $\widetilde H:M\to\R$ by
  $$
  \widetilde H(x)=\left\{\begin{array}{lc}
  \max H+ H_{\varepsilon}(x)\quad\;\; {\rm if}\; x\in
  M\setminus K,\\
H(x)\hspace{21mm} {\rm if}\;x\in K\setminus U, \\
-H^\varphi_{r,\varepsilon}(x)\hspace{15mm} {\rm if}\;x\in U,
\end{array}\right.
$$
and set $F=\widetilde H +\varepsilon$. Then $\max F\ge\max H$,
$\min F=0$ and  $\dot x=X_F(x)$ has no nonconstant fast periodic
solutions.  As in Case (i) one checks that $F\in{\mathcal
H}_{ad}(M,\omega;pt,pt)$ with $P(F)=\varphi(B^{2n}(r))$ and
$Q(F)=M\setminus{\rm Int}(W_{-1})$. So we have $\max H\le\max F\le
C_{HZ}(M,\omega)$ for any $H\in{\mathcal H}_{ad}(M,\omega)$, and
thus $c_{HZ}(M,\omega)\le C_{HZ}(M,\omega)$. As above we get that
$C_{HZ}(M,\omega)=c_{HZ}(M,\omega)$ and
$C_{HZ}^\circ(M,\omega)=c_{HZ}^\circ(M,\omega)$.
   \hfill$\Box$\vspace{2mm}

\noindent{\bf Proof of Theorem~\ref{t:1.5}.} (i) We take
$H\in{\cal H}_{ad}(M,\omega;\alpha_0,\alpha_\infty)$. Let $P=P(H)$
and $Q=Q(H)$ be the corresponding submanifolds in
Definition~\ref{def:1.2}, and $\alpha_0$, $\alpha_\infty$ the
chain representatives. Define $G=-H+\max H$. Then $0\le G\le\max
G=\max H$, $G|_P=\max G$, $G|_{M\setminus \Int (Q)}=0$ and
$X_G=-X_H$. Therefore, $G\in {\cal
H}_{ad}(M,\omega;\alpha_\infty,\alpha_0)$, and (i) follows.

(ii) is a special case of (v), and (iv) and (vi) are clear.

For (iii), note that $B^{2n}(1)$ and $Z^{2n}(1)$ are contractible.
One can slightly modify the proofs of Lemma~3 and Theorem~2 in
Chapter~3 of \cite{HZ2} to show that
$C_{HZ}^{(2)}(B^{2n}(1),\omega_0;\alpha_0,\alpha_\infty)\ge\pi$
and
$C_{HZ}^{(2\circ)}(Z^{2n}(1),\omega_0;\alpha_0,\alpha_\infty)\le\pi$.
Then (iii) follows from (v) and definitions:
\begin{eqnarray*}
\pi&\le& C_{HZ}^{(2)}(B^{2n}(1),\omega_0;\alpha_0,\alpha_\infty)\\
&\le& C_{HZ}^{(2\circ)}(B^{2n}(1),\omega_0;\alpha_0,\alpha_\infty)\\
&\le & C_{HZ}^{(2\circ)}(Z^{2n}(1),\omega_0;\alpha_0,\alpha_\infty)\\
&\le&\pi.
\end{eqnarray*}

For (v) we only prove the first claim. The second claim then
follows together with the argument in \cite{Lu1}. For $H\in{\cal
H}_{ad}(M_1,\omega_1;\alpha_0, \alpha_\infty)$ let the
submanifolds $P_1$ and $Q_1$ of $(M_1,\omega_1)$ be as in
Definition~\ref{def:1.2}. Set $P_2=\psi(P_1)$ and $Q_2=\psi(Q_1)$,
and define $\psi_\ast(H)\in C^\infty(M_2,\R)$ by
\[
\psi_\ast(H)(x)=\left\{\begin{array}{ll}
H\circ\psi^{-1}(x)&{\rm if}\; x\in\psi(M_1)\\
\max H &{\rm if}\; x\notin\psi(M_1).\end{array}\right.
\]
It is clear that $\psi_\ast(H)\in {\cal H}_{\rm
ad}(M_2,\omega_2;\psi_\ast(\alpha_0), \psi_\ast(\alpha_\infty))$,
and so (v) follows.

To prove (vii) we only need to show that ${\cal
H}_{ad}(M,\omega;\alpha_0,\alpha_\infty)$ is nonempty under the
assumptions there. Without loss of generality let $\alpha_0$ be
represented by a compact connected submanifold $S\subset \Int M$
without boundary. Since $\dim\alpha_0+\dim\alpha_\infty\le\dim
M-1$ it follows from intersection theory that there is a cycle
representative $\tilde\alpha_\infty$ of $\alpha_\infty$ such that
$S \cap \tilde\alpha_\infty =\emptyset$.

Choose a Riemannian metric $g$ on $M$. For $\epsilon >0$ let
${\mathcal N}_\epsilon$ be the closed $\epsilon$-ball bundle in
the normal bundle along $S$, and let $\exp \colon {\mathcal
N}_\epsilon \to M$ be the exponential map. For $\epsilon >0$ small
enough, $P = S_\epsilon = \exp ({\mathcal N}_\epsilon)$ and
$Q=S_{2 \epsilon} = \exp ({\mathcal N}_{2\epsilon})$ are smooth
compact submanifolds of $M$ of codimension zero,  and $S_{2
\epsilon}$ is still disjoint from $\tilde\alpha_\infty$. Since
$\dim S=\dim\alpha_0\le \dim M-2$, both $P$ and $Q$ have connected
boundary.

Take a smooth function $f \colon R\to\R$ such that $f(t)=0$ for
$t\le\varepsilon^2$, $f(t)=1$ for $t\ge 4\varepsilon^2$ and
$f^\prime(t)>0$ for $\varepsilon^2<t<4\varepsilon^2$. We define a
smooth function $F \colon M\to\R$ by $F(x)=0$ for $x\in P$,
$F(x)=1$ for $x\in M\setminus Q$ and $F(x)=f(\|v_x\|^2_g)$ for
$x=(s_x, v_x)\in S_{2\varepsilon}$. In view of Lemma~\ref{lem:2.2}
above, for $\delta>0$ sufficiently small the function
$F_\delta=\delta F$ belongs to ${\cal
H}_{ad}(M,\omega;\alpha_0,\alpha_\infty)$.
\hfill$\Box$\vspace{2mm}

\noindent{\bf Proof of Proposition~\ref{prop:1.7}}. Note that
every function $H$ in ${\cal H}_{ad}(W, \omega;\tilde\alpha_0,
pt)$ can be viewed as one in ${\cal H}_{ad}(M, \omega; \alpha_0,
\alpha_\infty)$ in a natural way, and so (\ref{e:1.6}) follows.

If the inclusion $W\hookrightarrow M$ induces an injective
homomorphism $\pi_1(W)\to\pi_1(M)$ then each function $H$ in
${\cal H}_{ad}^\circ(W, \omega;\tilde\alpha_0, pt)$ can be viewed
as one in ${\cal H}_{ad}^\circ(M, \omega; \alpha_0,
\alpha_\infty)$. Therefore we get (\ref{e:1.8}).

To prove (\ref{e:1.10}) let us take a function $H\in {\cal
H}_{ad}(M\setminus W, \omega; \tilde\alpha_\infty, pt)$. Suppose
that $P(H)\subset Q(H)\subset \Int (M\setminus W)$ are
submanifolds associated with $H$. Then $H=\max H$ on $(M\setminus
W)\setminus Q$. Therefore we can extend $H$ to $M$ by setting
$H=\max H$ on $W$. We denote this extension by $\bar H$. Since  we
have assumed that $\alpha_0$  has a cycle representative whose
support is contained in $\Int (W)\subset M\setminus Q$, $\bar H$
belongs to ${\cal H}_{ad}(M,\omega;\alpha_\infty,\alpha_0)$.

If $H\in {\cal H}_{ad}^\circ(M\setminus W, \omega;
\tilde\alpha_\infty, pt)$ and the inclusion $M\setminus
W\hookrightarrow M$ induces an injective homomorphism
$\pi_1(M\setminus W)\to\pi_1(M)$ then the above $\bar H$ belongs
to ${\cal H}_{ad}^\circ (M,\omega;\alpha_\infty,\alpha_0)$. This
implies (\ref{e:1.11}).

For (\ref{e:1.12}) we only need to prove that ${\cal
W}_G(M,\omega) \le C_{HZ}^{(2)}(M,\omega; pt, \alpha)$ since
$C_{HZ}^{(2)}(M,\omega; pt, \alpha)\le C_{HZ}^{(2\circ)}(M,\omega;
pt, \alpha)$. For any given symplectic embedding $\psi
\colon(B^{2n}(r),\omega_0)\to (\Int (M),\omega)$ and sufficiently
small $\epsilon>0$ we can choose a representative of $\alpha$ with
support in $M\setminus\psi(B^{2n}(r-\epsilon))$ because $\dim
\alpha \le \dim M -1$. By (\ref{e:1.5}) and (\ref{e:1.7}) we have
$$
\pi (r-\epsilon)^2 ={\cal W}_G(\psi(B^{2n}(r-\epsilon)),\omega)
\le C_{HZ}(\psi(B^{2n}(r-\epsilon)),\omega)\le
C_{HZ}^{(2)}(M,\omega; pt,\alpha).
$$
With $\epsilon\to 0$ we arrive at the desired conclusion.
\hfill$\Box$\vspace{2mm}

\noindent{\bf Proof of Theorem~\ref{t:1.8}.} \hspace{2mm}To prove
(\ref{e:1.13}) let $W$ and $\alpha_0,\alpha_\infty$ satisfy the
assumptions in Theorem~\ref{t:1.8}. For $H\in{\cal
H}_{ad}(W,\omega;\tilde\alpha_0, pt)$ and $G\in {\cal
H}_{ad}(M\setminus W,\omega;\tilde\alpha_\infty, pt)$ let
$P_1\subset \Int (Q_1)\subset Q_1\subset \Int (W)$ and $P_2\subset
\Int (Q_2)\subset Q_2\subset M\setminus W$ be corresponding
submanifolds as in Definition~\ref{def:1.2}.
  Then $H|_{P_1}=0$, $H|_{W\setminus \Int (Q_1)}=\max H$ and $G|_{P_2}=0$,
$G|_{(M\setminus W)\setminus \Int (Q_2)}=\max G$. Define $K \colon
M\to\R$ by
$$
K(x)=\left\{\begin{array}{lll}
H(x), &&{\rm if}\; x\in W,\\
\max H+ \max G- G(x), &&{\rm if}\; x\in M\setminus W.
\end{array}\right.
$$
    This is
a smooth function and belongs to ${\cal
H}_{ad}(M,\omega;\alpha_0,\alpha_\infty)$ with $P(K)=P_1$ and
$Q(K)=M\setminus \Int (P_2)$. But $\max K=\max H+\max G$. This
leads to (\ref{e:1.13}).\hfill$\Box$\vspace{2mm}

The following corollary of Theorem~\ref{t:1.8} will be useful
later on.
\begin{corollary}\label{cor:2.3}
Under the assumptions of Theorem~\ref{t:1.8}, let $(N,\sigma)$ be
another closed connected symplectic manifold and $\beta\in
H_\ast(N; \Q)\setminus\{0\}$. Then
\begin{eqnarray*}
& &C_{HZ}^{(2)}(N\times W, \sigma\oplus\omega;
\beta\times\tilde\alpha_0, pt)+ C_{HZ}^{(2)}(N\times (M\setminus
W), \sigma\oplus\omega;
\beta\times\tilde\alpha_\infty, pt)\\
& &\le C_{HZ}^{(2)}(N\times M,\sigma\oplus\omega;
\beta\times\alpha_0, \beta\times\alpha_\infty),
\end{eqnarray*}
and
\begin{eqnarray*}
& &C_{HZ}^{(2\circ)}(N\times W, \sigma\oplus\omega;
\beta\times\tilde\alpha_0, pt)+ C_{HZ}^{(2\circ)}(N\times
(M\setminus W), \sigma\oplus\omega;
\beta\times\tilde\alpha_\infty, pt)\\
& &\le C_{HZ}^{(2\circ)}(N\times M,\sigma\oplus\omega;
\beta\times\alpha_0, \beta\times\alpha_\infty)
\end{eqnarray*}
if both inclusions $W\hookrightarrow M$ and $M\setminus
W\hookrightarrow M$ also induce an injective homomorphisms
$\pi_1(W)\to\pi_1(M)$ and $\pi_1(M\setminus W)\to\pi_1(M)$.
\end{corollary}

\section{The proof of Theorem~\ref{t:1.10}}

We wish to reduce the proof of this theorem to the arguments in
\cite{LiuT}. Liu-Tian's approach is to introduce the Morse
theoretical version of Gromov-Witten invariants. In their work the
paper \cite{FHS} plays an important role. To show how the
arguments in \cite{LiuT} apply to our case we need to recall some
related material from \cite{FHS}.

  Consider the vector space ${\cal S}=\{S\in\R^{2n\times 2n} \mid
S^T=S\}$ of symmetric $(2n\times 2n)$-matrices. It has an
important subset ${\cal S}^{2n}_{\rm reg}$ consisting of all
matrices $S\in{\cal S}$ such that for any four real numbers $a,
b,\alpha,\beta$ the system of equations
\begin{equation}\label{e:3.1}
\left\{\begin{array}{ll}
(SJ_0-J_0S-aI_{2n}-bJ_0)\zeta=0\\
(SJ_0-J_0S-aI_{2n}-bJ_0)S\zeta-\alpha\zeta-\beta J_0\zeta=0
\end{array}\right.
\end{equation}
has no nonzero solution $\zeta\in\R^{2n\times 2n}$, where $I_n$
denotes the identity matrix in $\R^{n\times n}$ and
$$
J_0=\left(\begin{array}{cr} 0& -I_n\\I_n& 0\end{array}\right).
$$
It
has been proved in Theorem 6.1 of \cite{FHS} that for $n\ge 2$ the
set ${\cal S}^{2n}_{\rm reg}$ is open and dense in ${\cal S}$ and
$\tau\Phi^TS\Phi\in{\cal S}^{2n}_{\rm reg}$ for any $S\in{\cal
S}^{2n}_{\rm reg}$, any $\Phi\in GL(n,\C)\cap O(2n)$ and any real
number $\tau\ne 0$. In view of Definition~7.1 in \cite{FHS} and the
arguments in \cite{McSl} we introduce

\begin{definition}  \label{def:3.1}
{\rm A nondegenerate critical point $p$ of a smooth function $H$
on a symplectic manifold $(M,\omega)$ is called  {\bf strong
admissible} if it satisfies the following two conditions:
\begin{description}
\item[(i)] the spectrum of the linear transformation
            $DX_H(p): T_pM\to T_pM$ is contained in
            $\C\setminus\{\lambda i\,|\, 2\pi\le\pm\lambda<+\infty\}$;
\item[(ii)] there exists $J_p\in {\cal J}(T_pM,\omega_p)$
             such that for some (and hence every) unitary frame
            $\Phi \colon \R^{2n}\to T_pM$ (i.e.\ $\Phi J_0=J_p\Phi$ and
            $\Phi^\ast\omega_p=\omega_0$) we have
           $$S=J_0\Phi^{-1}DX_H(p)\Phi\in{\cal S}^{2n}_{\rm reg}.$$
\end{description}
}
\end{definition}

\begin{definition}  \label{def:3.2}
{\rm An $(\alpha_0,\alpha_\infty)$-admissible  $($resp.\
$(\alpha_0,\alpha_\infty)^\circ$-admissible) function $H$ in
Definition~\ref{def:1.2} is said to be {\bf
$(\alpha_0,\alpha_\infty)$-strong admissible} (resp.\ {\bf
$(\alpha_0,\alpha_\infty)^\circ$-strong admissible}) if instead of
condition {\bf (5)} it satisfies the stronger condition
\begin{description}
\item[{\bf ($5'$)}] $H$ has only finitely many critical points in
$\Int (Q) \setminus P$,  and each of them is strong admissible in
the sense of Definition~\ref{def:3.1}.
\end{description}
}
\end{definition}

Let us respectively denote by
\begin{equation}\label{e:3.2}
{\cal H}_{sad}(M,\omega;\alpha_0,\alpha_\infty)\quad{\rm and}\quad
{\cal H}_{sad}^\circ(M,\omega;\alpha_0,\alpha_\infty)
\end{equation}
the set of $(\alpha_0,\alpha_\infty)$-strong admissible and
$(\alpha_0,\alpha_\infty)^\circ$-strong admissible functions. They
are subsets of $ {\cal H}_{ad}(M,\omega;\alpha_0,\alpha_\infty)$
and ${\cal H}_{ad}^\circ(M,\omega;\alpha_0,\alpha_\infty)$
respectively. The following lemma is key to our proof.

\begin{lemma}  \label{lem:3.3}
If $\dim M\ge 4$, then ${\cal
H}_{sad}(M,\omega;\alpha_0,\alpha_\infty)$ {\rm (}resp.\ ${\cal
H}^\circ_{sad}(M,\omega;\alpha_0,\alpha_\infty)${\rm )} is
$C^0$-dense in ${\cal H}_{ad}(M,\omega;\alpha_0,\alpha_\infty)$
{\rm (}resp.\ ${\cal
H}^\circ_{ad}(M,\omega;\alpha_0,\alpha_\infty)${\rm )}.
\end{lemma}

\noindent{\bf Proof.}\quad Let $F\in {\cal
H}_{ad}(M,\omega;\alpha_0,\alpha_\infty)$ (resp.\ ${\cal
H}_{ad}^\circ(M,\omega;\alpha_0,\alpha_\infty)$). We shall prove
that for any small $\epsilon>0$ there exists a $G\in {\cal
H}_{sad}(M,\omega;\alpha_0,\alpha_\infty)$ (resp.\ ${\cal
H}_{sad}^\circ(M,\omega;\alpha_0,\alpha_\infty)$) such that
\begin{equation}\label{e:3.3}
\max F\ge\max G\ge\max F-\epsilon.
\end{equation}
Our proof is inspired by  the proof of Proposition 3.1 in
\cite{Schl}.

  Let $C_F$ (resp.\ $c_F$) be
the largest (resp.\ smallest) critical value of $F$ in $(0,\max F)$.
If there are no such critical values, there is nothing to show. If
$c_F=C_F$, then it is the only critical value of $F$ in $(0, \max
F)$, and this case can easily be proved by the following method. So
we now assume $c_F < C_F$. Then by Definition~\ref{def:1.2}(5) we
have
$$0<c_F<C_F<\max F.$$
Let $C(F)$ be the set of critical values of $F$. It is compact and
has zero Lebesgue measure, so that for small $\epsilon>0$ we can
choose regular values of $F$,
$$
b_0'<a_1'<b_1'<\cdots<a'_{k-1}<b'_{k-1}<a'_k,
$$
such that: \\
(i) $0<b'_0<c_F$ and $C_F<a'_k<\max F$.\\
(ii) $[a'_i, b'_i]\subset [c_F, C_F]\setminus C(F)$,
$i=1,\cdots,k-1$.\\
(iii) $\sum^{k-1}_{i=1}(b'_i-a'_i)+ b'_0+ \max F-a'_k>\max
F-\epsilon$.\\
Furthermore we may also take regular values of $F$,
$$
b_0<a_1<b_1<\cdots<a_{k-1}<b_{k-1}<a_k,
$$
such that
\begin{eqnarray*}
&&b_0<b'_0,\;a_k>a'_k, \;a'_i<a_i<b_i<b'_i,\;i=1,\cdots,k-1,\\
&&\sum^{k-1}_{i=1}(b_i-a_i)+ b_0+ \max F-a_k>\max F-2\epsilon.
\end{eqnarray*}
Consider the piecewise-linear function $f \colon \R\to\R$,
$$
f(t)=\left\{\begin{array}{ll}
  t\quad &{\rm for}\;t\le b_0,\\
b_0\quad &{\rm for}\; b_0\le t\le a_1,\\
t-a_1+ b_0\quad &{\rm for}\; a_1\le t\le b_1,\\
b_1-a_1+ b_0\quad &{\rm for}\; b_1\le t\le a_2,\\
t-a_2+ (b_1-a_1)+ b_0\quad &{\rm for}\; a_2\le t\le b_2,\\
\cdots\quad &{\rm for}\;\cdots,\\
t-a_{k-1}+\sum^{k-2}_{i=1}(b_i-a_i)+ b_0\quad &{\rm for}\;
a_{k-1}\le
t\le b_{k-1},\\
\sum^{k-1}_{i=1}(b_i-a_i)+ b_0\quad &{\rm for}\; b_{k-1}\le t\le a_k,\\
t-a_k+ \sum^{k-1}_{i=1}(b_i-a_i)+ b_0\quad &{\rm for}\;  t\ge a_k.
\end{array}\right.$$
Then $\min\{f(t)\,|\, t\in [0,\max F]\}=0$ and
$$
\max\{f(t)\,|\, t\in [0,\max F]\}=\max F-a_k+
\sum^{k-1}_{i=1}(b_i-a_i)+ b_0>\max F-2\epsilon.
  $$
  Note that
$b_0<a_1<b_1<\cdots<a_{k-1}<b_{k-1}<a_k$ are all nonsmooth points of
$f$ in $(0,\max F)$. By suitably smoothing $f$ near these points we
can get a smooth function $h \colon \R\to\R$ satisfying:
\begin{description}
\item[$(h)_1$] $0\le h^\prime(t)\le1$ for $t \in \R$;

\item[$(h)_2$] $0<h^\prime(t)\le 1$ for $t\in [0, b'_0)\cup (a'_k,
\max F]\cup(\cup^{k-1}_{i=1}(a'_i, b'_i))$;

  \item[$(h)_3$] $h(t)=f(t)$ for $t\in\cup^{k-1}_{i=0}[b'_i, a'_{i+1}]$;

\item[$(h)_4$] $h(t)=f(t)$ near $t=0$ and $t=\max F$.
\end{description}
  Set
$H=h\circ F$. Then $(h)_1$ and $(h)_4$ imply that $H\in{\cal
H}_{ad}(M,\omega;\alpha_0,\alpha_\infty)$ and
\begin{equation}\label{e:3.4}
\quad\max H=h(\max F)=\max F-a_k+ \sum^{k-1}_{i=1}(b_i-a_i)+
b_0>\max F-2\epsilon.
\end{equation}
Furthermore, one easily checks that
\begin{description}
\item[$(H)_1$] The critical values of $H$ in $(0,\max H)$ are
exactly $b_0$, $\sum^{j}_{i=1}(b_i-a_i)+b_0$, $j=1,\cdots, k-1$;

\item[$(H)_2$] The corresponding critical sets are respectively
$\{b'_0\le F\le a'_1\}$ and $\{b'_j\le F\le a'_{j+1}\}$,
$j=1,\cdots, k-1$;

  \item[$(H)_3$] $H=b_0$ on $\{b'_0\le F\le a'_1\}$;

\item[$(H)_4$] $H=\sum^{j}_{i=1}(b_i-a_i)+b_0$ on $\{b'_j\le F\le
a'_{j+1}\}$, $j=1,\cdots, k-1$.
\end{description}
For each
  $0\le s< \frac{1}{2}\min\{a'_{i+1}-b'_i,  b'_i-a'_i, b'_0, \max
F-a'_k\,|\, 0\le i\le k-1\}$  we set
$$N_s:=\bigcup^{k-1}_{i=0}\{b'_i-s\le F\le a'_{i+1}+ s\}.$$
Since the set of regular values of $F$ is open, both $N_0$ and
$N_s$ with sufficiently small $s>0$ are compact smooth
submanifolds with boundary. For any open neighborhood ${\cal O}$
of $N_0$ we have also $N_s\subset {\cal O}$ if  $s>0$ is small
enough.
  By
$(H)_3$ and $(H)_4$, $\nabla_g\nabla_g H=0$ on $N_0$ and thus we
can choose
$$0<\delta<\frac{1}{4}\min\{a'_{i+1}-b'_i,  b'_i-a'_i, b'_0, \max
F-a'_k\,|\, 0\le i\le k-1\}$$ so  small that
$$
\sup_{x\in N_{2\delta}}\|\nabla_g\nabla_g H(x)\|_g<\rho/2.
$$
Here $\rho$ is given by Lemma~\ref{lem:2.2}. Let us take a smooth
function $L \colon M \to \R$ such that
\begin{description}
\item[$(L)_1$] ${\rm supp}(L)\subset  N_\delta$;

\item[$(L)_2$] $\|L\|_{C^2} < \rho/2$ (and thus
            $\sup_{x\in N_{2\delta}}\|\nabla_g\nabla_g (H+L)(x)\|_g<\rho$);

\item[$(L)_3$] $h(b'_i-2\delta)< H(x)+L(x)<h(a'_{i+1}+2\delta)$
for $x\in \{b'_i-\delta\le F\le a'_{i+1}+\delta\}$, $i=0,\cdots,
k-1$;

\item[$(L)_4$]  $H+L$ has only finitely many critical points in
$N_\delta$ and each of them is strong admissible.
\end{description}
  The condition $(L)_4$ can be assured by Lemma 7.2 (i) in [FHS].
  To see that $(L)_3$ can be satisfied,  note that $(h)_1$ implies that
  $h(b'_i-\delta)\le H(x)\le
h(a'_{i+1}+\delta)$ as $b'_i-\delta\le F(x)\le a'_{i+1}+\delta$.
  By the choice of $\delta$ we have $b'_0-2\delta>0$,
$a'_k+2\delta<\max
  F$ and
  $$a'_i+\delta,\, a'_i+2\delta,\,
  \,b'_i-2\delta, \,b'_i-\delta\in (a'_i, b'_i),\; i=1,\cdots, k-1.$$
It follows from $(h)_2$ that for $i=0,\cdots, k-1$,
  \begin{equation}\label{e:3.5}
  \qquad h(b'_i-2\delta)<h(b'_i-\delta)<h(b'_i)\le
  h(a'_{i+1})<h(a'_{i+1}+\delta)<h(a'_{i+1}+2\delta).
  \end{equation}
Using these and $(L)_2$ we can easily choose $L$ satisfying
$(L)_3$. Set $G=H+ L$. Then
  $P(G)=P(H)$, $Q(G)=Q(H)$ and
\begin{equation}\label{e:3.6}
  \max G=\max H\quad{\rm and}\quad \min G=\min H=0.
  \end{equation}

Now we are in position to prove
$$
G\in {\cal H}_{sad}(M,\omega;\alpha_0,\alpha_\infty)\quad{\rm
(resp.}\; {\cal
H}^{\circ}_{sad}(M,\omega;\alpha_0,\alpha_\infty){\rm )}.$$

   Firstly, the above construction shows that all critical
values of $G$ in $(0,\max G)$ sit in
$$
\bigcup^{k-1}_{i=0}\bigr(h(b'_i-2\delta),
h(a'_{i+1}+2\delta)\bigl)
$$
and the corresponding critical points
sit in $N_\delta$. It follows that $G$ has only finitely many
critical points in ${\rm Int}(Q)\setminus P$ and each of them is
strong admissible.

Next we prove that $X_G$ has no nonconstant fast periodic orbits.
Assume that $\gamma$ is such an orbit. It cannot completely sit in
$M\setminus N_\delta$ because $G=H$ in $M\setminus N_\delta$.
Moreover, Lemma~\ref{lem:2.2} and $(L)_2$ imply that $\gamma$
cannot completely sit in $N_{2\delta}$. So there must exist two
points $\gamma(t_1)$ and $\gamma(t_2)$ such that $\gamma(t_1)\in
\partial N_\delta$ and $\gamma(t_2)\in
\partial N_{2\delta}$. Note that all possible values $G$ takes on
$\partial N_\delta$ (resp.\ $\partial N_{2\delta}$) are
  \begin{eqnarray*}
  h(b'_i-\delta), h(a'_{i+1}+\delta), i=0,\cdots, k-1.\\
  ({\rm resp.}\;h(b'_i-2\delta), h(a'_{i+1}+2\delta), i=0,\cdots,
  k-1.)
  \end{eqnarray*}
By (\ref{e:3.5}) any two of them are different. But
$G(\gamma(t_1))=G(\gamma(t_2))$. This contradiction shows that
$X_G$ has no nonconstant fast periodic orbit. Clearly, this
argument also implies that $X_G$ has no nonconstant contractible
fast periodic orbit if $F\in {\cal
H}^{\circ}_{ad}(M,\omega;\alpha_0,\alpha_\infty)$.

  Finally,  (\ref{e:3.4}) and
(\ref{e:3.6}) together gives
$$\max G\ge\max F-2\epsilon.$$
The desired conclusion is proved. \hfill$\Box$\vspace{2mm}

As direct consequences of Lemma~\ref{lem:3.3} and (\ref{e:1.2}) we
have
\begin{equation} \label{e:3.7}
 \qquad\left\{\begin{array}{ll}
C_{HZ}^{(2)}(M,\omega;\alpha_0,\alpha_\infty)=\sup \left\{\max
H\,|
\, H\in {\cal H}_{sad}(M,\omega;\alpha_0,\alpha_\infty) \right\},\\
C_{HZ}^{(2\circ)}(M,\omega;\alpha_0,\alpha_\infty)=\sup
\left\{\max H\,| \, H\in {\cal
H}_{sad}^\circ(M,\omega;\alpha_0,\alpha_\infty) \right\}.
\end{array}
\right.
\end{equation}
\vspace{2mm}

\noindent{\bf Proof of Theorem~\ref{t:1.10}}. We only prove
(\ref{e:1.15}). The proof of (\ref{e:1.16}) is similar. Without
loss of generality we assume that
$C_{HZ}^{(2)}(M,\omega;\alpha_0,\alpha_\infty)>0$ and ${\rm
GW}(M,\omega;\alpha_0,\alpha_\infty)<+\infty$. We need to prove
that if
\begin{equation}\label{e:3.8}
  \Psi_{A, g,
m+2}(C;\alpha_0,\alpha_\infty,\beta_1,\dots,\beta_m)\ne 0
\end{equation}
for homology classes $A\in H_2(M;\Z)$, $C\in H_\ast(\overline{\cal
M}_{g, m+2};\Q)$ and $\beta_1,\dots,\beta_m\in H_\ast(M;\Q)$ and
integers $m \ge 1$ and $g\ge 0$, then
\begin{equation}\label{e:3.9}
C_{HZ}^{(2)}(M,\omega;\alpha_0,\alpha_\infty)\le\omega(A).
\end{equation}

  Arguing by contradiction, we may assume by (\ref{e:3.7})
that there exists $H\in{\cal
H}_{sad}(M,\omega;\alpha_0,\alpha_\infty)$ such that $\max
H>\omega(A)$. Then we take $\eta>0$ such that
\begin{equation}\label{e:3.10}
\max H-2\eta>\omega(A).
\end{equation}
By the properties of $H$ there exist two smooth compact
submanifolds $P, Q\subset M$ with connected boundary and of
codimension zero such that the conditions (1),(2),(3),(4),(6) in
Definition~\ref{def:1.2} and ($5^\prime$) in
Definition~\ref{def:3.2} are satisfied.
Changing $H$ slightly near $\{ H =0 \}$ and near $\{ H = \max H \}$
in the class ${\mathcal H}_{sad}(M,\omega;\alpha_0,\alpha_\infty)$
and using Lemma~\ref{lem:2.1}, we can choose embeddings
$$
\Phi \colon [-2, 0]\times\partial Q\to Q\setminus \Int (P)
\quad{\rm and}\quad \Psi \colon [0, 2]\times\partial P\to
Q\setminus \Int (P) $$
  such that:
\begin{description}
\item[(i)] $\Phi(\{0\}\times\partial Q)=\partial Q$ and
           $\Psi(\{0\}\times\partial P)=\partial P$;
\item[(ii)] $\Phi([-2, 0]\times\partial Q)\cap\Psi([0,
             2]\times\partial P)=\emptyset$;
\item[(iii)] $H$ has no critical points in
             $\Phi([-2, 0)\times\partial Q)\cup\Psi((0,
             2]\times\partial P)$ and is constant $m_s$ on
             $\Phi(\{s\}\times\partial Q)$ and $n_t$ on
              $\Psi(\{t\}\times\partial P)$ for each $s\in [-2, 0]$
             and $t\in [0,2]$;
\item[(iv)] $H(x)<m_s$ for any $x\in M\setminus\widehat Q_s$ and
               $s\in [-2,0]$, and $n_t<H(x)$ for any $x\in M\setminus
\widehat
               P_t$ and $t\in [0,2]$, where
$$
\widehat Q_s=(M\setminus Q)\cup\Phi([s, 0]\times\partial Q)
\quad{\rm and}\quad \widehat P_t=P\cup\Psi([0, t]\times\partial
P).
$$
\end{description}
  Notice that the above assumptions imply
$$
m_s<m_{s'}<\max H\quad{\rm and}\quad 0<n_t<n_{t'}
  $$
  for $-2<s<s'<0$ and $0<t<t'<2$. Moreover $\widehat Q_s$ (resp.
$\widehat P_t$) is a smooth compact submanifold of $M$ with
boundary $\Phi(\{s\}\times\partial Q)$ (resp.\
$\Psi(\{t\}\times\partial P)$). Clearly, $\widehat Q_s\cap\widehat
P_t=\emptyset$. For $\tau\in [0, 2]$ we abbreviate
$$ B_\tau \,=\,\widehat P_\tau\cup\widehat Q_{-\tau}.$$
  By the properties of $H$ and (\ref{e:3.10}) we find
$\delta\in (0, 1)$ such that
\begin{equation}\label{e:3.11}
m_{-2\delta}>\max H-\eta,\quad n_{2\delta}<\eta\quad{\rm and}\quad
\sup_{x\in B_{2\delta}}\|\nabla_g\nabla_g H(x)\|_g<\rho/2,
  \end{equation}
where $\rho$ is as in  Lemma~\ref{lem:2.2}. As before we may
choose a smooth function $L \colon M \to \R$ such that
\begin{description}
\item[(a)] ${\rm supp}(L)\subset \Int (B_\delta)$;

\item[(b)] $\|L\|_{C^2} <\min\{\rho/2,\eta\}$ (and thus
$\sup_{x\in B_{2\delta}}\|\nabla_g\nabla_g (H+L)(x)\|_g< \rho$);

\item[(c)] $H+L$ has only finitely many critical points in
           $\Int (B_\delta)$, and each of them is also strong admissible;

\item[(d)] $m_{-2\delta}<H(x)+L(x)$ for
            $x\in \Int (\widehat Q_{-\delta})$.

\item[(e)] $H(x)+L(x)< n_{2\delta}$ for
            $x\in \Int (\widehat P_\delta)$.
\end{description}
As above, condition (c) is assured by Lemma 7.2 (i) in \cite{FHS}.
Set $F=H+L$. If $x\in B_\delta$ then either $F(x)>m_{-2\delta}$ or
$F(x)<n_{2\delta}$. On the other hand, the above (a) and (iv) imply
that $n_{2\delta}<F(x)<m_{-2\delta}$ if $x\in M\setminus
B_{2\delta}$. This means that a solution of
  $\dot x=X_F(x)$ cannot go to $B_\delta$ from $M\setminus
  B_{2\delta}$ because $F$ is constant along any solution of $\dot
x=X_F(x)$.
  So any nonconstant
solution of $\dot x=X_F(x)$ lies either in $B_{2\delta}$ or in
$M\setminus B_{\delta}$. It follows from (a) and (b) that $\dot
x=X_F(x)$ has no nonconstant fast periodic solutions. Using
(\ref{e:3.11}) and (a)-(e) again we get that $F$ is a smooth Morse
function on $M$ satisfying
\begin{description}
\item[$(F)_1$] each critical point of $F$ is strong admissible;
\item[$(F)_2$] $\lambda\cdot F$ has no non-trivial periodic
                 solution of period $1$ for any $\lambda\in (0, 1]$;
\item[$(F)_3$] $F(x)>\max H-\eta$ for $x\in\widehat Q_{-\delta}$,
and $F(x)<\eta$ for any $x\in\widehat P_\delta$;

\item[$(F)_4$] $\max F\le\max H+\eta$ and $\min F\ge -\eta$.
\end{description}

As a consequence of $(F)_1$ we get that ${\cal J}_{ad}(M,\omega,
X_F)$ is nonempty. From Lemma 7.2(iii) in \cite{FHS} we also know
that ${\cal J}_{ad}(M,\omega, X_F)$ is open in ${\cal J}(M,\omega)$
with respect to the $C^0$-topology.  Therefore, we may choose a
regular $J\in {\cal J}_{ad}(M,\omega, X_F)$ and then repeat the
arguments in \cite{LiuT} to define the Morse theoretical
Gromov-Witten invariants
$$
\Psi_{A, J_\lambda,\lambda F, g, m+2}(C;\alpha_0,\alpha_\infty,
\beta_1,\dots,\beta_m)
$$
and to prove
\begin{equation}  \label{e:3.12}
\quad\Psi_{A, J_\lambda,\lambda F, g,
m+2}(C;\alpha_0,\alpha_\infty, \beta_1,\dots,\beta_m)\equiv
\Psi_{A, g, m+2}(C;\alpha_0,\alpha_\infty, \beta_1,\dots,\beta_m)
\end{equation}
for each $\lambda\in [0, 1]$. As in Lemma 7.2 of \cite{LiuT} we
can prove the corresponding moduli space ${\cal FM}(c_0, c_\infty;
J_1, F, A)$ to be empty for any critical points $c_0\in\widehat
P_\delta$ and $c_\infty\in\widehat Q_{-\delta}$ of $F$.
  In fact, otherwise we may
choose an element $f$ in it. Then one easily gets the estimate
\begin{equation}  \label{e:3.13}
0\le E(f)=F(c_0)-F(c_\infty)+\omega(A).
\end{equation}
  (Note: from the proof of
Lemma 7.2 in \cite{LiuT} one may easily see that the energy
identity above their Lemma 3.2 should read $E(f)=\omega(A)+
H(c_-)-H(c_+)$.) From the above  $(F)_3$ and (\ref{e:3.13}) it
follows that
$$\max H-2\eta<F(c_\infty)-F(c_0)\le\omega(A).$$
  This contradicts (\ref{e:3.10}).
  So ${\cal FM}(c_0, c_\infty; J_1, F, A)$ is empty and thus
$$\Psi_{A, J_1, F, g, m+2}(C;\alpha_0,\alpha_\infty,
\beta_1,\dots,\beta_m)=0.$$ By (\ref{e:3.12}) we get $\Psi_{A, g,
m+2}(C;\alpha_0,\alpha_\infty, \beta_1,\dots,\beta_m)=0$. This
contradicts (\ref{e:3.8}). (\ref{e:3.9}) is proved.
\hfill$\Box$\vspace{2mm}

\section{Proofs of Theorems~\ref{t:1.15},
\ref{t:1.16}, \ref{t:1.17} and \ref{t:1.21}}

\noindent{\bf Proof of Theorem~\ref{t:1.15}}. We start with the
matrix definition of the Grassmannian manifold $G(k, n)=G(k, n;
\C)$. Let $n=k+m$, $M(k, n;\C)=\{ A\in\C^{k\times n}\,|\, {\rm
rank}A=k\,\}$ and ${\rm GL}(k;\C)=\{Q\in\C^{k\times k}\,|\,{\rm
det}Q\ne 0\}$. Then ${\rm GL}(k;\C)$ acts freely on $M(k, n;\C)$
from the left by  matrix multiplication. The quotient $M(k,
n;\C)/{\rm GL}(k;\C)$ is exactly $G(k, n)$. For $A\in M(k, n;\C)$
we denote by $[A]\in G(k, n)$ the ${\rm GL}(k;\C)$-orbit of $A$ in
$M(k, n;\C)$, and by
$$
{\rm Pr} \colon M(k, n;\C)\to G(k, n),\; A\mapsto [A]
$$
the quotient projection. Any representative matrix $B$ of $[A]$ is
called a {\bf homogeneous coordinate} of the point $[A]$. For
increasing integers $1\le\alpha_1<\cdots<\alpha_k\le n$ let
$\{\alpha_{k+1},\cdots, \alpha_n\}$ be the complement of
$\{\alpha_1,\cdots,\alpha_k\}$ in the set $\{1,2,\dots, n\}$. Let
us write $A\in M(k, n;\C)$ as $A=(A_1,\cdots, A_n)$ and
$$A_{\alpha_1\cdots\alpha_k}=(A_{\alpha_1},\cdots,
A_{\alpha_k})\in\C^{k\times k} \;{\rm and}\;
A_{\alpha_{k+1}\cdots\alpha_n}=(A_{\alpha_{k+1}},\cdots,
A_{\alpha_n}) \in\C^{k\times m},$$ where $A_1,\cdots, A_n$ are
$k\times 1$ matrices. Define a subset of $M(k, n;\C)$ by
$$V(\alpha_1,\cdots,\alpha_k)=\{A\in M(k, n;\C)\,|\,
{\rm det}A_{\alpha_1\cdots\alpha_k}\ne 0\,\}$$ and set
$U(\alpha_1,\cdots,\alpha_k)={\rm
Pr}(V(\alpha_1,\cdots,\alpha_k))$ and
\begin{eqnarray*}
\Theta(\alpha_1,\cdots,\alpha_k): U(\alpha_1,\cdots, \alpha_k)\to
\C^{k\times m}\equiv\C^{km},\\
 \quad [A]\to
Z=(A_{\alpha_1\cdots\alpha_k})^{-1}A_{\alpha_{k+1}\cdots\alpha_n}.
\end{eqnarray*}
It is easily checked that this is a homeomorphism.  $Z$ is called
the {\bf local coordinate} of $[A]\in G(k, n)$ in the canonical
coordinate neighborhood $U(\alpha_1,\cdots,\alpha_k)$. Note that
for any $Z\in\C^{k\times m}$ there must exist an $n\times n$
permutation matrix $P(\alpha_1,\cdots,\alpha_k)$ such that for the
matrix $A=(I^{(k)}, Z)P(\alpha_1,\cdots,\alpha_k)$ we have
\begin{equation}  \label{e:4.00}
A_{\alpha_1\cdots\alpha_k}=I^{(k)}\quad{\rm and}\quad
A_{\alpha_{k+1}\cdots\alpha_n}=Z.
\end{equation}
Hereafter $I^{(k)}$ denotes the unit $k\times k$ matrix. It
follows from this fact that for another set of increasing integers
$1\le\beta_1<\cdots <\beta_k\le n$ the transition function
$\Theta(\beta_1,\cdots,
\beta_k)\circ\Theta(\alpha_1,\cdots,\alpha_k)^{-1}$ from
$\Theta(\alpha_1,\cdots,\alpha_k)(U(\alpha_1,\cdots,\alpha_k))$ to
$\Theta(\beta_1,\cdots,\beta_k)(U(\beta_1,\cdots,\beta_k))$ is
given by
$$
Z\to
W=(W_{\beta_1\cdots\beta_k})^{-1}W_{\beta_{k+1}\cdots\beta_n},
$$
  where
$(W_{\beta_1\cdots\beta_k},
W_{\beta_{k+1}\cdots\beta_n})=(I,Z)P(\alpha_1,\cdots,\alpha_k)P^\prime(\beta_1,\cdots,\beta_k)$.
It is not hard to check that this transformation is biholomorphic.
Thus
\begin{equation}  \label{e:4.0}
  \Bigl\{\Bigl(U(\alpha_1,\cdots,\alpha_k),\,
\Theta(\alpha_1,\cdots,\alpha_k)\Bigr)\, \Bigm|\,
1\le\alpha_1<\cdots<\alpha_k\le n\Bigr\}
\end{equation}
gives an atlas of the natural complex structure on $G(k, n)$,
which is called the {\bf canonical atlas}. It is not hard to prove
that the canonical K\"ahler form $\sigma^{(k, n)}$  on $G(k, n)$
in such coordinate charts is given by
$$
\frac{\sqrt{-1}}{2}{\rm tr}[(I^{(k)}+ Z\overline
Z^\prime)^{-1}dZ\wedge (I^{(m)}+ \overline Z^\prime
Z)^{-1}d\overline Z^\prime]=\frac{\sqrt{-1}}{2}
\partial\bar\partial\log\det(I^{(k)}+ Z\overline Z^\prime),
$$
where $dZ=(dz_{ij})_{1\le i\le k, 1\le j\le m}$ and $\partial$,
$\bar\partial$ are the differentials with respect to the
holomorphic and antiholomorphic coordinates respectively (cf.\
\cite{L}).

On the other hand, it is easy to see that
\begin{eqnarray*}
 &&\tau_{k, n}=\frac{\sqrt{-1}}{2}\partial\bar\partial
  \log\det(A\overline A^\prime)\\
&&\qquad =\frac{\sqrt{-1}}{2}{\rm tr}[-(A\overline
A^\prime)^{-1}dA\wedge \overline A^\prime(A\overline A^\prime
)^{-1}Ad\overline A^\prime+ (A\overline A^\prime)^{-1}dA\wedge
d\overline A^\prime]
\end{eqnarray*}
is an invariant K\"ahler form on $M(k, n;\C)$  under the left action
of ${\rm GL}(k;\C)$. Thus it descends to a symplectic form
$\widehat\tau_{k, n}$ on $G(k, n;\C)$. If $A=(I^{(k)}, Z)$ it is
easily checked that
\begin{eqnarray*}
& &\frac{\sqrt{-1}}{2}{\rm tr}[-(A\overline A^\prime)^{-1}dA\wedge
\overline A^\prime(A\overline A^\prime )^{-1}Ad\overline A^\prime+
(A\overline A^\prime)^{-1}dA\wedge d\overline A^\prime]\\
&=&\frac{\sqrt{-1}}{2}{\rm tr}[-(I^{(k)}+ Z\overline
Z^\prime)^{-1}dZ\wedge \overline Z^\prime(I^{(k)}+ Z\overline
Z^\prime )^{-1}Zd\overline Z^\prime \\
&&\qquad +(I^{(k)}+ Z\overline Z^\prime)^{-1}dZ\wedge d\overline Z^\prime]\\
&=&\frac{\sqrt{-1}}{2}{\rm tr}[(I^{(k)}+ Z\overline
Z^\prime)^{-1}dZ\wedge (I^{(m)}+ \overline Z^\prime
Z)^{-1}d\overline Z^\prime].
\end{eqnarray*}
It follows that $\widehat\tau_{k, n}=\sigma^{(k, n)}$. Since
  ${\rm Pr}^\ast\widehat\tau_{k, n}=\tau_{k, n}$ we arrive at
\begin{equation}  \label{e:4.1}
{\rm Pr}^\ast\sigma^{(k, n)}=\tau_{k, n}.
\end{equation}
As usual, if we identify $z=(z_{11},\cdots, z_{1m}, z_{21},\cdots,
z_{2m},\cdots, z_{k1},\cdots, z_{km})\in\C^{km}$ with the matrix
$Z=(z_{ij})_{1\le i\le k, 1\le j\le m}$ the standard symplectic form
in $\C^{km}$ becomes
$$\omega^{(km)}=\frac{\sqrt{-1}}{2}{\rm tr}[dZ\wedge d\overline Z^\prime].
$$
Denote by
$$M^0(k, n;\C)=\{ A\in M(k, n;\C)\,|\, A\overline A^\prime=I^{(k)}\,\}.
$$
Then
\begin{equation}  \label{e:4.2}
\tau_{k, n}|_{M^0(k, n;\C)}=\omega^{(km)}|_{M^0(k, n;\C)}.
\end{equation}
In fact, since $A\overline A^\prime=I^{(k)}$ we have that $dA
\overline A^\prime+ Ad\overline A^\prime=0$ and thus
\begin{eqnarray*}
& &\frac{\sqrt{-1}}{2}{\rm tr}[-(A\overline A^\prime)^{-1}dA\wedge
\overline A^\prime(A\overline A^\prime )^{-1}Ad\overline A^\prime+
(A\overline A^\prime)^{-1}dA\wedge d\overline A^\prime]\\
&=&\frac{\sqrt{-1}}{2}{\rm tr}[dA\wedge\overline dA^\prime] +
\frac{\sqrt{-1}}{2}{\rm tr}[dA\overline A^\prime\wedge dA\overline
A^\prime].
\end{eqnarray*}
We want to prove the second term is zero. A direct computation
yields
\begin{eqnarray*}
{\rm tr}[dA\overline A^\prime\wedge dA\overline A^\prime]
&=&\sum^k_{i=1}\sum^k_{j=1}(\sum^n_{s=1}\bar a_{js} da_{is})\wedge
(\sum^n_{s=1}\bar a_{is} da_{js})\\
&=&\sum^k_{j=1}\sum^k_{i=1}(\sum^n_{s=1}\bar a_{is} da_{js})\wedge
(\sum^n_{s=1}\bar a_{js} da_{is})\;({\rm interchanging}\;i, j)\\
&=&-\sum^k_{i=1}\sum^k_{j=1}(\sum^n_{s=1}\bar a_{js}
da_{is})\wedge (\sum^n_{s=1}\bar a_{is} da_{js}).
\end{eqnarray*}
Hence ${\rm tr}[dA\overline A^\prime\wedge dA\overline
A^\prime]=0$. (\ref{e:4.2}) is proved.

\begin{lemma}  \label{lem:4.1}
For the classical domain of the first type {\rm (}cf.\
\cite{L}{\rm )}
$$
R_I(k, m)=\{ Z\in\C^{k\times m}\,|\, I^{(k)}-Z\overline
Z^\prime>0\},
$$
the map
$$
\Phi \colon (R_I(k, m), \omega^{(km)})\to (\C^{k\times n},
\omega^{(kn)}),\;\; Z\mapsto \bigl(\sqrt{I^{(k)}-Z\overline
Z^\prime}, Z\bigr)
$$
is a symplectic embedding with image in $M^0(k, n;\C)$, and
therefore we get a symplectic embedding $\widehat\Phi={\rm
Pr}\circ\Phi$ of $(R_I(k, m), \omega^{(km)})$ into $(G(k,
n;\C),\sigma^{(k, n)})$.
\end{lemma}

\noindent{\bf Proof}.\hspace{2mm} Differentiating
$$
\Phi(Z)\overline{\Phi(Z)}^\prime=\sqrt{I^{(k)}-Z\overline
Z^\prime} \sqrt{I^{(k)}-Z\overline Z^\prime}+ Z\overline
Z^\prime=I^{(k)}
$$
twice from both sides we get
$$
d\sqrt{I^{(k)}-Z\overline Z^\prime}\bigwedge
d\sqrt{I^{(k)}-Z\overline Z^\prime}=0.
$$
This leads to
\[
d\Phi(Z)\wedge d\overline{\Phi(Z)}^\prime=dZ\wedge d\overline
Z^\prime,\;\;i.e., \Phi^\ast\omega^{(kn)}=\omega^{(km)}.
\]
Using (\ref{e:4.1}) and (\ref{e:4.2}) we get that the composition
$\widehat\Phi={\rm Pr} \circ\Phi$ yields the desired symplectic
embedding from $(R_I(k, m), \omega^{(km)})$ to $(G(k,
n;\C),\sigma^{(k, n)})$. \hfill$\Box$\vspace{2mm}

\begin{lemma}  \label{lem:4.2}
The open unit ball $B^{2km}(1)$ is contained in $R_I(k, m)$.
\end{lemma}

\noindent{\bf Proof}.\hspace{2mm} It is well known that for any
$Z\in\C^{k\times m}$ with $k\le m$ (resp.\ $k>m$) there exist
unitary matrices $U$ of order $k$ and $V$ of order $m$  such that
$$UZV=({\rm diag}(\lambda_1,\cdots,\lambda_k), O)\; (\,{\rm resp.}\;
UZV=({\rm diag}(\mu_1,\cdots,\mu_m), O)^\prime)$$ for some
$\lambda_1\ge\cdots\ge\lambda_k\ge 0$ (resp.\ $\mu_1\ge\cdots\ge
\mu_m\ge 0$), where ${\rm diag}(\lambda_1,\cdots,\lambda_k)$
(resp. ${\rm diag}(\mu_1,\cdots,\mu_m)$) denote the diagonal
matrix of order $k$ (resp.\ $m$), and $O$ is the zero matrix of
order $k\times (m-k)$ (resp.\ $(k-m)\times m$).  Therefore, $Z\in
R_I(k, m)$, i.e., $I^{(k)}-Z\overline Z^\prime>0$, if and only if
$\lambda_j<1$, $j=1,\cdots, k$, (resp.\ $\mu_i<1$, $i=1,\cdots,
m$). Let $Z\in B^{2km}(1)$. Then
$$\|Z\|^2=\sum^k_{i=1}\sum^m_{j=1}|z_{ij}|^2={\rm tr}(Z\overline
Z^\prime)= \sum^m_{j=1} \lambda_j^2\; {\rm (}{\rm resp.\ }
\sum^n_{k=1} \mu_k^2{\rm )}<1,$$ and thus $\lambda_j<1$ (resp.\
$\mu_i<1$), i.e., $Z\in R_I(k, m)$. \hfill$\Box$\vspace{1mm}

Now Lemma~\ref{lem:4.1} and Lemma~\ref{lem:4.2} yield directly
\begin{equation}  \label{e:4.3}
{\cal W}_G(G(k, n), \sigma^{(k,n)})\ge {\cal W}_G(R_I(k, m),
\omega^{(km)})\ge\pi
\end{equation}
  for $m=n-k$. Moreover, for the
submanifolds $X^{(k,n)}$ and $Y^{(k,n)}$ of $G(k, n)$ the
computation in \cite{SieT,Wi} shows $\Psi_{L^{(k,n)},0,3}(pt;
[X^{(k,n)}], [Y^{(k,n)}], pt)=1$. Thus (\ref{e:1.12}) and
Theorem~\ref{t:1.13} lead to
\begin{equation}  \label{e:4.4}
\quad {\cal W}_G(G(k, n), \sigma^{(k,n)})\le C_{HZ}^{(2)}(G(k,
n),\sigma^{(k,n)}; pt, \alpha)\le\sigma^{(k,n)}(L^{(k, n)})=\pi
\end{equation}
for $\alpha=[X^{(k,n)}]$ or $\alpha=[Y^{(k,n)}]$ with $k\le n-2$.
Hence the conclusions follow from (\ref{e:4.3}) and (\ref{e:4.4}).
Theorem~\ref{t:1.15} is proved. \hfill$\Box$\vspace{2mm}

\noindent{\bf Proof of Theorem~\ref{t:1.16}}. \hspace{2mm} Since
$\Psi_{L^{(k,n)},0,3}(pt; [X^{(k,n)}], [Y^{(k,n)}], pt)=1$ it
follows from Proposition~\ref{prop:7.4} that
$$\Psi_{A,0,3}(pt; [M]\times [X^{(k,n)}], [M]\times [Y^{(k,n)}], pt)\ne 0$$
for $A=0\times L^{(k,n)}$, where $0$ denotes the zero class in
$H_2(M;\Z)$. Theorem~\ref{t:1.13} implies
$$
C_{HZ}^{(2\circ)}(M\times G(k,n), \omega\oplus (a\sigma^{(k,n)});
  pt, [M]\times\alpha)\le |a|\pi
$$
for $\alpha=[X^{(k, n)}]$ or $\alpha=[Y^{(k,n)}]$ with $k\le n-2$.
This implies (\ref{e:1.20}).

For (\ref{e:1.21}) we only prove the case $r=2$ for the sake of
simplicity. The general case is similar. Let us take
$A=\oplus^2_{i=1}L^{(k_i, n_i)}\in H_2(W,\Z)$. Then
$\Omega(A)=(|a_1|+|a_2|)\pi$. Note that
\begin{eqnarray*}
&&\Psi_{L^{(k_i,n_i)},0,3}(pt; pt, [X^{(k_i,n_i)}],
[Y^{(k_i,n_i)}])\\
&& =\Psi_{L^{(k_i,n_i)},0,3}(pt; pt, [Y^{(k_i,n_i)}],
[X^{(k_i,n_i)}])=1
\end{eqnarray*}
because the dimensions of $[X^{(k_i,n_i)}]$ and $[Y^{(k_i,n_i)}]$
are even for $i=1,2$. Proposition~\ref{prop:7.7} gives
\begin{eqnarray*}
&&\Psi_{A,0,3}(pt; pt, [X^{(k_1,n_1)}]\times [Y^{(k_2,n_2)}],
[Y^{(k_1,n_1)}]\times[X^{(k_2,n_2)}])\\
&& =\Psi_{L^{(k_1,n_1)},0,3}(pt; pt, [X^{(k_1,n_1)}],
[Y^{(k_1,n_1)}], pt)\\
&&\quad\cdot \Psi_{L^{(k_2,n_2)},0,3}(pt; pt,
[Y^{(k_2,n_2)}], [X^{(k_2,n_2)}])=1,\\
&&\Psi_{A,0,3}(pt; pt, [X^{(k_1,n_1)}]\times [X^{(k_2,n_2)}],
[Y^{(k_1,n_1)}]\times[Y^{(k_2,n_2)}])\\
&& =\Psi_{L^{(k_1,n_1)},0,3}(pt; pt, [X^{(k_1,n_1)}],
[Y^{(k_1,n_1)}])\\
&&\quad\cdot \Psi_{L^{(k_2,n_2)},0,3}(pt;pt, [X^{(k_2,n_2)}],
[Y^{(k_2,n_2)}])=1.
\end{eqnarray*}
As before it follows that
\begin{eqnarray*}
C_{HZ}^{(2\circ)}(W, \Omega; pt,[X^{(k_1,n_1)}]\times
[Y^{(k_2,n_2)}])&\le&\Omega(A)=(|a_1|+|a_2|)\pi,\\
C_{HZ}^{(2\circ)}(W, \Omega; pt,[X^{(k_1,n_1)}]\times
[X^{(k_2,n_2)}])&\le&\Omega(A)=(|a_1|+|a_2|)\pi,\\
C_{HZ}^{(2\circ)}(W, \Omega; pt,[Y^{(k_1,n_1)}]\times
[Y^{(k_2,n_2)}])&\le&\Omega(A)=(|a_1|+|a_2|)\pi,
\end{eqnarray*}
proving (\ref{e:1.21}).

To see (\ref{e:1.21.5}) we assume $r>1$ because of the result in
Theorem~\ref{t:1.15}. It immediately follows from (\ref{e:1.12})
and (\ref{e:1.20}) that
$$
{\cal W}_G(G(k_1, n_1)\times\cdots\times G(k_r,n_r),
\sigma^{(k_1,n_1)}\oplus\cdots\oplus\sigma^{(k_r,n_r)})\le\pi.
$$
On the another hand, by Lemma~\ref{lem:4.1} we have a symplectic
embedding from $(R_I(k_1, n_1)\times\cdots\times R_I(k_r, n_r),
\omega^{(k_1n_1)}\oplus\cdots\oplus\omega^{(k_rn_r)})$ to $(G(k_1,
n_1)\times\cdots\times G(k_r,n_r),
\sigma^{(k_1,n_1)}\oplus\cdots\oplus\sigma^{(k_r,n_r)})$.
Moreover, Lemma~\ref{lem:4.2} implies
that
\begin{eqnarray*}
B^{2k_1n_1+\cdots+ 2k_rn_r}(1)\!\!\!\!\!\!\!\!\!&&\subset
B^{2k_1n_1}(1)\times\cdots\times B^{2k_rn_r}(1)\\
&&\subset R_I(k_1, n_1)\times\cdots\times R_I(k_r, n_r).
\end{eqnarray*}
These give
$$
{\cal W}_G(G(k_1, n_1)\times\cdots\times G(k_r,n_r),
\sigma^{(k_1,n_1)}\oplus\cdots\oplus\sigma^{(k_r,n_r)})\ge\pi
$$
 and thus desired (\ref{e:1.21.5}).
\hfill$\Box$\vspace{2mm}

\noindent{\bf Proof of Theorem~\ref{t:1.17}}.\hspace{2mm} Without
loss of generality we may assume $a>0$. Firstly, as in the proof
of Theorem~\ref{t:1.16} one shows that
$$\Psi_{A,0,3}(pt; [M\times\CP^n], [M\times pt], pt)\ne 0$$
for $A=[pt\times\CP^1]$, and thus arrive at
\begin{equation}  \label{e:4.5}
  C_{HZ}^{(2\circ)}(M\times\CP^{n},\omega\oplus a\sigma_n; pt, [M\times
pt])\le a \pi.
\end{equation}
  Next we prove
\begin{equation}  \label{e:4.6}
\qquad C_{HZ}^{(2)}(M\times B^{2n}(r),\omega\oplus\omega_0; pt,
[M\times pt])= C_{HZ}^{(2)}(M\times
B^{2n}(r),\omega\oplus\omega_0; pt, pt).
\end{equation}
By Definition~\ref{def:1.2} it is clear that the left side in
(\ref{e:4.6}) is less than or equal to the right side in
(\ref{e:4.6}). To see the converse inequality we take $H\in{\cal
H}_{ad}(M\times B^{2n}(r), \omega\oplus\omega_0; pt, pt)$. Let
$P=P(H)$ and $Q=Q(H)$ be the corresponding submanifolds in
Definition~\ref{def:1.2}.
  Since
  $$
  P\subset Q\subset {\rm Int} (M\times B^{2n}(r))=M\times
{\rm Int}(B^{2n}(r))
$$
and $Q$ is compact there exists $\eta\in (0, r)$ such that
$Q\subset M\times B^{2n}(\eta)$. (Note that here we use $\partial
M=\emptyset$.) Therefore, $H$ may be viewed as an element of
${\cal H}_{ad}(M\times B^{2n}(r), \omega\oplus\omega_0; pt,
[M\times pt])$ naturally. This implies that the left side in
(\ref{e:4.6}) is more than or equal to the right side in
(\ref{e:4.6}).

Thirdly, as in \cite{HZ1,HZ2} one proves
\begin{equation}  \label{e:4.7}
  C_{HZ}^{(2)}(M\times B^{2n}(r),\omega\oplus\omega_0; pt,
pt)\ge\pi r^2
  \end{equation}
for any $r>0$. By (\ref{e:4.5}), Theorem~\ref{t:1.5} (v),
(\ref{e:4.6}) and (\ref{e:4.7})  we can obtain
\begin{eqnarray*}
a\pi&\ge &C_{HZ}^{(2\circ)}(M\times\CP^{n},\omega\oplus a\sigma_n;
pt,
[M\times pt])\\
&\ge&C_{HZ}^{(2)}(M\times\CP^{n},\omega\oplus a\sigma_n; pt, [M\times pt])\\
&\ge&C_{HZ}^{(2)}(M\times
B^{2n}(\delta\sqrt{a}),\omega\oplus\omega_0; pt, [M\times
pt])\\
&=&C_{HZ}^{(2)}(M\times
B^{2n}(\delta\sqrt{a}),\omega\oplus\omega_0; pt,
pt)\\
&\ge&\pi\delta^2 a
\end{eqnarray*}
for any $\delta\in (0, 1)$. Here we use the symplectic embedding
$(B^{2n}(\delta\sqrt{a}), \omega_0 )\hookrightarrow (\CP^n,
a\sigma_n)$ in the proof of Corollary 1.5 in \cite{HV2} for any
$0<\delta<1$. Taking $\delta \to 1$, we find that for $\delta =1$
the above inequalities are equalities. Together with
Lemma~\ref{lem:1.4} we obtain (\ref{e:1.22}) and $C(M\times
B^{2n}(r),\omega\oplus\omega_0)=\pi r^2$ in (\ref{e:1.23}).

To prove the other equality of (\ref{e:1.23}), i.e., $C(M\times
Z^{2n}(r),\omega\oplus\omega_0)=\pi r^2$, note that each
$H\in{\cal H}_{ad}(M\times Z^{2n}(r), \omega\oplus\omega_0;
pt,pt)$ can naturally be viewed as a function in ${\cal
H}_{ad}(M\times B^2(r)\times \R^{2n-2}/m\Z^{2n-2},
\omega\oplus\omega_0\oplus\omega_{st}; pt,pt)$ for sufficiently
large $m>0$. Here $\omega_{st}$ is the standard symplectic
structure on the tours $\R^{2n-2}/m\Z^{2n-2}$. It follows from the
equality just proved in (\ref{e:1.23}) that $\max H\le\pi r^2$ and
so
$$C_{HZ}^{(2\circ)}(M\times Z^{2n}(r),\omega\oplus\omega_0; pt,
pt)\le\pi r^2
$$
for any $r>0$. The desired conclusions easily follow.
\hfill$\Box$\vspace{2mm}

In order to prove Theorem~\ref{t:1.21} we need the following lemma
told to me by Professor Dusa McDuff and Dr.\ Felix Schlenk.

\begin{lemma}  \label{lem:4.3}
  For any two closed symplectic manifolds  $(M,\omega)$ and
$(N,\sigma)$ it holds that
$$c(M\times N,\omega\oplus\sigma) \ge c(M,\omega)+
c(N,\sigma)
$$
for $c=c_{HZ}$, $c_{HZ}^\circ$ and $C_{HZ}$, $C_{HZ}^\circ$.
\end{lemma}

According to Lemma~\ref{lem:1.4} it suffices to prove
Lemma~\ref{lem:4.3} for $c_{HZ}$ and $c_{HZ}^\circ$. Let $F$ and
$G$ be admissible functions on $M$ and $N$, respectively. Since
the Hamiltonian system for $F+G$ splits, we see that $F+G$ is an
admissible function on $M \times N$. From this
Lemma~\ref{lem:4.3} follows at once.

\medskip
\noindent{\bf Proof of Theorem~\ref{t:1.21}}. We denote by
$(W,\omega)$ the product manifold in Theorem~\ref{t:1.21}. Without
loss of generality we may assume $a_i>0$, $i=1,\cdots,k$. Let
$A_i=[\CP^1]$ be the generators of $H_2(\CP^{n_i};\Z)$,
$i=1,\cdots, k$. They are indecomposable classes. Since
$[Y^{(1,n_i)}]=pt$ it follows from the proof of
Theorem~\ref{t:1.16} that
$$\Psi_{A_i,0,3}(pt; pt, pt, [X^{(1,n_i)}])=1$$
for $i=1,\cdots,k$.  Set $A=A_1\times\cdots\times A_k$. Note that
each $(\CP^{n_i}, a_i\sigma_{n_i})$ is monotone. By
Proposition~\ref{prop:7.7} in the appendix we have
$$\Psi_{A, 0, 3}(pt; pt, pt,\beta) = 1$$
for some class $\beta\in H_\ast(W,\Q)$. Thus by
Corollary~\ref{cor:1.19} we get that
\begin{equation}\label{e:4.8}
c(W,\omega)\le\omega(A)=(a_1+\cdots + a_k)\pi
\end{equation}
for $c=c_{HZ}$, $c_{HZ}^\circ$. On the other hand,
Lemma~\ref{lem:4.3} yields
\begin{eqnarray*}
c(W,\omega) \,\ge\, \sum^k_{i=1}c(\CP^{n_i}, a_i\sigma_{n_i})
\,=\, (a_1+\cdots +a_k)\pi
\end{eqnarray*}
for $c=c_{HZ}$, $c_{HZ}^\circ$. \hfill$\Box$

\section{Proof of Theorems~\ref{t:1.22} and \ref{t:1.24}}

\noindent{\bf Proof of Theorem~\ref{t:1.22}}.\quad Under the
assumptions of Theorem~\ref{t:1.22} it follows from
Remark~\ref{rem:1.11} that the Gromov-Witten invariant
\[
\Psi_{A, g, m+2}(\pi^\ast C;\alpha_0, PD([\omega]),
\alpha_1,\cdots,\alpha_m)\ne 0,
\]
and thus Theorem~\ref{t:1.10} leads to
$$
C_{HZ}^{(2)}(M,\omega;\alpha_0, PD([\omega]))<+\infty.
$$
For a sufficiently small $\epsilon>0$ the well-known Lagrangian
neighborhood theorem due to Weinstein \cite{We1} yields a
symplectomorphism $\phi$ from $(U_\epsilon, \omega_{\rm can})$ to
a neighborhood  of $L$ in $(M,\omega)$ such that $\phi|_L=id$.
Since $L$ is a Lagrange submanifold one can, as in \cite{Lu3,V6},
use the Poincar\'e-Lefschetz duality theorem to prove that there
exists a cycle representative of $PD([\omega])$ whose support is
contained in $M\setminus \phi(U_\epsilon)$ because $\omega$ is
exact near $L$. By (\ref{e:1.6})  we get that
\begin{eqnarray}
  C_{HZ}^{(2)}(U_\epsilon,\omega_{\rm can};\tilde\alpha_0, pt)
  \!\!\!\!\!\!&&=C_{HZ}^{(2)}(\phi(U_\epsilon),\omega;\tilde\alpha_0,
pt)\label{e:5.1}\\
  &&\le C^{(2)}_{HZ}(M,\omega;\alpha_0,
  PD([\omega]))<+\infty.\nonumber
\end{eqnarray}
  Here we still denote by $\tilde\alpha_0$ the images in
$H_\ast(U_\varepsilon, \Q)$ and $H_\ast(\phi(U_\varepsilon), \Q)$
of $\tilde\alpha_0$ under the maps induced by the inclusions
$L\hookrightarrow U_\varepsilon$ and
$L\hookrightarrow\phi(U_\varepsilon)$. Note that for any
$\lambda\ne 0$ the map
$$\Phi_\lambda: T^\ast L\to T^\ast L,\; (q, v^\ast)\mapsto (q,\lambda
v^\ast),$$ satisfies $\Phi_\lambda^\ast\omega_{\rm
can}=\lambda\omega_{\rm can}$. Theorem~\ref{t:1.5} (iv),
(\ref{e:5.1}) and this fact imply  that
$$
C_{HZ}^{(2)}(U_c,\omega_{\rm can};\tilde\alpha_0, pt)<+\infty
$$
for any $c>0$.

In the case $g=0$, since the inclusion $L\hookrightarrow M$
induces an injective homomorphism $\pi_1(L)\to\pi_1(M)$ and thus
$\phi(U_\epsilon)\hookrightarrow M$ also induces an injective
homomorphism $\pi_1(\phi(U_\varepsilon))\to\pi_1(M)$ it follows
from (\ref{e:1.8}) that
\begin{eqnarray*}
C_{HZ}^{(2\circ)}(U_\epsilon,\omega_{\rm can};\tilde\alpha_0,
pt)\!\!\!\!\!\!\!\!&&=C_{HZ}^{(2\circ)}(\phi(U_\epsilon),\omega;\tilde\alpha_0,
pt)\\
&&\le C_{HZ}^{(2\circ)}(M,\omega;\alpha_0, PD([\omega]))<+\infty,
\end{eqnarray*}
and thus that $C_{HZ}^{(2\circ)}(U_c,\omega_{\rm
can};\tilde\alpha_0, pt)<+\infty$ for any $c>0$.

In particular, if $L$ is a Lagrange submanifold of a $g$-symplectic
uniruled manifold $(M,\omega)$, then we can take $\alpha_0=pt$ and
derive from (\ref{e:1.7})
$$
c_{HZ}(U_c,\omega_{\rm can})=C_{HZ}(U_c,\omega_{\rm can})<+\infty
$$
for any $c>0$, and from (\ref{e:1.9})
$$
c^\circ_{HZ}(U_c,\omega_{\rm can})=C^\circ_{HZ}(U_c,\omega_{\rm
can})<+\infty
$$
for all $c>0$ if $g=0$ and the inclusion $L\hookrightarrow M$
induces an injective homomorphism $\pi_1(L)\to\pi_1(M)$. Here we
use Lemma~\ref{lem:1.4} and the fact that $U_c$ is a compact
smooth manifold with connected boundary and of codimension zero
because $\dim L\ge 2$.

To see the final claim note that $(M,-\omega)$ is also strong
$g$-symplectic uniruled. It follows from
Proposition~\ref{prop:7.5} that the product $(M\times
M,(-\omega)\oplus\omega)$ is strong $0$-symplectic uniruled. By
the Lagrangian neighborhood theorem there exists a neighborhood
${\cal N}(\triangle)\subset M\times M$ of the diagonal
$\triangle$, a fiberwise convex neighborhood ${\cal N}(M_0)\subset
T^\ast M$ of the zero section $M_0$, and a symplectomorphism
$\psi:({\cal N}(\triangle), (-\omega)\oplus\omega)\to (T^\ast M,
\omega_{\rm can})$ such that $\psi(x,x)=(x,0)$ for $x\in M$. Note
also that the inclusion $\triangle\hookrightarrow M\times M$
induces an injective homomorphism
$\pi_1(\triangle)\to\pi_1(M\times M)$. The desired conclusion
follows immediately. \hfill$\Box$\vspace{2mm}

\noindent{\bf Proof of Theorem~\ref{t:1.24}}.\quad The case $\dim
M=2$ is obvious. So we assume that $\dim M\ge 4$. We follow
\cite{Bi1}. Let $p\, \colon L={\nu}(N)\to N$ be the symplectic
normal bundle of $N$ in $(M,\omega)$. It may naturally be viewed
as a complex line bundle with an obvious $S^1$-action
$$
t\cdot (b,v)=(b, e^{2\pi it}v),\quad (b,v)\in L \;{\rm
and}\;t\in S^1=\R/\Z.
$$
Consider the projectivized bundle $\pi:{\bf P}(L\oplus \C)\to N$
whose fiber at $b\in N$ is the complex projective space ${\bf
P}(L_b\oplus \C)$. This bundle has a natural $S^1$-action induced
by the action $t\cdot z=e^{-2\pi it}z$
of $S^1$ on each fiber summand of $\C$, i.e.,
$t\cdot(b,[v:z])=(b, [v:e^{-2\pi it}z])$. It has also two special
sections, the zero section $Z_0={\bf P}(\{0\}\oplus\C)$ and the
infinity section $Z_\infty={\bf P}(L\oplus\{0\})$.  One can
construct an $S^1$-invariant symplectic form on ${\bf P}(L\oplus
\C)$. Roughly speaking, fix any Hermitian metric $\|\cdot\|$ on
$L$ and denote by $p_N:S(L)=\{(b,v)\in L\,|\,\|v\|=1\}\to N$ the
associated unit circle bundle of $L$. The latter is a principal
$S^1$-bundle. Let $S^1=\R/\Z$ act on $\CP^1$ and $S(L)\times\CP^1$
by
\begin{eqnarray*}
&&t\cdot [z_0:z_1]=[z_0:e^{-2\pi it}z_1],\quad [z_0:z_1]\in\CP^1,\\
&&t\cdot \bigl((b,v), [z_0:z_1]\bigr)=\bigl((b, e^{-2\pi
it}v),[z_0:e^{-2\pi it}z_1]\bigr)
\end{eqnarray*}
for $(b,v)\in S(L)$ and $t\in S^1$. Then the quotient manifold
$S(L)\times_{S^1}\CP^1$ and ${\bf P}(L\oplus\C)$  can be
identified via the diffeomorphism induced by the projection
$$
\widetilde\Phi:S(L)\times\CP^1\to {\bf P}(L\oplus\C),\; \bigl((b,
v), [z_0:z_1]\bigr)\mapsto\bigl(b, [z_0v:z_1]\bigr).
$$
Under this identification one has
$Z_0=S(L)\times_{S^1}\{[0:1]\}$ and
$Z_\infty=S(L)\times_{S^1}\{[1:0]\}$. If $R^\nabla$ is the
curvature of the Hermitian connection $\nabla$ on $L$, then
$\rho_N:=\frac{1}{2\pi i}R^\nabla$ is a representing $2$-form of
the Chern class $c_1(L)$. Choose $0<\lambda_0<\varepsilon$ so that
$$
\tau_\lambda:=\omega|_N + \lambda\rho_N
$$
are symplectic forms on $N$ for all $0<\lambda\le\lambda_0$. Let
$h:\CP^1\to [0,1]$ be given by
$h([z_0:z_1])=|z_0|/(|z_0|^2+|z_1|^2)$. Define a map
$$
H_{\lambda_0}:S(L)\times_{S^1}\CP^1\to [0,\lambda_0],\quad [(b,v),
[z_0:z_1]]\mapsto\lambda_0h([z_0:z_1]).
$$
Then all level sets $H^{-1}_{\lambda_0}(\lambda)$, $\lambda \notin
\{ 0,\lambda_0 \}$, are diffeomorphic to $S(L)$, and the only critical
submanifolds of $H_{\lambda_0}$ are $H_{\lambda_0}^{-1}(0)=Z_0$ and
$H_{\lambda_0}^{-1}(\lambda_0)=Z_\infty$. As in Example~5.10 in
\cite{McSa1} (see also \cite{MWo}) one gets an $S^1$-invariant
symplectic form $\omega_{\lambda_0}$ on ${\bf P}(L\oplus\C)$ such
that $Z_0$, $Z_\infty$ and all fibers are symplectic submanifolds.
More precisely, $\omega_{\lambda_0}|_{Z_0}=\omega|_N$,
$\omega_{\lambda_0}|_{Z_\infty}=\omega|_N+ \lambda_0\rho_N$ and
each fiber ${\bf P}(L\oplus\C)_b\equiv\CP^1$ is equipped with an
$S^1$-invariant symplectic form with corresponding moment map
$\lambda_0h$, i.e., $\lambda_0\omega_{\rm FS}$. Here $\omega_{\rm
FS}$ is the standard Fubini-Study form on $\CP^1$ with
$\int_{\CP^1}\omega_{\rm FS}=1$. Furthermore, $Z_0$ has normal bundle in $({\bf
P}(L\oplus\C),\omega_{\lambda_0})$ with first Chern class
$[\rho_N]=c_1(L)$, see the appendix in \cite{MWo}.

As in \cite{Bi1}, from a transgression $1$-form
$\alpha^\nabla$ of the connection $\nabla$ on $L\setminus 0$ one can
get a $1$-form $\alpha\in\Omega^1(S(L))$ such that
$d\alpha=-p_N^\ast\rho_N$ and $\alpha(X)\equiv 1$, where
$X:S(L)\to TS(L)$ is the fundamental vector field of the above
$S^1$-action on $S(L)$. The first condition means that ${\cal
H}={\rm Ker}(\alpha)$ is the horizontal distribution of the
connection on $S(L)$ induced by $\nabla$. By Exercise 5.11 in
\cite{McSa1}, under the above identification ${\bf
P}(L\oplus\C)=S(L)\times_{S^1}\CP^1$, the symplectic form
$\omega_{\lambda_0}$ is induced by the
$S^1$-invariant closed $2$-form
$$
\Omega_{\lambda_0}:=p_N^\ast(\omega|_N)-
\lambda_0d(h\alpha)+\lambda_0\omega_{\rm FS}
$$
on $S(L)\times\CP^1$.

Set $X_L={\bf P}(L\oplus\C)$. Take an almost complex structure
$J_N\in{\cal J}(N,\omega|_N)$. By shrinking $\lambda_0>0$ we can
assume that $J_N$ is $(\omega|_N+\lambda\rho_N)$-tame for all
$0<\lambda\le\lambda_0$, i.e., $J_N\in{\cal J}_\tau(N,
\omega|_N+\lambda\rho_N)$. (This is not needed in the case of
\cite{Bi1} because $\rho_N$ can be taken as $\omega|_N$ there.)
Notice that the above horizontal distribution ${\cal H}$ over
$S(L)$ naturally induces a horizontal distribution
$\widetilde{\cal H}=\widetilde\Phi_\ast({\cal H}\times 0)$ on
$X_L$. So $TX_L=\widetilde{\cal H}\oplus\widetilde{\cal V}$, where
$\widetilde{\cal V}\subset TX_L$ is the vertical subbundle whose
fiber at $q\in X_L$ is $\widetilde{\cal V}_q={\rm
Ker}(d\pi(q))=T_q(X_L)_{\pi(q)}$. Actually $\widetilde{\cal H}$ is
exactly the horizontal distribution  on ${\bf P}(L\oplus\C)$
induced by the sum of the connection $\nabla$ on $L$ and the
trivial connection on $\C\to N$. Since
$\Omega_{\lambda_0}=p_N^\ast(\omega|_N+ \lambda_0h\rho_N)+
\lambda_0(\omega_{\rm FS}- dh\wedge\alpha)$, it is not hard to
check that for any $q\in X_L$ the subspaces $\widetilde{\cal H}_q$
and $\widetilde{\cal V}_q$ are
$(\omega_{\lambda_0})_q$-orthogonal. Similar to \cite{Bi1} we
construct an almost complex structure $J_X$ on $X_L$ as follows;
for any $q\in X_L$, ${J_X}|_{\widetilde{\cal H}_q}$  is the
horizontal lift of $(J_N)_q$ by the linear isomorphism
$d\pi(q)|_{\widetilde{\cal H}_q}: \widetilde{\cal H}_q\to
T_{\pi(q)}N$, and the restriction of ${J_X}$ to the fiber
$(X_L)_{\pi(q)}={\bf P}(L_{\pi(q)}\oplus\C)$ is the sum of the
complex structure determined by the Hermitian metric $\|\cdot\|$
on $L$ and of the standard one on $\C$. This $J_X$ is
$\omega_{\lambda_0}$-tame because $J_N\in{\cal J}_\tau(N,
\omega|_N+\lambda\rho_N)$ for all $0<\lambda\le\lambda_0$. One
easily sees that the almost complex structure $J_X$ is {\bf
fibred} on $X_L$ in the sense of Definition~2.2 of \cite{Mc2}.
Hence with $J_X$ we can prove as in Lemma~2.3 of \cite{Lu6} that
for
 the homology class $F\in H_2(X_L;\Z)$ of a fibre of
$X_L\to N$  the Gromov-Witten invariant
$$
\Psi^{(X_L,\omega_{\lambda_0})}_{F,0,3}(pt;[Z_0],[Z_\infty],pt)=1.
$$
That is, $(X_L,\omega_{\lambda_0})$ is a strong $0$-symplectic
uniruled manifold in the sense of Definition~1.14. By Theorem~1.10
we have
$$C_{HZ}^{(2\circ)}(X_L,\omega_{\lambda_0}; pt, [Z_\infty])\le
{\rm GW}_0(X_L,\omega_{\lambda_0}; pt, [Z_\infty])\le
\omega_{\lambda_0}(F)=\lambda_0.
$$
Note that for any $0<\delta<\lambda_0$ the set
$\{H_{\lambda_0}\le\delta\}$ is a smooth compact submanifold in
$X_L$ with connected boundary and of codimension zero that is a
neighborhood of $Z_0$ in $X_L$. It is easily seen that
 the inclusion
$\{H_{\lambda_0}\le\delta\}\subset X_L$ induces an injective
homomorphism $\pi_1(\{H_{\lambda_0}\le\delta\})\hookrightarrow
\pi_1(X_L)$. It follows from (\ref{e:1.9})  that
\[
c_{HZ}^\circ(\{H_{\lambda_0}\le\delta\},\omega_{\lambda_0})\le
C_{HZ}^{(2\circ)}(X_L,\omega_{\lambda_0}; pt,
[Z_\infty])\le\lambda_0.
\]
Identifying $N$ with the zero section $0_L$ and thus $Z_0$
 it follows from  the symplectic neighborhood theorem that for $\delta>0$ sufficiently small,
$(\{H_{\lambda_0}\le\delta\},\omega_{\lambda_0})$ is
symplectomorphic to a smooth compact submanifold $W\subset M$ with
connected boundary and of codimension zero that is a neighborhood
of $N$ in $M$.
Together with Lemma~\ref{lem:1.4} we therefore get
$$
c_{HZ}^\circ(W,\omega)=C_{HZ}^\circ(W,\omega)\le\lambda_0<\varepsilon.
$$
   The desired conclusion  is proved.
\hfill$\Box$\vspace{2mm}

\section{Proof of Theorem~\ref{t:1.35}}
The idea is the same as in \cite{Ka}. We can assume that $n/k\ge
2$. Following the notations in the proof of Theorem~\ref{t:1.15},
notice that the canonical atlas on $G(k, n)$ given by
(\ref{e:4.0}) has $n \choose k$ charts, and that for each chart
$$
(\Theta(\alpha_1,\cdots, \alpha_k), U(\alpha_1,\cdots, \alpha_k))
$$
Lemma~\ref{lem:4.1} yields a symplectic embedding
$\widehat\Phi_{\alpha_1\cdots\alpha_k}$ of $(R_I(k, m),
\omega^{(km)})$ into $(G(k, n), \sigma^{(k, n)})$ given by
$$Z\mapsto [(\sqrt{I^{(k)}-Z\overline Z^\prime},
Z)P(\alpha_1,\cdots,\alpha_k)],
$$ where
$P(\alpha_1,\cdots,\alpha_k)$ is the $n\times n$ permutation
matrix such that (\ref{e:4.00}) holds for the matrix $B=(I^{(k)},
Z)P(\alpha_1,\cdots,\alpha_k)$. Moreover, for the matrix
$A=(\sqrt{I^{(k)}-Z\overline Z^\prime},
Z)P(\alpha_1,\cdots,\alpha_k)$ we have
$$
A_{\alpha_1\cdots\alpha_k}=\sqrt{I^{(k)}-Z\overline
Z^\prime}\quad{\rm and}\quad A_{\alpha_{k+1}\cdots\alpha_n}=Z.
$$
Note that
\begin{eqnarray*}
\|A\|^2&=&\|A_{\alpha_1\cdots\alpha_k}\|^2+\|A_{\alpha_{k+1}\cdots\alpha_n}\|^2\\
&=&{\rm tr}(A_{\alpha_1\cdots\alpha_k}\overline
A^\prime_{\alpha_1\cdots\alpha_k})+ {\rm
tr}(A_{\alpha_{k+1}\cdots\alpha_n}\overline
A^\prime_{\alpha_{k+1}\cdots\alpha_n})
\\
&=&{\rm tr}(I^{(k)}-Z\overline Z^\prime)+ {\rm tr}(Z\overline
Z^\prime)=k,
\end{eqnarray*}
and therefore
$$
\|A_{\alpha_1\cdots\alpha_k}\|^2=k-\|A_{\alpha_{k+1}\cdots\alpha_n}\|^2=k-\|Z\|^2.
$$
By Lemma~\ref{lem:4.2} these show that
$\widehat\Phi_{\alpha_1\cdots\alpha_k}(B^{2km}(r))$ is contained
in
\begin{equation}  \label{e:6.1}
\left.\begin{array}{ll} & \Lambda(\alpha_1,\cdots,\alpha_k; r)=
\{[B]\in G(k, n) \mid \mbox{ for all } A\in [B]\cap M^0(k, n;\C),\\
&\qquad\qquad\qquad\qquad\qquad\qquad\qquad\qquad\;
\|A_{\alpha_1\cdots\alpha_k}\|^2>k-r^2\}
\end{array}\right.
\end{equation}
for any $0<r\le 1$. Note that $k>1$ and $n/k\ge 2$. There must be
two disjoint subsets of $\{1,\cdots, n\}$, say
$\{\alpha_1,\cdots,\alpha_k\}$ and $\{\beta_1,\cdots,\beta_k\}$,
such that $\alpha_1<\cdots<\alpha_k$ and $\beta_1<\cdots<\beta_k$.
For any two such subsets we claim that
$$\widehat\Phi_{\alpha_1\cdots\alpha_k}(B^{2km}(1))\cap
\widehat\Phi_{\beta_1\cdots\beta_k}(B^{2km}(1))=\emptyset.
$$
In fact,  $\Lambda(\alpha_1,\cdots,\alpha_k; 1)$ and
$\Lambda(\beta_1,\cdots,\beta_k; 1)$ are disjoint. Otherwise, let
$[B]$ belong to their intersection and take a representative $A$
of $[B]$ in $M^0(k,n;\C)$. Then
$$k\ge \|A_{\alpha_1\cdots\alpha_k}\|^2+
\|A_{\beta_1\cdots\beta_k}\|^2>2k-2$$ by (\ref{e:6.1}). This
contradicts the assumption that $k\ge 2$. Now the conclusion
follows from the fact that there exist exactly $[n/k]$ mutually
disjoint subsets of $\{1,\cdots, n\}$  consisting of $k$ numbers.
\hfill$\Box$\vspace{2mm}

\noindent{\bf Proof of (\ref{e:1.26}).}\hspace{2mm}Notice that
$G(k,n)$ can be embedded into the complex projective space $\CP^N$
with $N=n!/(n-k)!k!-1$ by the Pl\"ucker map $p$ (\cite{GH}), and
that for any $l$-dimensional subvariety $X$ of $\CP^N$ one has
$${\rm Vol}(X)=\deg(X)\cdot {\rm Vol}(L)$$
with respect to the Fubini-Study metric, where $L$ is an
$l$-dimensional linear subspace of $\CP^N$ (cf.\ \cite[p.\
384]{Fu}). But it was shown in Example 14.7.11 of \cite{Fu} that
$$\deg(p(G(k,n)))=\frac{1!\cdot 2!\cdots (k-1)!\cdot (k(n-k))!}
{(n-k)!\cdot (n-k+1)!\cdots (n-1)!}.$$ It is well-known that the
volume of a $k(n-k)$-dimensional linear subspace $L$
  of $\CP^N$ is
$${\rm Vol}(\CP^{k(n-k)})=\frac{\pi^{k(n-k)}}{(k(n-k))!}.$$
These give (\ref{e:1.26}). \hfill$\Box$\vspace{2mm}

\section{Appendix: The Gromov-Witten invariants of product manifolds}

In this appendix we collect some results on Gromov-Witten
invariants needed in this paper. They either are easily proved or
follow from the references given.

Let $(V,\omega)$ be a closed symplectic manifold of dimension
$2n$. Recall that for a given class $A\in H_2(V;\Z)$ the
Gromov-Witten invariant of genus $g$ and with $k$ marked points is
a homomorphism
$$
\Psi^V_{A, g, k}: H_\ast(\overline{\cal M}_{g, k};\Q)\times
H_\ast(V;\Q)^k\to \Q,
$$
where $2g+k\ge 3$ and $\overline{\cal M}_{g, k}$ is the space of
isomorphism classes of genus $g$ stable curves with $k$ marked
points, which is a connected K\"ahler orbifold of complex
dimension $3g-3+k$. In \cite{Lu8} we used the cohomology
$H^\ast_c(V;\Q)$ with compact support and the different notation
${\mathcal G}{\mathcal W}^{(\omega,\mu,J)}_{A,g,k}$ to denote the
GW-invariants since we also considered noncompact symplectic
manifolds for which the dependence on further data needs to be
indicated.  For closed symplectic manifolds we easily translate
the composition law and reduction formulas in \cite{Lu8} into the
homology version, which is the same as the ones in \cite{RT2}. Let
integers $g_i\ge 0$ and $k_i>0$ satisfy $2g_i+ k_i\ge 3$, $i=1,
2$. Set $g=g_1+ g_2$ and $k=k_1+ k_2$ and fix a decomposition
$S=S_1\cup S_2$ of $\{1,\cdots, k\}$ with $|S_i|=k_i$. Then there
is a canonical embedding
\begin{equation}\label{e:7.1}
\theta_S:\overline{\cal M}_{g_1, k_1+1}\times\overline{\cal
M}_{g_2, k_2+1}\to \overline{\cal M}_{g, k},
\end{equation}
which assigns to marked curves $(\Sigma_i;x_1^i,\cdots,
x^i_{k_i+1})$, $i=1,2$, their union $\Sigma_1\cup\Sigma_2$ with
$x^1_{k_1+1}$ and $x^2_{k_2+1}$ identified and the remaining
points renumbered by $\{1,\cdots, k\}$ according to $S$. Let
$$\mu_{g,k}: \overline{\cal M}_{g-1, k+2}\to\overline{\cal M}_{g, k}
$$
be the map corresponding to gluing together the last two marked
points. It is continuous. Suppose that $\{\beta_b\}^L_{b=1}$ is a
homogeneous basis of $H_\ast(V;\Z)$ modulo torsion, $(\eta_{ab})$
its intersection matrix and
$(\eta^{ab})=(\eta_{ab})^{-1}$.\vspace{2mm}

\noindent{\bf Composition law}.\hspace{2mm}{\it Let
$$
[K_i]\in H_\ast(\overline{\cal M}_{g_i, k_i+1};\Q),\quad
i=1,2,\quad [K_0]\in H_\ast(\overline{\cal M}_{g-1, k+2};\Q)
$$
and $A\in H_2(V;\Z)$. Then for any $\alpha_1,\cdots,\alpha_k$ in
$H_\ast(V;\Q)$ we have
\begin{eqnarray*}
& &\Psi^V_{A,g,k}({\theta_S}_\ast([K_1\times
K_2]);\alpha_1,\cdots, \alpha_k) = (-1)^{{\rm
cod}(K_2)\sum^{k_1}_{i=1}{\rm cod}(\alpha_i)}\\
\sum_{A=A_1+A_2}\sum_{a,b}\!\!\!\!\!\!\!\!&&\Psi^V_{A_1,
g_1,k_1+1}([K_1]; \{\alpha_i\}_{i\le k_1},
\beta_a)\eta^{ab}\Psi^V_{A_2, g_2,k_2+1}([K_2]; \beta_b,
\{\alpha_j\}_{j>k_1}),
\end{eqnarray*}
$$\Psi^V_{A,g,k}((\mu_{g,k})_\ast([K_0]);\alpha_1,\cdots,\alpha_k)=\sum_{a,b}\Psi^V_{A,
g-1,k+2}([K_0];\alpha_1,\cdots,\alpha_k,\beta_a,
\beta_b)\eta^{ab}.
$$}

Remark that $(-1)^{{\rm cod}(K_2)\sum^{k_1}_{i=1}{\rm
cod}(\alpha_i)}=(-1)^{{\rm dim}(K_2)\sum^{k_1}_{i=1}{\rm
dim}(\alpha_i)}$ because the dimensions of $\overline{\cal
M}_{g_i, k_i+1}$ and $V$ are even.
  Denote the map forgetting the last marked point
by
$$\pi_k: \overline{\cal M}_{g, k}\to\overline{\cal M}_{g, k-1}.$$

\noindent{\bf Reduction formula}.\hspace{2mm}{\it Suppose that
$(g,k)\ne (0,3), (1,1)$. Then
\begin{description}
\item[(i)] for any $\alpha_1,\cdots,\alpha_{k-1}$ in
$H_\ast(V;\Q)$ and $[K]\in H_\ast(\overline{\cal M}_{g, k};\Q)$ we
have
\begin{equation}\label{e:7.2}
\quad\Psi^V_{A,g,k}([K];\alpha_1,\cdots,\alpha_{k-1},
[V])=\Psi^V_{A,g,k-1}((\pi_k)_\ast([K]);\alpha_1,\cdots,\alpha_{k-1})
\end{equation}
\item[(ii)] if $\alpha_k\in H_{2n-2}(V;\Q)$ we have
\begin{equation}\label{e:7.3}
\Psi^V_{A,g,k}([\pi_k^{-1}(K)];\alpha_1,\cdots,\alpha_k)=PD(\alpha_k)(A)\Psi^V_{A,g,k-1}([K];\alpha_1,\cdots,\alpha_{k-1}).
\end{equation}
\end{description}}

\begin{lemma}\label{lem:7.1}
Let $(V,\omega)$ be a closed symplectic manifold,
$\{\beta_b\}^L_{b=1}$ a homogeneous basis of $H_\ast(V;\Z)$ modulo
torsion as in the composition law above. Suppose that there exist
homology classes $A\in H_2(V;\Z)$, $\alpha_1,\cdots,\alpha_m\in
H_\ast(V;\Q)$ and $g>0$ such that
\begin{equation}  \label{e:7.4}
\Psi^V_{A,g,m}(pt; \alpha_1,\cdots,\alpha_m)\ne 0.
\end{equation}
Then for each nonnegative integer $g'<g$ we have
$$
\Psi^V_{A, g', m+2s}(pt;
\alpha_1,\cdots,\alpha_m,\beta_{a_1},\beta_{b_1},
\cdots,\beta_{a_s},\beta_{b_s})\ne 0
$$
for $s=g-g'$ and some $\beta_{a_i},\beta_{b_i}$ in
$\{\beta_b\}^L_{b=1}$, $i=1,\cdots,s$.
\end{lemma}

\noindent{\bf Proof.} By the composition law for Gromov-Witten
invariants
  we have
\begin{eqnarray*}
\Psi^V_{A,g,m}(pt; \alpha_1,\cdots,\alpha_m)&=&
\Psi^V_{A,g,m}((\mu_{g,m})_\ast (pt); \alpha_1,\cdots,\alpha_m)\\
&=&\sum_{a,b}\Psi^V_{A, g-1,m+2}(pt;
\alpha_1,\cdots,\alpha_m,\beta_a,\beta_b) \eta^{ab}.
\end{eqnarray*}
By (\ref{e:7.4}),  the left side is not equal to zero. So there
exists a pair $(a,b)$ such that
$$\Psi^V_{A, g-1,m+2}(pt;
\alpha_1,\cdots,\alpha_m,\beta_a,\beta_b)\ne 0.
$$
If $g-1>g'$ we can repeat this argument to reduce $g-1$. After
$s=g-g'$ steps the lemma follows. \hfill$\Box$\vspace{2mm}

\begin{lemma}\label{lem:7.2}
Let $(V,\omega)$ and $\{\beta_b\}^L_{b=1}$ be as in
Lemma~\ref{lem:7.1}.  Suppose that there exist homology classes
$$
A\in H_2(V;\Z),\; \xi_1,\cdots,\xi_k\in H_\ast(V;\Q)\quad{\rm
and}\quad [K] \in H_\ast(\overline{\cal M}_{g,k};\Q)
$$
such that
\begin{equation}  \label{e:7.5}
\Psi^V_{A,g,k}([K]; \xi_1,\cdots,\xi_k)\ne 0
\end{equation}
for some integer $g\ge 0$. Then for each integer $m>k$ we have
$$
\Psi^V_{A,g, m}([K_1]; \xi_1,\cdots,\xi_k,
\overbrace{PD([\omega]),\cdots, PD([\omega])}^{m-k})\ne 0.
$$
  Here $K_1=(\pi_m)^{-1}\circ\cdots\circ (\pi_{k+1})^{-1}(K)$.
\end{lemma}

\noindent{\bf Proof.} Using the definition of the GW-invariants,
it follows from (\ref{e:7.5}) that $2g+k\ge 3$ and that the space
$\overline{\mathcal M}_{g,k}(V, J, A)$ of $k$-pointed stable
$J$-maps of genus $g$ and of class $A$ in $V$ is nonempty for
generic $J\in{\mathcal J}(V,\omega)$. In particular, this implies
$\omega(A)\ne 0$. Applying the reduction formula (\ref{e:7.3}) to
(\ref{e:7.5}) we have
$$
\Psi^V_{A,g, k+1}([\pi_{k+1}^{-1}(K)]; \xi_1,\cdots,\xi_k,
PD([\omega]))=\omega(A)\cdot \Psi^V_{A,g, k}([K];
\xi_1,\cdots,\xi_k)\ne 0.
$$
Continuing this process $m-k-1$ times again we get the desired
conclusion. \hfill$\Box$\vspace{2mm}

\begin{proposition}  \label{prop:7.3}
For a closed symplectic manifold $(V,\omega)$, if there exist
homology classes
  $A\in H_2(V;\Z)$ and $\alpha_i\in H_\ast(V;\Q)$,
$i=1,\cdots, k$, such that the Gromov-Witten invariant
\begin{equation}  \label{e:7.7}
\Psi_{A, g, k+1}(pt; pt,\alpha_1,\cdots, \alpha_k)\ne 0
\end{equation}
for some integer $g\ge 0$, then there exist homology classes $B\in
H_2(V;\Z)$ and $\beta_1, \beta_2 \in H_\ast(V;\Q)$ such that
\begin{equation}  \label{e:7.8}
\Psi_{B, 0, 3}(pt; pt,\beta_1,\beta_2)\ne 0.
\end{equation}
Consequently, every strong symplectic uniruled manifold is strong
$0$-symplectic uniruled.
\end{proposition}

  (\ref{e:7.8}) implies that $B$ is  spherical. In fact, in this case
   there exists a $3$-pointed stable $J$-curve of genus zero and
in class $B$. By the gluing arguments we can get a $J$-holomorphic
sphere $f:\CP^1\to M$ which represents the class $B$. That is, $B$
is $J$-effective. So $B$ is necessarily spherical, cf.\ Page 67 in
\cite{McSa2}.

\noindent{\bf Proof of Proposition~\ref{prop:7.3}}.\hspace{2mm} By
Lemma~\ref{lem:7.1}, we can assume $g=0$ in (\ref{e:7.7}), i.e.\
\begin{equation}  \label{e:7.9}
\Psi_{A, 0, k+1}(pt; pt,\alpha_1,\cdots, \alpha_k)\ne 0.
\end{equation}
This implies that $k+1\ge 3$ or $k\ge 2$. If $k=2$ then the
conclusion holds. If $k=3$ we can use the reduction formula
(\ref{e:7.2}) to get
$$\Psi_{A, 0, 5}(pt; pt, \alpha_1,\cdots,\alpha_3, [V])
=\Psi_{A, 0, 4}(pt;pt, \alpha_1,\cdots,\alpha_3)\ne 0.
$$
  Therefore we can actually assume
that $k\ge 4$ in (\ref{e:7.9}).  Since $\overline{\cal M}_{0,m}$
is connected for every integer $m\ge 3$, $H_0(\overline{\cal
M}_{0,m},\Q)$ is generated by pt.   For the canonical embedding
$\theta_S$ as in (\ref{e:7.1}) we have ${\theta_S}_\ast(pt\times
pt)=pt$. Hence it follows from the composition law  that
\begin{eqnarray*}
&& \Psi_{A, 0, k+1}(pt;pt,\alpha_1,\cdots,\alpha_k)\\
&&=\sum_{A=A_1+ A_2}\sum_{a,b}\Psi_{A_1, 0, 4}(pt; pt, \alpha_1,
\alpha_2,\beta_a) \eta^{ab}\Psi_{A_2, 0,
k-1}(pt;\beta_b,\alpha_3,\cdots,\alpha_k)
\end{eqnarray*}
because ${\rm cod}(K_2)={\rm cod}(pt)$ is even. This implies that
\begin{equation}\label{e:7.10}
\Psi_{A_1, 0, 4}(pt; pt, \alpha_1, \alpha_2,\beta_a)\ne 0
\end{equation}
for some $A_1\in H_2(V;\Z)$ and $1\le a\le L$. By the
associativity of the quantum multiplication,
\begin{eqnarray*}
&&\Psi_{A_1, 0, 4}(pt; pt, \alpha_1,\alpha_2,\beta_a)= \\
&&\qquad\qquad\pm\sum_{A_1=A_{11}+ A_{12}}\sum_l\Psi_{A_{11}, 0,
3}(pt; pt, \alpha_1, e_l) \Psi_{A_{12}, 0, 3}(pt;
f_l,\alpha_2,\beta_a)
\end{eqnarray*}
where $\{e_l\}_l$ is a basis for the homology $H_\ast(M;\Q)$ and
$\{f_l\}_l$ is the dual basis with respect to the intersection
pairing, see (6) in \cite{Mc2}. It follows from this identity and
(\ref{e:7.10}) that
$$ \Psi_{A_{11}, 0, 3}(pt;  pt, \alpha_1, e_l)\ne 0$$
for some $l$. Taking $B=A_{11}$ we get (\ref{e:7.8}).
\hfill$\Box$\vspace{2mm}

\begin{proposition}  \label{prop:7.4}
Let $(M,\omega)$ and $(N,\sigma)$ be two closed symplectic
manifolds. Then for every integer $k\ge 3$ and homology classes
$A_2\in H_2(N;\Z)$ and $\beta_i\in H_\ast(N;\Q)$,
   $i=1,\cdots, k$,
$$\Psi^{M\times N}_{0\oplus A_2, 0, k}(pt; [M]\otimes\beta_1,
\cdots, [M]\otimes\beta_{k-1}, pt\otimes\beta_k)=\Psi^N_{A_2, 0,
k}(pt;\beta_1,\cdots,\beta_k),$$ where $0\in H_2(M;\Z)$ denotes
the zero class.
\end{proposition}

\noindent{\bf Proof}.  Take $J_M\in{\cal J}(M,\omega)$,
$J_N\in{\cal J}(N,\sigma)$ and set $J=J_M\times J_N$. Note that
the product symplectic manifold $(M\times N,\omega\oplus\sigma)$
is a special symplectic fibre bundle over $(M,\omega)$ with fibres
$(N,\sigma)$. Moreover, the almost complex structure $J=J_M\times
J_N$ on $M\times N$ is {\bf fibred} in the sense of Definition 2.2
in \cite{Mc2}. So for a fibre class $0\oplus A_2$ we can, as in
the proof of Proposition 4.4 of \cite{Mc2}, construct a virtual
moduli cycle ${\overline M}^\nu_{0,3}(M\times N, J, 0\oplus A_2)$
of ${\overline M}_{0,3}(M\times N, J, 0\oplus A_2)$ such that the
$M$-components of each element in ${\overline M}^\nu_{0,3}(M\times
N, J, 0\oplus A_2)$ are $J_M$-holomorphic, and thus constant. This
shows that the virtual moduli cycle ${\overline
M}^\nu_{0,3}(M\times N, J, 0\oplus A_2)$ may be chosen as $M\times
{\overline M}^\nu_{0,3}(N, J_N, A_2)$. The desired conclusion
follows. These techniques were also used in the proof of Lemma 2.3
in \cite{Lu6}. We refer to there and $\S4.3$ in \cite{Mc2} for
more details. \hfill$\Box$\vspace{2mm}

As a direct consequence of Proposition~\ref{prop:7.3} and
Proposition~\ref{prop:7.4} we get

\begin{proposition}  \label{prop:7.5}
The product of a closed symplectic manifold and a strong symplectic
uniruled manifold is strong $0$-symplectic uniruled. In particular,
the product of finitely many strong symplectic uniruled manifolds is
also strong $0$-symplectic uniruled.
\end{proposition}

Actually we can generalize Proposition~\ref{prop:7.4} to a
symplectic fibre bundle over a closed symplectic manifold with a
closed symplectic manifold as fibre. Therefore, a symplectic fibre
bundle over a closed symplectic manifold with   a strong symplectic
uniruled fibre is also strong symplectic uniruled.

In the proof of Theorem~\ref{t:1.21} we need a product formula for
Gromov-Witten invariants. Such a formula was given for algebraic
geometry GW-invariants of two projective algebraic manifolds in
\cite{B}. However it is not clear whether the GW-invariants used
in this paper agree with those of \cite{B} for projective
algebraic manifolds. For the sake of simplicity we shall give a
product formula for a special case, which is sufficient
  for the proof of Theorem~\ref{t:1.21}. Recall
that a symplectic manifold $(M,\omega)$ is said to be {\bf
monotone} if there exists a number $\lambda>0$ such that
$\omega(A)=\lambda c_1(A)$ for $A\in\pi_2(M)$. The {\bf minimal
Chern number} $N\ge 0$ of a symplectic manifold $(M,\omega)$ is
defined by $\langle c_1,\pi_2(M)\rangle=N\Z$.
  For $J\in{\mathcal
J}(M,\omega)$, a homology class $A\in H_2(M,\Z)$ is called J-{\bf
effective} if it can be represented by a $J$-holomorphic sphere
$u:\CP^1\to M$. Such a homology class must be spherical. Moreover,
a class $A\in H_2(M,\Z)$ is called {\bf indecomposable} if it
cannot be decomposed as a sum $A=A_1+\cdots + A_k$ of classes
which are spherical and satisfy $\omega(A_i)>0$ for
$i=1,\cdots,k$.

\begin{proposition}  \label{prop:7.6}
{\it Let the closed symplectic manifold $(M,\omega)$ either be
monotone or have minimal Chern number $N\ge 2$. Then for each
indecomposable class $A\in H_2(M,\Z)$ and classes $\alpha_i\in
H_\ast(M,\Z)$, $i=1,2,3$ the Gromov-Witten invariant
$\Psi^M_{A,0,3}(pt;\alpha_1,\alpha_2,\alpha_3)$ adopted in this
paper agrees with the invariant
$\Psi^M_{A,3}(\alpha_1,\alpha_2,\alpha_3)$ in \S7.4 of
\cite{McSa2}.}
\end{proposition}

\noindent{\bf Proof.}\quad Let $J\in{\mathcal J}(M,\omega)$.
Consider the space $\overline{\mathcal M}_{0,3}(M, A, J)$ of
equivalence classes of all $3$-pointed stable $J$-maps of genus
zero and of class $A$ in $M$. For $[{\bf f}]\in\overline{\mathcal
M}_{0,3}(M, A, J)$, since $A$ is indecomposable it follows from
the definition of stable maps that ${\bf
f}=(\Sigma;z_1,z_2,z_3;f)$ must be one of the following four
cases:
\begin{description}
\item[(a)] The domain $\Sigma=\CP^1$, $z_i$, $i=1,2,3$ are three
distinct marked points on $\Sigma$, and $f:\Sigma\to M$ is a
$J$-holomorphic map of class $A$.

\item[(b)] The domain $\Sigma$ has exactly two components
$\Sigma_1=\CP^1$ and $\Sigma_2=\CP^1$ which have a unique
intersection point. $f|_{\Sigma_1}$ is nonconstant and $\Sigma_1$
only contains one marked point. $f|_{\Sigma_2}$ is constant and
$\Sigma_2$ contains two marked points.

\item[(c)] The domain $\Sigma$ has exactly two components
$\Sigma_1=\CP^1$ and $\Sigma_2=\CP^1$ which have a unique
intersection point. $f|_{\Sigma_1}$ is nonconstant and $\Sigma_1$
contains no marked point. $f|_{\Sigma_2}$ is constant and
$\Sigma_2$ contains three marked points.

\item[(d)] The domain $\Sigma$ has exactly three components
$\Sigma_1=\CP^1$, $\Sigma_2=\CP^1     $ and $\Sigma_3=\CP^1$.
$\Sigma_1$ and $\Sigma_2$ (resp.\ $\Sigma_2$ and $\Sigma_3$) have
only one intersection point, and $\Sigma_1$ and $\Sigma_3$ have no
intersection point.  $f|_{\Sigma_1}$ is nonconstant and $\Sigma_1$
contains no marked point. $f|_{\Sigma_2}$ is constant and
$\Sigma_2$ contains one marked point. $f|_{\Sigma_3}$ is constant
and $\Sigma_3$ contains two marked points.
\end{description}

Let $\overline{\mathcal M}_{0,3}(M, A, J)_i$, $i=1,2,3,4$ be the
subsets of the four kinds of stable maps. It is easily proved that
for generic $J\in{\mathcal J}(M,\omega)$ they are smooth manifolds
of dimensions
\begin{eqnarray*}
&&\dim \overline{\mathcal M}_{0,3}(M, A, J)_1=\dim M+ 2c_1(A),\\
&& \dim \overline{\mathcal M}_{0,3}(M, A, J)_2=\dim M+ 2c_1(A)-4, \\
&&\dim\overline{\mathcal M}_{0,3}(M, A, J)_3=\dim M+ 2c_1(A)-6, \\
&&\dim\overline{\mathcal M}_{0,3}(M, A, J)_4=\dim M+ 2c_1(A)-6.
  \end{eqnarray*}
So $\overline{\mathcal M}_{0,3}(M, A,
J)=\cup^4_{i=1}\overline{\mathcal M}_{0,3}(M, A, J)_i$ is a
stratified smooth compact manifold. Note that each stable map in
$\overline{\mathcal M}_{0,3}(M, A, J)$ has no free components. The
construction of the virtual moduli cycle in \cite{Lu8} with
Liu-Tian's method in \cite{LiuT} is thus trivial or not needed:
The virtual moduli cycle of $\overline{\mathcal M}_{0,3}(M, A, J)$
may be taken as
$$
\overline{\mathcal M}_{0,3}(M, A, J)\to {\mathcal
B}^M_{0,3,A},\;[{\bf f}]\mapsto [{\bf f}],
$$
where ${\mathcal B}^M_{0,3,A}$ is the space of equivalence classes
of all $3$-pointed stable $L^{k,p}$-maps of genus zero and of
class $A$ in $M$. Therefore for homology classes $\alpha_i\in
H_2(M,\Z)$, $i=1,2,3$, satisfying the dimension condition
$$\deg(\alpha_1)+\deg(\alpha_2)+\deg(\alpha_3)=2n+ 2c_1(A)$$
the Gromov-Witten invariant
\begin{equation}\label{e:7.11}
\begin{array}{cr}
\hspace{-10mm}\Psi^M_{A,0,3}(pt;\alpha_1,\alpha_2,\alpha_3)=({\rm
EV}^{J,A}_{0,3})\cdot(\overline\alpha_1\times\overline\alpha_2\times\overline\alpha_3)\\
\hspace{42mm}=({\rm EV}^{J,A}_{0,3})|_{\overline{\mathcal
M}_{0,3}(M, A, J)_1}
\cdot(\overline\alpha_1\times\overline\alpha_2\times\overline\alpha_3)
\end{array}
\end{equation}
because the intersections can only occur in the top strata. Here
\begin{equation}\label{e:7.12}
{\rm EV}^{J,A}_{0,3}:\overline{\mathcal M}_{0,3}(M, A, J)\to
M^3,\;[{\bf f}]\mapsto (f(z_1), f(z_2), f(z_3)),
\end{equation}
and $\overline\alpha_i: U_i\to M$ are generic pseudocycle
representatives of the classes $\alpha_i$, $i=1,2,3$, cf.
\cite{McSa2} for details. Note that each element $[{\bf f}]$ in
$\overline{\mathcal M}_{0,3}(M, A, J)_1$ has a unique
representative of the form $(\CP^1; 0,1,\infty; f)$. So
$\overline{\mathcal M}_{0,3}(M, A, J)_1$ may be identified with
the space ${\mathcal M}(M, A, J)$ of all $J$-holomorphic curves
which represent the class $A$. Fix marked points ${\bf
z}=(0,1,\infty)\in (\CP^1)^3$ and define the evaluation map
\begin{equation}\label{e:7.13}
E_{A,J,{\bf z}}:{\mathcal M}(M, A, J)\to M^3,\;f\mapsto (f(0),
f(1), f(\infty)).
\end{equation}
 From the above arguments one easily
checks that it is a pseudocycle in the sense of \cite{McSa2}. Then
(\ref{e:7.11}) gives rise to
\begin{equation}\label{e:7.14}
\Psi^M_{A,0,3}(pt;\alpha_1,\alpha_2,\alpha_3)=E_{A,J,{\bf z}}
\cdot(\overline\alpha_1\times\overline\alpha_2\times\overline\alpha_3)=\Psi^M_{A,3}(\alpha_1,\alpha_2,\alpha_3)
\end{equation}
because we can require that
$\overline\alpha_1\times\overline\alpha_2\times\overline\alpha_3$
is also transverse to $E_{A,J,{\bf z}}$.
  \hfill$\Box$\vspace{2mm}

\begin{proposition}  \label{prop:7.7}
{\it Consider closed symplectic manifolds $(M_k,\omega_k)$ as in
Proposition~\ref{prop:7.6} and indecomposable classes $A_k\in
H_2(M_k,\Z)$, $k=1,\cdots, m$. Then for $\alpha^{(k)}_i\in
H_\ast(M_k,\Z)$, $i=1,2,3$ and $k=1,\cdots, m$ we have the
Gromov-Witten invariant
\begin{equation}\label{e:7.15}
\hspace{5mm}\Psi^M_{A,0,3}(pt;\times^m_{k=1}\alpha^{(k)}_1,
\times^m_{k=1}\alpha^{(k)}_2,\times^m_{k=1}\alpha^{(k)}_3)=\prod^m_{k=1}\Psi^{M_k}_{A_k,0,3}(pt;\alpha^{(k)}_1,
\alpha^{(k)}_2,\alpha^{(k)}_3),
\end{equation}
where $A=\oplus^m_{k=1}A_k$.}
\end{proposition}

\noindent{\bf Proof.}\quad Set
$(M,\omega)=(\times^m_{k=1}M_k,\times^m_{k=1}\omega_k)$. Take
$J_k\in{\mathcal J}(M_k,\omega_k)$, $k=1,\cdots,m$ and set
$J=\times^m_{k=1}J_k$. Then $J\in {\mathcal J}(M,\omega)$.
  It is not hard to prove that for generic $J_k\in{\mathcal
J}(M_k,\omega_k)$ the space  $\overline{\mathcal M}_{0,3}(M, A,
J)$ is still a stratified smooth compact manifold. We still denote
by $\overline{\mathcal M}_{0,3}(M, A, J)_1$ its top stratum, which
consists of elements $[{\bf f}]\in \overline{\mathcal M}_{0,3}(M,
A, J)$ whose domain has only one component $\CP^1$. It is a smooth
noncompact manifold of dimension $\dim M+ 2c_1(A)=\sum^m_{k=1}
\dim M_k+ 2c_1(A_k)$, and each element $[{\bf
f}]\in\overline{\mathcal M}_{0,3}(M, A, J)_1$ has a unique
representative of the form
$${\bf f}=(\CP^1; 0, 1, \infty; f=(f_1,\cdots, f_m)),$$
where $f_k:\CP^1\to M_k$ are $J$-holomorphic maps in the homology
classes $A_k$, $k=1,\cdots,m$. Note that the other strata of
$\overline{\mathcal M}_{0,3}(M, A, J)$ have at least codimension
two. For homology classes $\alpha^{(k)}_i\in H_\ast(M_k,\Z)$,
$i=1,2,3$ and $k=1,\cdots, m$, satisfying
  the dimension condition
$$\deg(\alpha^{(k)}_1)+\deg(\alpha^{(k)}_2)+\deg(\alpha^{(k)}_3)=\dim
M_k+ 2c_1(A_k),$$
   we may choose the pseudo-cycle representatives
$\overline\alpha_i^{(k)}: U_i^{(k)}\to M$, $i=1,2,3$ and
$k=1,\cdots,m$ such that:
  \begin{description}
\item[(i)] $(\times^m_{k=1}\overline\alpha_1^{(k)})
\times(\times^m_{k=1}\overline\alpha_2^{(k)})\times(\times^m_{k=1}\overline\alpha_3^{(k)})$
is transverse to the evaluations ${\rm EV}^{J,A}_{0,3}$ in
(\ref{e:7.12}) and $E_{A,J,{\bf z}}$ in (\ref{e:7.13}),

\item[(ii)] each
$\overline\alpha_1^{(k)}\times\overline\alpha_2^{(k)}\times\overline\alpha_3^{(k)}$
is transverse to the evaluations $E_{A_k,J_k,{\bf z}}$ and
$${\rm EV}^{J_k,A_k}_{0,3}:\overline{\mathcal M}_{0,3}(M_k, A_k, J_k)\to
M_k^3,\;[{\bf f}_k]\mapsto (f_k(0), f_k(1), f_k(\infty))
$$
for $k=1,\cdots,m$.
\end{description}
Then as above we get that the Gromov-Witten invariant
\begin{eqnarray}\label{e:7.16}
&&\Psi^M_{A,0,3}(pt;\times^m_{k=1}\alpha^{(k)}_1,\times^m_{k=1}\alpha^{(k)}_2,\times^m_{k=1}\alpha^{(k)}_3)\\
&&=({\rm
EV}^{J,A}_{0,3})\cdot\bigl((\times^m_{k=1}\overline\alpha^{(k)}_1)
\times(\times^m_{k=1}\overline\alpha^{(k)}_2)\times(\times^m_{k=1}\overline\alpha^{(k)}_3)\bigr)\nonumber\\
&&=({\rm EV}^{J,A}_{0,3})|_{\overline{\mathcal M}_{0,3}(M, A, J)_1}
\cdot\bigl((\times^m_{k=1}\overline\alpha^{(k)}_1)
\times(\times^m_{k=1}\overline\alpha^{(k)}_2)\times(\times^m_{k=1}\overline\alpha^{(k)}_3)\bigr)\nonumber\\
&&=E_{A,J,{\bf z}}\cdot\bigl((\times^m_{k=1}\overline\alpha^{(k)}_1)
\times(\times^m_{k=1}\overline\alpha^{(k)}_2)\times(\times^m_{k=1}\overline\alpha^{(k)}_3)\bigr)\nonumber
\end{eqnarray}
because of (\ref{e:7.14}). Note that ${\mathcal
M}(M,A,J)=\prod^m_{k=1}{\mathcal M}(M_k, A_k, J_k)$.
  It easily follows from the above (i) and (ii)  that
\begin{eqnarray*}
&&E_{A,J,{\bf z}}\cdot\bigl((\times^m_{k=1}\overline\alpha^{(k)}_1)
\times(\times^m_{k=1}\overline\alpha^{(k)}_2)\times(\times^m_{k=1}\overline\alpha^{(k)}_3)\bigr)\\
&&=\prod^m_{k=1}E_{A_k, J_k, {\bf z}}\cdot(\overline\alpha^{(k)}_1
\times\overline\alpha^{(k)}_2\times\overline\alpha^{(k)}_3)\\
&&=\prod^m_{k=1}\Psi^{M_k}_{A_k, 3}(\alpha^{(k)}_1,
\alpha^{(k)}_2,\alpha^{(k)}_3)\\
&&=\prod^m_{k=1}\Psi^{M_k}_{A_k,0,3}(pt;\alpha^{(k)}_1,
\alpha^{(k)}_2,\alpha^{(k)}_3).
\end{eqnarray*}
The final step comes from Proposition~\ref{prop:7.6}.
  This and
(\ref{e:7.16}) lead to (\ref{e:7.15}). \hfill$\Box$
\end{document}